\documentclass[12pt,a4paper]{article}

\usepackage{amsbsy,amsfonts,amsmath,amssymb,amsthm}
\usepackage{array}
\usepackage{bm}
\usepackage{color}
\usepackage{enumerate}
\usepackage{mathrsfs}
\usepackage{latexsym}
\usepackage[hidelinks]{hyperref}
\definecolor{wineRed}{rgb}{0.7,0,0.3}
\definecolor{grandBleu}{rgb}{0,0,0.8}
\definecolor{darkGreen}{rgb}{0,0.4,0}
\definecolor{blueViolet}{rgb}{0.4,0,1.0}
\definecolor{bloodOrange}{rgb}{0.85,0.05,0}
\definecolor{mycolor}{rgb}{0.8,0,0.2}
\definecolor{}{rgb}{0.8,0,0.2}
\usepackage[textsize=small]{todonotes}
\usepackage{cite}

\usepackage{comment}
 
\DeclareMathAlphabet{\mathpzc}{OT1}{pzc}{m}{it}

\usepackage[textsize=small]{todonotes}
\numberwithin{equation}{section}

\theoremstyle{plain}
\newtheorem{mTh}{Main Theorem} 
\newtheorem{lem}{Lemma}%[section]

\newtheorem{prop}{Proposition}%[section]
\theoremstyle{definition}
\newtheorem{defn}{Definition}%[section]
\newtheorem{rem}{Remark}

\newtheorem{ex}{Example}

\def\N{\mathbb{N}}
\def\R{\mathbb{R}}

\def\ds{\displaystyle}

\def\Sgn{\mathop{\mathrm{Sgn}}\nolimits}

\begin{document}
\thispagestyle{plain}
\begin{center}
    \textbf{\Large One-dimensional optimal control problems for time-discrete constrained quasilinear diffusion equations of Allen--Cahn types}\footnotemark[1]
\end{center}
    \bigskip
%%
%\dedication{A line for dedication}
%%
%\vspace{-0.5ex}
\begin{center}
    \textsc{Shodai Kubota}
    %\footnotemark[3]
    \\[1ex]
    {Department of Mathematics and Informatics, \\ Graduate School of Science and Engineering, Chiba University, \\ 1-33, Yayoi-cho, Inage-ku, 263-8522, Chiba, Japan}
    \\[0ex]
    ({\ttfamily skubota@chiba-u.jp})
\end{center}
%\vspace{-1ex}

\footnotetext[1]{
AMS Subject Classification: 
35K59,  % Quasilinear parabolic equations
35K67,  % Singular parabolic equations 
35K87,  % Systems of parabolic variational inequalities
49J20,  % Optimal control problems involving partial differential equations
\\[1ex]
Keywords: optimal control problems, time-discrete, Allen--Cahn type, quasilinear diffusion equations, regularized approximating problems, necessary optimality conditions; limiting optimality condition
}
\bigskip

\noindent
{\bf Abstract.}
In this paper, we consider a class of optimal control problems for a one-dimensional time-discrete constrained quasilinear diffusion state-systems of singular Allen--Cahn types and its regularized approximating problems. 
We note that the control parameter for each system is given by physical temperature. 
The principal part of this paper is started with the verification of a Key-Theorem dealing with the decompositions of the subdifferentials of the governing convex energies of the state-systems. 
On this basis, we will prove five Main Theorems, concerned with: the solvability and precise regularity results of state-systems; the continuous-dependence of the solutions to state-systems including convergences in spatially $C^1$-topologies; the existence and parameter-dependence of optimal controls; the necessary optimality conditions for approximate optimal controls; precise characterizations of the approximating limit of the optimality conditions.  
\newpage

\section*{Introduction}
Let $0 < L < \infty$ and $0 < T < \infty$ be fixed positive constants, let $\Omega := (-L, L) \subset \mathbb{R}$ be a one-dimensional spatial domain. 
Besides, we set $X := L^2(\Omega)$ and $ Y := H^1(\Omega) $.%$\mathbb{X} := [X]^n$ with norm $|\cdot|_{\mathbb{X}}$. 

The aim of this work is to give some advanced observations for the optimal control problem, governed by the following Allen--Cahn type equation with singularities: 
%In this paper, we consider a class of optimal control problems for the following time-discrete Allen--Cahn type equation, denoted by \hyperlink{(AC)$^{(0, 0)}$}{(AC)$^{(0, 0)}$}: 
%
%\noindent
%\hypertarget{(AC)$^{(0, 0)}$}{(AC)$^{(0, 0)}$}
%\vspace{-1ex}
%\begin{equation}\label{0}
%\left\{ \parbox{14.0cm}{
%    $\ds \frac{1}{\tau} (w_i - w_{i-1}) -\partial_x \left( \frac{\partial_x w_i}{|\partial_x w_i|}  + \nu^2 \partial_x w_i \right) +\partial I_{[-1, 1]}(w_i) +g(w_i) \ni M_u u_i $ \quad in $ \Omega $,
%\\[1ex]
%$\partial_x w_i(\pm L) = 0 $, \quad $ i = 1, 2, 3, \ldots, n $, 
%\\[1ex]
%$ w_0 \in H^1(\Omega) $.
%}\right. 
%\end{equation}
\begin{align}\label{0}
    & \begin{cases}
        \ds \partial_t w -\partial_x  \left( \frac{\partial_x w}{|\partial_x w|}  + \nu^2 \partial_x w \right) +\partial I_{[-1, 1]}(w) +g(w) \ni M_u u, ~~\mbox{in } (0, T) \times \Omega,
        \\[1ex]
        \partial_x w(t, \pm L) = 0, ~ t \in (0, T), 
        \\[0ex]
        w(0, x) = w_0(x), ~ x \in \Omega.
    \end{cases}
\end{align}
The equation \eqref{0} is based on the modeling method of Visintin \cite{MR1423808} as a possible mathematical model of solid-liquid phase transitions in a mesoscopic length scale. 
In this context, the unknown $ w = w(t, x) $ is the nonconserved order parameter that indicates the physical phase of material: $w = 1 $ (resp. $w = -1$) corresponds to pure liquid (resp. solid). 
%$ w = [w_1, \ldots , w_n] $ is the nonconserved order parameter that indicates the physical phase of material: $w_i(x) = 1 $ (resp. $w_i(x) = -1$) corresponds to pure liquid (resp. solid), for $x \in \Omega$, $i = 1, 2, 3, \ldots , n$. 
Besides, $w_0 \in Y$ is a given initial data of $w$. 
%The forcing term $u = [u_1, \ldots , u_n] \in \mathbb{X}$ 
The forcing term $ u = u(t, x) $ denotes the control variable that controls the profile of the solution $w$, and physically, $u$ is the relative temperature. 
Furthermore, $g$ is a semi-monotone $C^1$-function on $\mathbb{R}$, and $\nu > 0$ and $M_u \geq 0$ are fixed constants. 
In addition, $\partial I_{[-1, 1]}$ is the subdifferential of the indicator function $ I_{[-1,1]} $ on the closed interval $[-1,1]$, that is
defined as:
\begin{equation*} %\label{indicator}
I_{[-1,1]} (z) := \left\{
\begin{array}{cl}
0,    & \mbox{ if } z \in [-1,1], \\[1mm]
+ \infty,  & \mbox{ otherwise. }
\end{array}
\right.
\end{equation*}

Recently, Allen--Cahn type equations kindred to \eqref{0} have been studied by a lot of mathematicians (cf. \cite{MR2836557,MR2459669,MR2509574,MR2101878,MR3661429}), and in many previous works, the singular term $ -\partial_x \bigl( \frac{\partial_x w}{|\partial_x w|} \bigr) $ is approximated by quasilinear diffusions:
\begin{align}\label{sing.eps}
    & -\partial_x \left( \frac{\partial_x w}{\sqrt{\varepsilon^2 +|\partial_x w|^2}} +\nu^2 \partial_x w \right) ~\mbox{with a constant}~ \varepsilon > 0,
\end{align}
and the singular term $ \partial I_{[-1, 1]}(w) $ of the set-valued subdifferential is approximated by means of the \emph{Yosida regularizations,} or the following $ C^1 $-regularizations:
\begin{align}\label{Ind.dlt}
    & \partial \widetilde{I}_{[-1, 1]}^\delta(r) := \begin{cases}
        0, ~\mbox{if $ -1 \leq r \leq 1 $,}
        \\[1ex]
        \ds \frac{1}{2 \delta^2} \frac{r}{|r|}(r -1)^2, ~\mbox{if $ 1 < |r| \leq 1 +\delta $,}
        \\[2ex]
        \ds \frac{1}{\delta} \frac{r^2}{|r|} -\left( \frac{1}{\delta} +\frac{1}{2} \right), ~\mbox{if $ |r| > 1 +\delta $,}
    \end{cases}
    \mbox{for $ r \in \R $, and $ \delta > 0 $.}
\end{align}
The setting \{\eqref{sing.eps},\eqref{Ind.dlt}\} was adopted as a representative effective range of approximating method for the singular terms. Indeed, the previous researches of \cite{MR2836557,MR2459669,MR2509574} adopted the setting \{\eqref{sing.eps},\eqref{Ind.dlt}\} to deal with the optimal control problems governed by \eqref{0}, and as a consequence, they obtained a distributional characterizations of the first necessary optimality condition for the optimal control. However, in the previous results, we still have two unfinished issues. 

The first unfinished issue is to make clear the effective range of approximating methods, besides \{\eqref{sing.eps},\eqref{Ind.dlt}\}. From the viewpoint of application, it is beneficial if we could choose a matching approximating method for specific property of each scientific/technological background. In this light, it could be said that the current mathematical results would not be sufficient to respond to the requirements for such flexibility. 

In the meantime, the second unfinished one is concerned with the precise observation for the distribution, which appears in the optimality condition, obtained in \cite{MR2836557,MR2459669,MR2509574} ,and is roughly expressed as:
\begin{align}\label{illposed}
    & \zeta^\circ \sim \partial_x \bigl[ \mbox{\large$\mathfrak{d}$}(\partial_x w^\circ) \partial_x u^\circ \bigr] -\bigl[ \partial I_{[-1, 1]} \bigr]\mbox{\bf\large'}(w^\circ)u^\circ ~\mbox{in $ \mathfrak{D}'\bigl( (0, T) \times \Omega \bigr) $,}
\end{align}
with use of Dirac's delta {\large$ \mathfrak{d} $}, some operator $ \bigl[ \partial I_{[-1, 1]} \bigr]\mbox{\bf\large'} $ corresponding to derivative of the subdifferential $ \partial I_{[-1, 1]} $, an optimal control $ u^\circ $, and a solution $ w^\circ $ to the singular Allen--Cahn type equation \eqref{0} when $ u = u^\circ $. The expression \eqref{illposed} looks ill-posed as a variational formula. Nevertheless, it would be expectable that the support $ \mathrm{spt}(\zeta^\circ) $ of the distribution $ \zeta^\circ $ would be somehow associated with $ \{ \partial_x w^\circ = 0 \} \cup \{ |w^\circ| = 1 \}$ $ ( \approx \{ \partial_x w^\circ = 0 \} )$, and $ \zeta^\circ $ would have some functional expression on the outside of the support:
\begin{center}
    $ \mathrm{spt} (\zeta^\circ)^\mathrm{C} \approx \{ \partial_x w^\circ \ne 0 \} $ $ \bigl( =  \{ \partial_x w^\circ \ne 0 \} \cap \{ |w^\circ| < 1 \} \bigr) $. 
\end{center}

For developments of the unfinished issues, we now let $ n \in \N $ be a fixed number, define $ \mathbb{X} := [X]^n $, and consider the time-discrete scheme for the singular Allen--Cahn type equation \eqref{0}, with the time step size $ \tau := T/n $, as a simplified state-equation of our optimal control problem. On this basis, we set our goal to study a class of optimal control problems, denoted by \hyperlink{(OP)$^{(\varepsilon, \delta)}$}{(OP)$^{(\varepsilon, \delta)}$}, with constant parameters $ \varepsilon, \delta \in [0, 1] $. 
\begin{description}
    \item[\textmd{\hypertarget{(OP)$^{(\varepsilon, \delta)}$}{(OP)$^{(\varepsilon, \delta)}$:}}]to find a sequence of functions $ u^* = [u_1^*, \ldots , u_n^*] \in \mathbb{X} $, called \emph{optimal control}, which minimizes the following cost functional $ \mathcal{J}^{(\varepsilon, \delta)}  = \mathcal{J}^{(\varepsilon, \delta)} (u) $ on $ \mathbb{X} $, defined as
        \begin{align}\label{J}
            \mathcal{J}^{(\varepsilon, \delta)} : u \in \mathbb{X} \mapsto \mathcal{J}^{(\varepsilon, \delta)}(u) := \frac{M_w }{2} |(w -w^\mathrm{ad})|_{\mathbb{X}}^2 +\frac{M_u }{2} |u|_{\mathbb{X}}^2  \in [0, \infty), 
        \end{align}
        i.e.
        \begin{align*}
            \ds \mathcal{J}^{(\varepsilon, \delta)}(u^*) = 
            %\min_{[u, v] \in \mathscr{U}_\mathrm{ad}^K} \mathcal{J}_{\varepsilon}(u, v),
            \min \left\{ \begin{array}{l|l} 
               \mathcal{J}^{(\varepsilon, \delta)}(u) & u \in \mathbb{X}
        \end{array} \right\},
    \end{align*}
    %where $ \mathcal{J}^{(\varepsilon, \delta)} = \mathcal{J}^{(\varepsilon, \delta)}(u) $ is a cost functional on $ \mathbb{X} $, defined as follows:
        where for any $ u = [u_1, \dots, u_n] \in \mathbb{X} $, the sequence of functions  $w = [w_1, \ldots , w_n] \in \mathbb{X}$ is a solution to the following time-discrete quasilinear equations, denoted by \hyperlink{(AC)$^{(\varepsilon, \delta)}$}{(AC)$^{(\varepsilon, \delta)}$}.
\end{description}
\begin{description}
    \item[\textmd{\hypertarget{(AC)$^{(\varepsilon, \delta)}$}{(AC)$^{(\varepsilon, \delta)}$}:}]to find a sequence of functions $ w = [w_1, \cdots, w_n] \in \mathbb{X} $, which fulfills
\begin{align*}%\label{2}
    & \begin{cases}
        \ds \frac{1}{\tau} (w_i - w_{i-1}) -\partial_x \bigl( (f^\varepsilon)'(\partial_x w_i) + \nu^2 \partial_x w_i \bigr) +K^{\delta}(w_i) +g(w_i) \ni M_u u_i ~~\mbox{in $ \Omega $,}
\\[1ex]
\partial_x w_i(\pm L) = 0, ~~ i = 1, 2, 3, \ldots, n, 
    \end{cases} 
\end{align*}
starting from the initial value $ w_0 \in Y $.
\end{description}
In this context, the equation (AC)$^{(0, 0)}$, i.e. the case when $ \varepsilon = \delta = 0 $, corresponds to the time-discrete scheme of the Allen--Cahn type equation \eqref{0}, while the equations \hyperlink{(AC)$^{(\varepsilon, \delta)}$}{(AC)$^{(\varepsilon, \delta)}$}, for positive $ \varepsilon, \delta $, are its approximating equations. In this regard, the effective range of approximation, brought in the first unfinished issue, will be presented in forms of assumptions for:
\begin{itemize}
    \item[--]the regularization sequence $ \{ f^\varepsilon \}_{\varepsilon \in (0, 1]} \subset C^2(\R) $ for the function $ f^0 : r \in \R \mapsto |r| \in [0, \infty)  $ of the absolute value; 
    \item[--]the regularization sequence $ \{ K^\delta \}_{\delta \in (0, 1]} \subset C^1(\R) $ for the set-valued function $ K^0 : r \in \R \mapsto \partial I_{[-1, 1]}(r) \in 2^\R $ of the subdifferential of indicator function.
\end{itemize}
In addition, we will obtain a positive answer for the second unfinished issue, as a consequence of five  Main Theorems. Here, we give their assertions, briefly.

%For each $\varepsilon \in [0, 1]$ and $\delta \in [0, 1]$, the optimal control problem, denoted by \hyperlink{(OP)$^{(\varepsilon, \delta)}$}{(OP)$^{(\varepsilon, \delta)}$}, is prescribed as follows:
%
%In the optimal control problems \hyperlink{(OP)$^{(\varepsilon, \delta)}$}{(OP)$^{(\varepsilon, \delta)}$} for $\varepsilon \in [0, 1]$ and $\delta \in [0, 1]$, the function $w^\mathrm{ad} = [w_1^\mathrm{ad}, \ldots , w_n^\mathrm{ad}] \in \mathbb{X}$ is given \emph{admissible target profile} of the order parameter $ w = [w_1, \ldots , w_n]$, and the coefficient $M_w \geq 0$ is a fixed constant.
%\bigskip
%
%
%Now, based on these, we set the goal of this paper to prove five Main Theorems, summarized as follows:

\begin{description}
    \item[{\boldmath Main Theorem \ref{mainTh01}:}]Existence, uniqueness, and $H^2$-regularity of the solution $w = [w_1, \ldots,$ $ w_n]$ to the state-system \hyperlink{(AC)$^{(\varepsilon, \delta)}$}{(AC)$^{(\varepsilon, \delta)}$}, for every $ \varepsilon \in [0, 1]$, $\delta \in [0, 1]$, initial data $ w_0 \in Y $, and forcing term $ u = [u_1, \ldots , u_n] \in \mathbb{X} $.
    
\item[{\boldmath Main Theorem \ref{MainTh02}:}]Continuous dependence of solutions to the state-systems \hyperlink{(AC)$^{(\varepsilon, \delta)}$}{(AC)$^{(\varepsilon, \delta)}$}, with respect to 
 the constants $ \varepsilon \in [0, 1]$, $\delta \in [0, 1]$, initial data $ w_0 \in Y $, and forcing term $ u = [u_1, \ldots , u_n] \in \mathbb{X} $, including the (strong) convergence in $[C(\overline{\Omega})]^n$. 
    
\item[{\boldmath Main Theorem \ref{mainTh02}:}]Results concerned with the control problems.
    \vspace{-0.5ex}
\item[{\boldmath~~(\hyperlink{III-A}{III-A})(Solvability of optimal control problems):}]Existence for the optimal control problem \hyperlink{(OP)$^{(\varepsilon, \delta)}$}{(OP)$^{(\varepsilon, \delta)}$}, for every constants $ \varepsilon \in [0, 1] $, $\delta \in [0, 1]$, and initial data $ w_0 \in Y $ with $\hat{K}^{\delta}(w_0) \in L^1(\Omega)$.
    \vspace{-0.5ex}
\item[{\boldmath~~(\hyperlink{III-B}{III-B})(Parameter dependence of optimal controls):}]Some semi-continuous dependence of  the optimal controls, with respect to the constants $ \varepsilon \in [0, 1] $, $\delta \in [0, 1]$, and initial data $ w_0 \in Y$ with $\hat{K}^{\delta}(w_0) \in L^1(\Omega)$.

\item[{\boldmath Main Theorem \ref{mainTh03}:}] Derivation of the first order necessary optimality conditions for \linebreak \hyperlink{(OP)$^{(\varepsilon, \delta)}$}{(OP)$^{(\varepsilon, \delta)}$}, via adjoint method, under $\varepsilon \in (0, 1]$ and $\delta \in (0, 1]$. %$ \varepsilon > 0 $, $ [\eta_0, \theta_0] \in (V \cap L^\infty(\Omega)) \times V_0 $,         and $ K = \jump{\kappa^0, \kappa^1} \in \mathfrak{K}_0 $.

\item[{\boldmath Main Theorem \ref{mainTh04}:}] The optimality conditions which are obtained as approximation limits of the necessary conditions as $\varepsilon, \delta \downarrow 0$, and a precise characterization of the limiting optimality conditions.
\end{description}

Details of the Main Theorems are stated in Section 3, after the preliminaries in Section 1, and the auxiliary results in Section 2. 
In particular, the Main Theorem 1 is proved by means of the Key-Theorem which is stated in Section 3. 
The proof of Key-Theorem is given in Section 4. 
After Section 4, the remaining Sections 5--9 will be devoted to the proofs of the respective five Main Theorems 1--5.

\section{Preliminaries} 

We begin by prescribing the notations used throughout this paper. 
\medskip

\noindent
\underline{\textbf{\textit{Basis notations.}}} 
For arbitrary $ r_0 $, $ s_0 \in [-\infty, \infty]$, we define:
\begin{equation*}
r_0 \vee s_0 := \max\{r_0, s_0 \}\ \mbox{and}\ r_0 \wedge s_0 := \min\{r_0, s_0 \},
\end{equation*}
%for arbitrary $-\infty \leq r \leq s \leq \infty$, we define the truncation function (operator) $\mathcal{T}_r^s : \mathbb{R} \longrightarrow [r, s] $ by letting:
%\begin{equation*}
%\xi \in \mathbb{R} \mapsto \mathcal{T}_r^s \xi := r \vee(s \wedge \xi) \in [r, s],
%\end{equation*}
and in particular, we set:
\begin{equation*}
    [r]^+ := r \vee 0 \ \mbox{and}\ [r]^- :=  -(r \wedge 0), \mbox{ for any $ r \in \R $.}
\end{equation*}    

\noindent
\underline{\textbf{\textit{Abstract notations.}}}
For an abstract Banach space $ E $, we denote by $ |\cdot|_{E} $ the norm of $ E $, and denote by $ \langle \cdot, \cdot \rangle_E $ the duality pairing between $ E $ and its dual $ E^* $. Let $I_d : E \longrightarrow E $ be the identity map from $E$ onto $E$. In particular, when $ H $ is a Hilbert space, we denote by $ (\cdot,\cdot)_{H} $ the inner product of $ H $. 

For any subset $ A $ of a Banach space $ E $, let $ \chi_A : E \longrightarrow \{0, 1\} $ be the characteristic function of $ A $, i.e.:
    \begin{equation*}
        \chi_A: w \in E \mapsto \chi_A(w) := \begin{cases}
            1, \mbox{ if $ w \in A $,}
            \\[0.5ex]
            0, \mbox{ otherwise.}
        \end{cases}
    \end{equation*}

For two Banach spaces $ E $ and $ \Xi $,  we denote by $  \mathscr{L}(E; \Xi)$ the Banach space of bounded linear operators from $ E $ into $ \Xi $, and in particular, we let $ \mathscr{L}(E) := \mathscr{L}(E; E) $. 

For Banach spaces $ E_1, \dots, E_N $, with $ 1 < N \in \mathbb{N} $, let $ E_1 \times \dots \times E_N $ be the product Banach space endowed with the norm $ |\cdot|_{E_1 \times \cdots \times E_N} := |\cdot|_{E_1} + \cdots +|\cdot|_{E_N} $. However, when all $ H_1, \dots, H_N $ are Hilbert spaces, $ H_1 \times \dots \times H_N $ denotes the product Hilbert space endowed with the inner product $ (\cdot, \cdot)_{H_1 \times \cdots \times H_N} := (\cdot, \cdot)_{H_1} + \cdots +(\cdot, \cdot)_{H_N} $ and the norm $ |\cdot|_{H_1 \times \cdots \times H_N} := \bigl( |\cdot|_{H_1}^2 + \cdots +|\cdot|_{H_N}^2 \bigr)^{\frac{1}{2}} $. In particular, when all $ E_1, \dots,  E_N $ coincide with a Banach space $ \Xi $, we write:
\begin{equation*}
    [\Xi]^N := \overbrace{\Xi \times \cdots \times \Xi}^{\mbox{$N$ times}}.
\end{equation*}
Additionally, for any transform (operator) $ \mathcal{T} : E \longrightarrow \Xi $, we let:
\begin{equation*}
    \mathcal{T}[w_1, \dots, w_N] := \bigl[ \mathcal{T} w_1, \dots, \mathcal{T} w_N \bigl] \mbox{ in $ [\Xi]^N $, \quad for any $ [w_1, \dots, w_N] \in [E]^N $.}
\end{equation*}

\noindent
\underline{\textbf{\textit{Specific notations of this paper.}}} 
As is mentioned in the previous section, let $ (0, T) \subset \R$ be a bounded time-interval with a finite constant $ T > 0 $. Let $\Omega := (-L, L) \subset \R $ be a fixed spatial bounded domain with a finite constant $ L > 0 $. Especially, we denote by $ \partial_x $ the distributional spatial derivative. 
\medskip

On this basis, we define  
\begin{align*} 
    X := L^2(\Omega), &\,\,  X_\Gamma := \left\{ \begin{array}{l|l}
            \tilde{z} & \tilde{z} : \Gamma \longrightarrow \mathbb{R}
        \end{array} \right\} ~ (\sim \mathbb{R}^2),\, Y := H^1(\Omega) ~(\subset C(\overline{\Omega})),
        \\
        &\, \mathbb{X} := [X]^n, \mbox{ and } \mathbb{Y} := [Y]^n, \mbox{ for any } n\in \mathbb{N}.
\end{align*} 
\bigskip

\noindent
\underline{\textbf{\textit{Notations in convex analysis. (cf. \cite[Chapter II]{MR0348562})}}} 
    Let $ H $ be an abstract Hilbert space. For a proper, lower semi-continuous (l.s.c.), and convex function $ \psi : H \longrightarrow (-\infty, \infty] $, we denote by $ D(\psi) $ the effective domain of $ \psi $. Also, we denote by $\partial \psi$ the subdifferential of $\psi$. The subdifferential $ \partial \psi $ corresponds to a weak differential of convex function $ \psi $, and it is known as a maximal monotone graph in the product space $ H \times H $. The set $ D(\partial \psi) := \bigl\{ z \in H \ |\ \partial \psi(z) \neq \emptyset \bigr\} $ is called the domain of $ \partial \psi $. We often use the notation ``$ [z_{0}, z_{0}^{*}] \in \partial \psi $ in $ H \times H $\,'', to mean that ``$ z_{0}^{*} \in \partial \psi(z_{0})$ in $ H $ for $ z_{0} \in D(\partial\psi) $'', by identifying the operator $ \partial \psi $ with its graph in $ H \times H $.
\medskip

\begin{ex}[Examples of the subdifferential]\label{exConvex}
    As one of the representatives of the subdifferentials, we exemplify the following set-valued signal function $ \Sgn^N: \R^N \longrightarrow 2^{\mathbb{R}^N} $, with $ N \in \mathbb{N} $, which is defined as:
\begin{align}\label{Sgn^d}
    \xi = [\xi_1, & \dots, \xi_N] \in \mathbb{R}^N \mapsto \Sgn^N(\xi) = \Sgn^N(\xi_1, \dots, \xi_N) 
    \nonumber
    \\
    & := \left\{ \begin{array}{ll}
            \multicolumn{2}{l}{
                    \ds \frac{\xi}{|\xi|} = \frac{[\xi_1, \dots, \xi_N]}{\sqrt{\xi_1^2 +\cdots +\xi_N^2}}, ~ } \mbox{if $ \xi \ne 0 $,}
                    \\[3ex]
            \mathbb{D}^N, & \mbox{otherwise,}
        \end{array} \right.
    \end{align}
where $ \mathbb{D}^N $ denotes the closed unit ball in $ \mathbb{R}^N $ centered at the origin. Indeed, the set-valued function $ \Sgn^N $ coincides with the subdifferential of the Euclidean norm $ |{}\cdot{}| : \xi \in \mathbb{R}^N \mapsto |\xi| = \sqrt{\xi_1^2 + \cdots +\xi_N^2} \in [0, \infty) $, i.e.:
\begin{equation*}
\partial |{}\cdot{}|(\xi) = \Sgn^N(\xi), \mbox{ for any $ \xi \in D(\partial |{}\cdot{}|) = \mathbb{R}^N $,}
\end{equation*}
and furthermore, it is observed that:
\begin{equation*}
    \partial  |{}\cdot{}|(0) = \mathbb{D}^N \begin{array}{c} \subseteq_{\hspace{-1.25ex}\mbox{\tiny$_/$}}  
    \end{array} [-1, 1]^N 
        = \bigl[ \partial_{\xi_1}  |{}\cdot{}| \times \cdots \times \partial_{\xi_N}  |{}\cdot{}| \bigr](0).
\end{equation*}
\end{ex}
\medskip

Finally, we mention about a notion of functional convergence, known as ``Mosco-convergence''. 
 
\begin{defn}[Mosco-convergence: cf. \cite{MR0298508}]\label{Def.Mosco}
    Let $ H $ be an abstract Hilbert space. Let $ \psi : H \longrightarrow (-\infty, \infty] $ be a proper, l.s.c., and convex function, and let $ \{ \psi_n \}_{n \in \mathbb{N}} $ be a sequence of proper, l.s.c., and convex functions $ \psi_n : H \longrightarrow (-\infty, \infty] $, $ n = 1, 2, 3, \dots $.  Then, it is said that $ \psi_n \to \psi $ on $ H $, in the sense of Mosco, as $ n \to \infty $, iff. the following two conditions are fulfilled:
\begin{description}
    \item[(\hypertarget{M_lb}{M1}) The condition of lower-bound:]$ \ds \varliminf_{n \to \infty} \psi_n(\check{w}_n) \geq \psi(\check{w}) $, if $ \check{w} \in H $, $ \{ \check{w}_n  \}_{n \in \mathbb{N}} \subset H $, and $ \check{w}_n \to \check{w} $ weakly in $ H $, as $ n \to \infty $; 
    \item[(\hypertarget{M_opt}{M2}) The condition of optimality:]for any $ \hat{w} \in D(\psi) $, there exists a sequence \linebreak $ \{ \hat{w}_n \}_{n \in \mathbb{N}}  \subset H $ such that $ \hat{w}_n \to \hat{w} $ in $ H $ and $ \psi_n(\hat{w}_n) \to \psi(\hat{w}) $, as $ n \to \infty $.
\end{description}
    As well as, if the sequence of convex functions $ \{ \widehat{\psi}_\varepsilon \}_{\varepsilon \in \Lambda} $ is labeled by a continuous argument $\varepsilon \in \Lambda$ with a range $\Lambda \subset \mathbb{R}$ , then for any $\varepsilon_{0} \in \Lambda$, the Mosco-convergence of $\{ \widehat{\psi}_\varepsilon \}_{\varepsilon \in \Lambda}$, as $\varepsilon \to \varepsilon_{0}$, is defined by those of subsequences $ \{ \widehat{\psi}_{\varepsilon_n} \}_{n \in \mathbb{N}} $, for all sequences $\{ \varepsilon_n \}_{n \in \mathbb{N}} \subset \Lambda$, satisfying $\varepsilon_{n} \to \varepsilon_{0}$ as $n \to \infty$.
\end{defn}

\begin{rem}\label{Rem.MG}
    Let $ H $, $ \psi $, and $ \{ \psi_n \}_{n \in \mathbb{N}} $ be as in Definition~\ref{Def.Mosco}. Then, the following hold.
\begin{description}
    \item[(\hypertarget{Fact1}{Fact\,1})](cf. \cite[Theorem 3.66]{MR0773850} and \cite[Chapter 2]{Kenmochi81}) Let us assume that
    \begin{equation}\label{Mosco01}
    \psi_n \to \psi \mbox{ on $ H $, in the sense of  Mosco, as $ n \to \infty $,}
    \vspace{-1ex}
\end{equation}
and
\begin{equation*}
\left\{ ~ \parbox{10cm}{
$ [w, w^*] \in H \times H $, ~ $ [w_n, w_n^*] \in \partial \psi_n $ in $ H \times H $, $ n \in \N $,
\\[1ex]
$ w_n \to w $ in $ H $ and $ w_n^* \to w^* $ weakly in $ H $, as $ n \to \infty $.
} \right.
\end{equation*}
Then, it holds that:
\begin{equation*}
[w, w^*] \in \partial \psi \mbox{ in $ H \times H $, and } \psi_n(w_n) \to \psi(w) \mbox{, as $ n \to \infty $.}
\end{equation*}
    \item[(\hypertarget{Fact2}{Fact\,2})](cf. \cite[Lemma 4.1]{MR3661429} and \cite[Appendix]{MR2096945}) Let $ N \in \mathbb{N} $ denote the dimension constant, and let $  S \subset \R^N $ be a bounded open set. Then, under the Mosco convergence as in \eqref{Mosco01}, a sequence $ \{ \widehat{\psi}_n^S \}_{n \in \mathbb{N}}$ of proper, l.s.c., and convex functions on $ L^2(S; H) $, defined as:
        \begin{equation*}
            w \in L^2(S; H) \mapsto \widehat{\psi}_n^S(w) := \left\{ \begin{array}{ll}
                    \multicolumn{2}{l}{\ds \int_S \psi_n(w(t)) \, dt,}
                    \\[1ex]
                    & \mbox{ if $ \psi_n(w) \in L^1(S) $,}
                    \\[2.5ex]
                    \infty, & \mbox{ otherwise,}
                \end{array} \right. \mbox{for $ n = 1, 2, 3, \dots $;}
        \end{equation*}
        converges to a proper, l.s.c., and convex function $ \widehat{\psi}^S $ on $ L^2(S; H) $, defined as:
        \begin{equation*}
            z \in L^2(S; H) \mapsto \widehat{\psi}^S(z) := \left\{ \begin{array}{ll}
                    \multicolumn{2}{l}{\ds \int_S \psi(z(t)) \, dt, \mbox{ if $ \psi(z) \in L^1(S) $,}}
                    \\[2ex]
                    \infty, & \mbox{ otherwise;}
                \end{array} \right. 
        \end{equation*}
        on $ L^2(S; H) $, in the sense of Mosco, as $ n \to \infty $. 
\end{description}
\end{rem}
\medskip

\section{Auxiliary results}

In this Section, we recall the previous work \cite{MR4279540}, and show some auxiliary results. 
Let $0 \leq \alpha \in Y$ and $0 < \beta \in Y$ be fixed functions, and let us set the following convex function on $X$: 
\begin{align}\label{Phi}
    z \in X \mapsto \Phi_{\alpha, \beta}(z) := V_\alpha(z) + W_\beta(z);
\end{align}
which is defined as a sum of two convex functions on $X$, defined as follows:
%\vspace{-1ex}
    \begin{align}\label{V}
    \displaystyle z \in X \mapsto V_\alpha(z) := 
\sup \left\{ \begin{array}{l|l}
        \displaystyle \int_\Omega z \partial_x \varphi \, dx, & \parbox{4cm}{$ \varphi \in Y \cap C_\mathrm{c}(\Omega) $, such that $ |\varphi| \leq \alpha $ on $ \overline{\Omega} $}
    \end{array} \right\},
    \end{align}
and
%\vspace{-1ex}
\begin{align}\label{W}
    \displaystyle z \in X \mapsto W_\beta(z) := \left\{
\begin{array}{l}
\displaystyle \frac{1}{2}\int_{\Omega} \beta |\partial_x z|^2 dx, \, \mbox{ if } z \in Y, \\[1.5ex]
\infty, \mbox{ otherwise.}
\end{array}
\right.
\end{align}
%with fixed nonhomogeneous coefficients $ 0 \leq \alpha \in Y $ and $ 0 < \beta \in Y $. 
    The functional $ V_\alpha $, defined in \eqref{V}, is a kind of generalized total variation, so that the functional $ \Phi_{\alpha, \beta} $, defined in \eqref{Phi}, can be called a \emph{regularized total variation} with nonhomogeneous coefficients $ \alpha $ and $ \beta $. 
    
\begin{rem}(cf. \cite{MR1857292,MR1736243})\label{Rem.V_alp}
The functional $ V_\alpha $ coincides with the so-called \emph{lower semi-continuous envelope} of the following convex function:
\begin{equation*}
    z \in W^{1, 1}(\Omega) \mapsto \widetilde{V}_{\alpha}(z) := \int_\Omega \alpha |\partial_x z| \, dx \in [0, \infty), 
\end{equation*}
more precisely,
%
%Also, this functional has another expression, as \emph{weighted total variation,} which is based on the general theory, developed in \cite{MR1259102,MR1857292,MR1736243}:
\begin{align*}%\label{V_another}
    & V_\alpha(z) = 
\inf \left\{ \begin{array}{l|l}
    \displaystyle \varliminf_{i \to \infty} \widetilde{V}_\alpha(\tilde{z}_i) & \parbox{8.5cm}{$ \{ \tilde{z}_i \}_{i \in \mathbb{N}} \subset W^{1, 1}(\Omega) $, and $ \tilde{z}_i \to z $ in $ X $, as $ i \to \infty $}
\end{array} \right\}. 
    %\sup \left\{ \begin{array}{l|l}
    %    \displaystyle \int_\Omega \theta \partial_x \varphi \, dx & \parbox{4cm}{$ \varphi \in Y \cap C_\mathrm{c}(\Omega) $, such that $ |\varphi| \leq \alpha $ on $ \overline{\Omega} $} \mbox{ for any $ \theta \in X $.}
\end{align*}
\end{rem}
\medskip 

\begin{rem}\label{Rem.W_bt}
    The functional $ W_\beta $, defined in \eqref{W}, is a proper, l.s.c., and convex function on $ X $, such that $ D(W_\beta) = Y $. Moreover, the subdifferential $ \partial W_\beta \subset X \times X $ is a single valued operator, such that
    \begin{equation*}
        [z, z^*] \in \partial W_\beta \mbox{ in $ X \times X $, iff. $ \beta \partial_x z \in H_0^1(\Omega) $, and $ z^* = -\partial_x (\beta \partial_x z) $ in $ X $.}
    \end{equation*}
\end{rem}
\medskip
    
Now, we refer to the previous work \cite{MR4279540}, to recall the key-properties of $\Phi_{\alpha, \beta}$, in forms of Proposition. 

\begin{prop}\label{decom}\cite[Main Theorem]{MR4279540}
 The subdifferential $ \partial \Phi_{\alpha, \beta} \subset X \times X  $ of the convex function $ \Phi_{\alpha, \beta}$ is decomposed as follows:
\begin{align*}%\label{target}
    \partial \Phi_{\alpha, \beta} = \partial V_\alpha + \partial W_\beta \mbox{ in } X \times X,
\end{align*}
i.e. $\partial \Phi_{\alpha, \beta}$ is represented as the sum of the subdifferentials $ \partial V_\alpha \subset X \times X $ and $ \partial W_\beta \subset X \times X $ of the respective convex functions $ V_\alpha $ and $ W_\beta $. 
\end{prop}

\begin{comment}
\begin{prop}(cf. \cite{MR4279540})
 Let us define a set-valued map $\mathcal{A}^0 :D(\mathcal{A}^0) \subset X \longrightarrow 2^X$, by putting:
\begin{equation*}%\label{A.0}
 D(\mathcal{A}^0) := \left\{\begin{array}{l|ll} z \in Y & \begin{array}{lll}\multicolumn{2}{l}{\mbox{ there exists } \varpi^*\in L^\infty(\Omega) \mbox{ such that } } 
\\[0.25ex]
& \bullet ~ \varpi^* \in \Sgn(\partial_x z) \mbox{ a.e. in } \Omega
\\[0.25ex]
& \bullet ~ \alpha \varpi^* + \beta \partial_x z \in H_0^1(\Omega)  \end{array}
\end{array}\right\},
\end{equation*}
and
\begin{align*}%\label{A.0f}
 z &\, \in D(\mathcal{A}^0) \subset X \nonumber
 \\
   &\, \mapsto \mathcal{A}^0 z := \left\{\begin{array}{l|l} z^* \in X & \begin{array}{lll}\multicolumn{2}{l}{ z^* = -\partial_x \bigl( \alpha\varpi^* + \beta \partial_x z\bigr) \mbox{ in } X, }
   \\[0.25ex]
   & \mbox{for some } \varpi^*\in L^\infty(\Omega), \mbox{ satisfying } 
   \\[0.25ex]
   & \varpi^*\in \Sgn(\partial_x z) \mbox{ a.e. in } \Omega \end{array}
 \end{array} \right\}.
\end{align*}
Then, the following equation holds that: 
\begin{align*}%\label{A.sub}
\mathcal{A}^0 = \partial \Phi_{\alpha, \beta} \mbox{ in } X \times X.
\end{align*}
\end{prop}
\end{comment}

\section{Main Theorems}

We begin by setting up the assumptions needed in our Main Theorems. 
All Main Theorems are discussed under the following assumptions.  
\begin{description}
\item[\textmd{(\hypertarget{A1l}{A1})}]
    Let $\nu>0$ be a fixed constant. 
\item[\textmd{(\hypertarget{A2l}{A2})}]
    We denote by $f^0$ the absolute value function on ${\mathbb R}$, i.e., $f^0(r):=|r|$ for all $r \in {\mathbb R}$. In addition, let $\{ f^\varepsilon \}_{\varepsilon \in (0,1]} \subset C^2({\mathbb R}) $ be a sequence of convex $C^2$-regularizations of $f^0(\cdot ):=|\cdot | $, such that:
\[
\mbox{$f^\varepsilon(0) = 0$ and $f^\varepsilon (r)\geq 0$ \ for any $r\in {\mathbb R}$ and any $\varepsilon \in (0,1]$,}
\]
\begin{equation*}
\left\{
\begin{array}{l}
f^\varepsilon(r) \to f^{\varepsilon_0}(r) \  \mbox{ for any } r \in {\mathbb R}, \\[1mm]
f^\varepsilon (\cdot ) \to f^{\varepsilon_0}(\cdot ) \ \mbox{  on } {\mathbb R}, \mbox{ \  in the sense of Mosco},
\end{array}
\mbox{ as $\varepsilon \to \varepsilon_0$, for any $\varepsilon_0 \in [0, 1],$} \right.
\end{equation*}
and there exists a positive constant $C_0 > 0$, independent of $\varepsilon \in (0,1]$, satisfying:
\[
|(f^\varepsilon)'(r)| \le C_0(1 + |r|) \mbox{ \ for any $r\in {\mathbb R}$ and any $\varepsilon \in (0,1]$,}
\]
and 
\begin{align}\label{A2.loc}
(f^\varepsilon)'' \to 0 \mbox{ uniformly on } \{ |r| \geq \lambda \}, \mbox{ for any } \lambda > 0, \mbox{ as } \varepsilon \to 0.
\end{align}
    \item[\textmd{(\hypertarget{A3l}{A3})}]
        Let $g : \mathbb{R} \longrightarrow  \mathbb{R}$ be a $C^{1}$-function, which is a semi-monotone function on $\mathbb{R}$, i.e. there exists a positive constant $C_g > 0$ such that: 
        \begin{align*}%\label{g.semi}
        g(r) + C_g r \mbox{ is monotone in } r \in \mathbb{R}.
        \end{align*}         
Also, $g$ has a nonnegative primitive $ 0 \leq G \in C^{2}(\mathbb{R})$, i.e. the derivative $ G' $ coincides with $ g $ on $ \R $. 
\item[\textmd{(\hypertarget{A4l}{A4})}]
    We denote by $\hat{K}^0$, and $K^0$ the indicator function $I_{[-1, 1]}$, and the subdifferential, i.e. 
    \begin{align*}
     \hat{K}^0 = I_{[-1, 1]} \mbox{ on } \mathbb{R}, \mbox{ and } K^0 := \partial I_{[-1, 1]} \mbox{ in } \mathbb{R} \times \mathbb{R},
    \end{align*}
respectively. 
In addition, let $\{\hat{K}^\delta \}_{\delta \in (0, 1]} \subset C^2(\mathbb{R})$, and $\{K^\delta \}_{\delta \in (0, 1]} \subset C^1(\mathbb{R})$ be, respectively, the approximating sequence of $\hat{K}^0$, and $K^0$, such that: 
        \begin{itemize}
            \item $K^\delta = (\hat{K}^\delta)'$ on $\mathbb{R}$, and $\hat{K}^\delta \equiv 0$ on $[-1, 1]$, for all $\delta \in (0, 1]$;
            \item for any $\delta \in (0, 1]$, there exists a constant $C_K^\delta > 0$ satisfying $0 \leq (K^\delta)' \leq C_K^\delta$ on $\mathbb{R}$;
            \item for any $\delta_0 \in [0, 1]$, $\hat{K}^\delta \to \hat{K}^{\delta_0}$ on $\mathbb{R}$, in the sense of Mosco, as $\delta \to \delta_0$, and in particular, if $\delta_0 \in (0, 1]$, then $\hat{K}^\delta \to \hat{K}^{\delta_0}$ in $C_{\mathrm{loc}}(\mathbb{R})$, as $\delta \to \delta_0$.
        \end{itemize}      
\item[\textmd{(\hypertarget{A5l}{A5})}]
    Let $T > 0$ be a fixed constant, and let $\tau^*$ be a small positive constant, such that: 
    \begin{align*}%\label{tau*}
         \tau^* := \frac{1}{8 \bigl( C_g + 1 \bigr)}.
    \end{align*}     
    On this basis, we fix constants of the time-step number $n \in \mathbb{N}$ and the time-step size $\tau > 0$, to satisfy that: 
    \begin{align*}%\label{tau}
     0 < \tau :=  \frac{T}{n} < \tau^* .
    \end{align*}
\item[\textmd{(\hypertarget{A6l}{A6})}] 
Let $ w^\mathrm{ad} = [w_1^\mathrm{ad}, \ldots , w_n^\mathrm{ad} ] \in \mathbb{X} $ be a fixed \emph{admissible target profile}.
\end{description}

\begin{rem}\label{rem.main00}(cf. \cite{MR4228004, MR3661429, MR3670006})
 The assumption similar to (\hyperlink{A2l}{A2}) was introduced in \cite[Section 3]{MR3661429}, \cite[Remark 3.1]{MR4228004}, and \cite[Definition 3.1]{MR3670006}. 
 In this context, the typical examples of $f^\varepsilon$ are as follows:
 \begin{itemize}
\item (Hyperbola type) \ $f^\varepsilon (r )=\sqrt{ r^2+\varepsilon ^2} - \varepsilon $ \ for any $r\in {\mathbb R}$ and any $\varepsilon \in (0,1]$.
\item (Hyperbolic-tangent type) \ $\displaystyle f^\varepsilon (r )=\varepsilon \log\left(\cosh\left(\frac{r}{\varepsilon} \right)\right) $ \ for any $r\in {\mathbb R}$ and any $\varepsilon \in (0,1]$.
\item (Arctangent type) \ $\displaystyle f^\varepsilon (r )=\frac{2\varepsilon}{\pi} \left[\frac{r}{\varepsilon} \tan^{-1}\left(\frac{r}{\varepsilon}\right)-\frac12\log\left(1+\left(\frac{r}{\varepsilon}\right)^2\right)\right] $ \ for any $r\in {\mathbb R}$ and any $\varepsilon \in (0,1]$.
\end{itemize}
\end{rem}
\medskip

\begin{rem}\label{f.sub}
    For any $\varepsilon \in [0, 1]$, and let us define a functional $\tilde{V}^\varepsilon : D(\tilde{V}^\varepsilon) \subset X \longrightarrow [0, \infty]$, as follows: 
    \begin{align}\label{VV^eps}
    \displaystyle z \in X \mapsto &\,  \tilde{V}^\varepsilon(z) :=  \int_\Omega f^\varepsilon(z) \, dx \in [0, \infty].
\end{align}
    Then, from (\hyperlink{A2l}{A2}) and (\hyperlink{Fact1}{Fact\,1}), it is immediately seen that:
    \begin{description}
        \item[($ \spadesuit $1)]$ \tilde{V}^\varepsilon $ is continuous and convex on $ X $, with $ D(\tilde{V}^\varepsilon) = X $;
        \item[($ \spadesuit $2)]for any $ \varepsilon_0 \in [0, 1] $, $ \tilde{V}^\varepsilon \to \tilde{V}^{\varepsilon_0}$ on $ X $, in the sense of Mosco, as $ \varepsilon \to \varepsilon_0 $.
    \end{description}
\end{rem}
\medskip

\begin{rem}\label{rem.f.lip}
By the definition of $f^\varepsilon$, the function $(1/\nu^2)(f^\varepsilon)'$ is a maximal monotone graph on $\mathbb{R} \times \mathbb{R}$. 
Hence, its resolvent $((1/\nu^2)(f^\varepsilon)' + I_d)^{-1}$ is non-expansive. 
Therefore, we can verify that $ \left((f^\varepsilon)' + \nu^2 I_d \right)^{-1}$ is a Lipschitz continuous function with a Lipschitz constant $1/\nu^2$. 

In fact, let us fix $z_i \in \mathbb{R}$, for $i = 1, 2$, and let us set $z^*_i \in \mathbb{R} $, for $i = 1, 2$, as follows: 
\begin{align*}
 z^*_i := \left((f^\varepsilon)' + \nu^2 I_d \right)^{-1}(z_i), \mbox{ for any } i = 1, 2.
\end{align*} 
Then, we compute that: 
\begin{align*}
 z_i = \left((f^\varepsilon)' + \nu^2 I_d \right)(z^*_i) = \nu^2 \left(\frac{1}{\nu^2}(f^\varepsilon)' + I_d \right)(z^*_i), 
\end{align*}
i.e.
\begin{align*}
 \frac{1}{\nu^2}z_i =  \left(\frac{1}{\nu^2}(f^\varepsilon)' + I_d \right)(z^*_i), \mbox{ for any } i = 1, 2.
\end{align*}
Thus, we obtain that: 
\begin{align*}
 z^*_i = \left(\frac{1}{\nu^2}(f^\varepsilon)' + I_d \right)^{-1}\left(\frac{1}{\nu^2} z_i \right), \mbox{ for any } i = 1, 2.
\end{align*}
Based on these estimates, one can see that: 
\begin{align*}%\label{step1-04}
 &\, \left|\left((f^\varepsilon)' + \nu^2 I_d \right)^{-1}(z_1) - \left((f^\varepsilon)' + \nu^2 I_d \right)^{-1}(z_2) \right| \nonumber
 \\
 & \quad = \left| \left(\frac{1}{\nu^2}(f^\varepsilon)' + I_d \right)^{-1}\left(\frac{1}{\nu^2}z_1 \right) - \left(\frac{1}{\nu^2}(f^\varepsilon)' + I_d \right)^{-1}\left(\frac{1}{\nu^2}z_2 \right) \right| \nonumber
 \\
 &\quad \leq \left| \frac{z_1}{\nu^2} - \frac{z_2}{\nu^2}\right| = \frac{1}{\nu^2}|z_1 - z_2|, \mbox{ for any } z_1, z_2 \in \mathbb{R}. 
\end{align*}
\end{rem}
\medskip

\begin{rem}\label{g.prime}
The assumption (\hyperlink{A3l}{A3}) leads to: 
\begin{align*}%\label{g.prime}
 g' \geq -C_g \mbox{ on } \mathbb{R}, \mbox{ i.e. } |[g']^-|_{L^\infty(\mathbb{R})} \leq C_g.
\end{align*}
\end{rem}
\medskip

\begin{rem}\label{K.rem}
The assumption (\hyperlink{A4l}{A4}) guarantees that: 
\begin{align*}
 r \cdot K^\delta(r) \geq 0, \mbox{ for any } r \in \mathbb{R},
\end{align*}
and hence, 
\begin{align*}
 (K^\delta(z), z)_X \geq 0, \mbox{ for any } z \in X.
\end{align*} 
\end{rem}
\bigskip

Next, let us set the convex function on $X$ as follows: 
\begin{align}\label{Main.Phi}
    z \in X \mapsto \Phi^0(z) := V^0(z) + V^D(z);
\end{align}
which is defined as a sum of two convex functions on $X$, defined as follows:
%\vspace{-1ex}
    \begin{align}\label{Main.V}
    \displaystyle z \in X \mapsto V^0(z) := 
\sup \left\{ \begin{array}{l|l}
        \displaystyle \int_\Omega z \partial_x \varphi \, dx, & \parbox{4cm}{$ \varphi \in Y \cap C_\mathrm{c}(\Omega) $, such that $ |\varphi| \leq 1 $ on $ \overline{\Omega} $}
    \end{array} \right\},
    \end{align}
and
%\vspace{-1ex}
\begin{align}\label{Main.Vd}
    \displaystyle z \in X \mapsto V^D(z) := \left\{
\begin{array}{l}
\displaystyle \frac{\nu^2}{2}\int_{\Omega} |\partial_x z|^2 dx, \, \mbox{ if } z \in Y, \\[1.5ex]
\infty, \mbox{ otherwise.}
\end{array}
\right.
\end{align}
Also, for any $\varepsilon \in (0, 1]$, let $\Phi^\varepsilon$ be a proper, l.s.c, and convex function 
\begin{align}\label{Phi.eps}
    \displaystyle z \in X \mapsto \Phi^\varepsilon(z) := \left\{
\begin{array}{l}
\displaystyle \int_\Omega f^\varepsilon(\partial_x z)\, dx + \frac{\nu^2}{2}\int_{\Omega} |\partial_x z|^2 dx, \, \mbox{ if } z \in Y, \\[2.0ex]
\infty, \mbox{ otherwise.}
\end{array}
\right.
\end{align}
Moreover, for any $\delta \in [0, 1]$, let us denote by $\hat{\mathcal{K}}^\delta$ a proper, l.s.c, and convex function on $X$, defined as: 
\begin{align}\label{K.del}
 \ds \hat{\mathcal{K}}^\delta : z \in X \mapsto \hat{\mathcal{K}}^\delta(z) := \int_\Omega \hat{K}^\delta(z)\, dx, 
\end{align}
and let us denote by $\mathcal{K}^\delta$ the subdifferential $\partial \hat{\mathcal{K}}^\delta$ of $\hat{\mathcal{K}^\delta}$ in $X \times X$. 
As is easily seen, 
\begin{align}\label{K.del1}
    \mathcal{K}^\delta(z) &\, = \left\{
\begin{array}{l}
\left\{z^* \in X ~\Bigl|~ z^* \in \partial I_{[-1, 1]}(z) \mbox{ a.e. in } \Omega\right\} \mbox{ in } X, \mbox{ if } \delta = 0,\\[2.5ex]
K^\delta(z), \mbox{ if } \delta \in (0, 1],
\end{array}
\right. 
\\
&\, \mbox{ for all } z \in X \mbox{ and } \delta \in [0, 1]. \nonumber
\end{align}

Based on these, let us denote by $\Psi^{(\varepsilon, \delta)}$ a function on $X$, defined as follows: 
\begin{align}\label{Psi.eps}
    \displaystyle z \in X \mapsto &\,  \Psi^{(\varepsilon, \delta)}(z) := \left\{
\begin{array}{l}
\displaystyle  \Phi^\varepsilon(z) + \hat{\mathcal{K}}^\delta(z), \, \mbox{ if } z \in Y, \\[1.5ex]
\infty, \mbox{ otherwise,}
\end{array}
\right.
\end{align}
and this definition implies that $\Psi^{(\varepsilon, \delta)}$ is a proper, l.s.c, and convex on $X$, and
\begin{align*}
    D(\Psi^{(\varepsilon, \delta)}) &\, : = D(\Phi^\varepsilon) \cap D(\hat{\mathcal{K}}^\delta)= \left\{
\begin{array}{l}
\left\{\tilde{z}\in Y ~\Bigl|~ |\tilde{z}| \leq 1 \mbox{ on } \overline{\Omega}\right\}, \mbox{ if } \delta = 0,\\[2.5ex]
Y, \mbox{ if } \delta \in (0, 1],
\end{array}
\right. 
\\
&\, \mbox{ for all } \varepsilon, \delta \in [0, 1].
\end{align*}

The principal part of this paper is the verification of the following Key-Theorem, which is concerned with the decomposition property of the subdifferential $\partial \Psi^{(\varepsilon, \delta)} \subset X \times X$ of the convex function $\Psi^{(\varepsilon, \delta)}$. 
\paragraph{\boldmath Key-Theorem.}{
Let us fix $\varepsilon \in [0, 1]$ and $\delta \in [0, 1]$. 
Then, $[z, z^*] \in \partial \Psi^{(\varepsilon, \delta)}$ in $X \times X$ if and only if there exist $\varpi^* \in Y \cap L^\infty(\Omega)$ and $\xi \in X$ such that 
\begin{itemize}
\item $z \in H^2(\Omega)$, with $(\varpi^* + \nu^2 \partial_x z)(\pm L) = 0$. 
\item $z^* = -\partial_x \varpi^* - \nu^2 \partial_x^2 z + \xi$ in $X$, where $\varpi^* \in \partial f^\varepsilon(\partial_x z)$ and $\xi \in K^\delta(z)$ a.e. in $\Omega$.
\end{itemize}

Key-Theorem plays an important role to lead the $H^2$-regularity of the solution $w = [w_1, \ldots , w_n]$ to the state-system \hyperlink{(AC)$^{(\varepsilon, \delta)}$}{(AC)$^{(\varepsilon, \delta)}$} in Main Theorem \ref{mainTh01}.
}

\begin{rem}
In the proof of Key-Theorem, we will also use the following functional: 
\begin{align}\label{V.eps}
    \ds V^\varepsilon : z \in X \mapsto V^\varepsilon(z) := \left\{
\begin{array}{l}
\displaystyle 
    %\int_{\Omega}  f^\varepsilon(\partial_x z)\, dx, \, 
    \tilde{V}^\varepsilon(\partial_x z), \mbox{ if } z \in Y, \\[1.5ex]
\infty, \mbox{ otherwise,}
\end{array}
\right. 
\mbox{ for all } \varepsilon \in (0, 1],
\end{align}
    as extra notations, where $ \tilde{V}^\varepsilon $ is the continuous and convex function on $ X $, given in \eqref{VV^eps}.  
\linebreak
    As is easily seen from (\hyperlink{A2l}{A2}), $V^\varepsilon$ is proper and convex on $X$. 
But while, the assumption (\hyperlink{A2l}{A2}) does not guarantee the lower semi-continuity of $V^\varepsilon$ on $X$. 
Nevertheless, we can obtain: 
\begin{align*}
         \partial \Phi^\varepsilon = \partial V^\varepsilon +  \partial V^D \mbox{ in } X \times X, \mbox{ for any } \varepsilon \in (0, 1],
\end{align*}
as a consequence of Key-Theorem (the Step 1 of the proof).
\end{rem}
\medskip

\begin{rem}\label{lem.op.ap}(cf. \cite[Remark 5]{MR4279540})
Let us fix any $\varepsilon \in (0, 1]$, and let us define a map $\mathcal{A}^\varepsilon :D(\mathcal{A}^\varepsilon) \subset X \longrightarrow X$, by putting: 
\begin{equation*}%\label{A.eps}
 D(\mathcal{A}^\varepsilon) := \left\{ z \in Y ~\Bigl|~ \alpha(f^\varepsilon)'(\partial_x z) + \beta \partial_x z \in H_0^1(\Omega)\right\},
\end{equation*}
and
\begin{equation*}
 z \in D(\mathcal{A}^\varepsilon) \subset X \mapsto \mathcal{A}^\varepsilon (z) := -\partial_x \bigl( \alpha(f^\varepsilon)'(\partial_x z) + \beta \partial_x z\bigr).
\end{equation*}
Then, by applying the standard variational technique, we can observe that:
\begin{align*}
\mathcal{A}^\varepsilon = \partial \Phi^\varepsilon \mbox{ in } X \times X.
\end{align*}
\end{rem}
\bigskip

Now, the Main Theorems of this paper are stated as follows.

\begin{mTh}\label{mainTh01}
    Under the assumptions (\hyperlink{A1l}{A1})--(\hyperlink{A6l}{A6}), let us fix a constant $ \varepsilon \in [0, 1] $ and $\delta \in [0, 1]$, a forcing term $ u = [u_1, \ldots , u_n] \in \mathbb{X}$, an initial data $w_0\in Y$ satisfying $\hat{K}^\delta(w_0) \in L^1(\Omega)$.
    Then, the state-system \hyperlink{(AC)$^{(\varepsilon, \delta)}$}{(AC)$^{(\varepsilon, \delta)}$} admits a unique solution $w = [w_1, \ldots , w_n] \in \mathbb{Y}$, and moreover, the following items holds:
    \begin{description}
        \item[\textmd{(\hypertarget{I}{I})}]$w_i \in H^2(\Omega)$ with $\partial_x w_i (\pm L) = 0$, for any $i = 1, 2, 3, \ldots, n$.
        \item[\textmd{(\hypertarget{II}{II})}]There exists $\varpi^* = [\varpi_i^{*}, \ldots , \varpi_n^{*}] \in \mathbb{Y} \cap [L^\infty(\Omega)]^n$, $\xi = [\xi_1, \ldots , \xi_n]\in \mathbb{X}$ such that:
        \begin{align}\label{MaTh01-01}
         \displaystyle \frac{1}{\tau}(w_i - w_{i-1}) -\partial_x \varpi_i^{*} - \nu^2 \partial_x^2 w_i + \xi_i + g(w_i) = M_u u_i \mbox{ in } X, i = 1, 2, 3, \ldots, n,
        \end{align}
        with
        \begin{equation*}
\left\{
\begin{array}{l}
\bullet \, \varpi_i^{*} \in \partial f^\varepsilon(\partial_x w_i) \mbox{ in } \mathbb{R}, \mbox{ a.e. in } \Omega, \\[3mm]
\bullet \, \xi_i \in K^\delta(w_i) \mbox{ in } \mathbb{R}, \mbox{ a.e. in } \Omega,
\end{array}
\, i = 1, 2, 3, \ldots, n. \right.
\end{equation*}
         \item[\textmd{(\hypertarget{III}{III})(Energy inequality)}]
         \begin{align*}%\label{MaTh01-02}
          \displaystyle \frac{1}{2\tau}|w_i - w_{i-1}|_X^2 + \mathscr{F}^{(\varepsilon, \delta)}(w_i) - \mathscr{F}^{(\varepsilon, \delta)}(w_{i-1}) \leq \tau M_u^2|u_i|_X^2, i = 1, 2, 3, \ldots, n,
         \end{align*}
         where $\mathscr{F}^{(\varepsilon, \delta)}$ is a functional, called free energy, which is defined as follows: 
         \begin{align}\label{free.en}
\ds z \in X \mapsto \mathscr{F}^{(\varepsilon, \delta)} (z) := &\, \Psi^{(\varepsilon, \delta)}(z) + \int_\Omega G(z)\, dx, \nonumber
\\
= &\, \Phi^\varepsilon(z) + \hat{\mathcal{K}}^\delta(z) + \int_\Omega G(z) dx, \mbox{ for any } z \in X.
\end{align}
    \end{description}
\end{mTh}

\begin{rem}\label{MaTh.rem}
By the definition of subdifferentials, we observe that the equation \eqref{MaTh01-01} is equivalent to the following variational inequality: 
\begin{align*}%\label{MaTh.rem01}
 \displaystyle &\, \frac{1}{\tau}(w_i - w_{i-1}, w_i - z)_X + (g(w_i) - M_u u_i, w_i - z)_X \nonumber 
 \\
 &\,\quad  + \Psi^{(\varepsilon, \delta)}(w_i) - \Psi^{(\varepsilon, \delta)}(z) \leq 0, \mbox{ for any } z \in D(\Psi^{(\varepsilon, \delta)}) \mbox{ and } i = 1, 2, 3, \ldots, n.
\end{align*}
\end{rem}
\medskip

\begin{mTh}\label{MainTh02}
Let $\varepsilon \in [0, 1]$, $\delta \in [0, 1]$, $\{\varepsilon_m \}_{m \in \mathbb{N}} \subset (0, 1]$, $\{\delta_m \}_{m \in \mathbb{N}} \subset (0, 1]$, $u = [u_1, \ldots , u_n] \in \mathbb{X}$, $\{u^m \}_{m \in \mathbb{N}} = \{[u_1^m, \ldots , u_n^m ]\}_{m \in \mathbb{N}} \subset \mathbb{X}$, $w_0 \in Y$, and $\{w_0^m \}_{m \in \mathbb{N}} \subset Y$ be given sequences such that: 
\begin{align}\label{w.i}
\ds \begin{cases}
 \varepsilon_m \to \varepsilon,\, \delta_m \to \delta,\, u^m \to u \mbox{ weakly in } \mathbb{X}, 
 \\[0.75ex]
 \quad \mbox{ and } w_0^m \to w_0 \mbox{ weakly in } Y, \mbox{ as } m \to \infty,
 \\[1.0ex]
 \ds \hat{K}_* := \sup_{m \in \mathbb{N}} \int_\Omega \hat{K}^{\delta_m}(w_0^m)\, dx < \infty.
 \end{cases}
\end{align}
In addition, let $w = [w_1, \ldots , w_n] \in \mathbb{X}$ be the unique solution to \hyperlink{(AC)$^{(\varepsilon, \delta)}$}{(AC)$^{(\varepsilon, \delta)}$} for the forcing term $u$ and the initial data $w_0$, and for any $m \in \mathbb{N}$, let $w^m = [w_1^m, \ldots , w_n^m] \in \mathbb{X}$ be the unique solution to (AC)$^{(\varepsilon_m, \delta_m)}$ for the forcing term $u^m$ and the initial data $w_0^m$. 
Then, it holds that: 
\begin{align}\label{mThConv}
 w^m \to w \mbox{ in } \mathbb{Y}, \mbox{ in } [C^1(\overline{\Omega})]^n, \mbox{ and weakly in } [H^2(\Omega)]^n, \mbox{ as } m \to \infty,
\end{align}
and in particular, 
\begin{align}\label{mThConv00}
 f^{\varepsilon_m}(\partial_x w_i^m) & \to f^\varepsilon(\partial_x w_i) \mbox{ in } X, \mbox{ and in the pointwise sense on } \overline{\Omega}, \nonumber 
 \\
 &\, \quad \mbox{ for any } i = 1, 2, 3, \ldots, n, \mbox{ as } m \to \infty.
\end{align}
\end{mTh}

\begin{mTh}\label{mainTh02}
    Let us assume (\hyperlink{A1l}{A1})--(\hyperlink{A6l}{A6}). Let us fix the constants $ \varepsilon \in [0, 1] $ and $\delta \in [0, 1]$, and fix the initial data $w_0 \in Y$ satisfying $\hat{K}^\delta(w_0) \in L^1(\Omega)$.  Then, the following two items hold.
    \begin{description}
        \item[\textmd{(\hypertarget{III-A}{III-A})}]The problem \hyperlink{(OP)$^{(\varepsilon, \delta)}$}{(OP)$^{(\varepsilon, \delta)}$} has at least one optimal control $ u^{*} = [u_1^{*}, \ldots , u_n^{*}] \in \mathbb{X}$, so that:
\begin{equation*}
    \mathcal{J}^{(\varepsilon, \delta)}(u^{*}) = \min \left\{ \begin{array}{l|l}
        \mathcal{J}^{(\varepsilon, \delta)}(u) & u = [u_1, \ldots , u_n] \in \mathbb{X} 
    \end{array} \right\}.
\end{equation*}
        \item[\textmd{(\hypertarget{III-B}{III-B})}]Let us take the sequences $\{\varepsilon_{m} \}_{m \in \mathbb{N}} \subset (0, 1]$, $\{ \delta_m\}_{m \in \mathbb{N}} \subset (0, 1]$, and take the sequence of initial data $\{w_0^m \}_{m \in \mathbb{N}} $ $ \subset Y $ as in \eqref{w.i}. In addition, for any $m \in \mathbb{N}$, let $ u^{(*, m)} = [u_1^{(*, m)}, \ldots , u_n^{(*, m)}] \in \mathbb{X}$ be the optimal control of (OP)$^{(\varepsilon_m, \delta_m)}$ in the case when the initial pair of corresponding state-system (AC)$^{(\varepsilon_m, \delta_m)}$ is given by $ w_0^m $. Then, there exist a subsequence $ \{ m_k \}_{k \in \mathbb{N}} \subset \{ m \} $ and a function $u^{**} = [u_1^{**}, \ldots , u_n^{**}] \in \mathbb{X}$, such that: 
\begin{equation*}
    \left\{ \hspace{-2ex} \parbox{12cm}{
        \vspace{-1ex}
        \begin{itemize}
            \item $ M_u u^{(*, m_{k})} \to M_u u^{**} $  weakly in $ \mathbb{X} $,  as $ k \to \infty $,

            \item $ u^{**} $ is an optimal control of \hyperlink{(OP)$^{(\varepsilon, \delta)}$}{(OP)$^{(\varepsilon, \delta)}$}.
        \vspace{-1ex}
        \end{itemize}
    } \right.
\end{equation*}
\end{description}
\end{mTh}
\begin{mTh}
    \label{mainTh03}(Necessary condition for \hyperlink{(OP)$^{(\varepsilon, \delta)}$}{(OP)$^{(\varepsilon, \delta)}$} under positive $ \varepsilon, \delta$) 
    Let us assume (\hyperlink{A1l}{A1})--(\hyperlink{A6l}{A6}). 
    Let us fix the constants $ \varepsilon \in (0, 1] $ and $\delta \in (0, 1]$, and fix the initial data $w_0 \in Y$ satisfying $\hat{K}^\delta(w_0) \in L^1(\Omega)$. 
    Let $u^{*} = [u_1^*, \ldots , u_n^*] \in \mathbb{X}$ be an optimal control of \hyperlink{(OP)$^{(\varepsilon, \delta)}$}{(OP)$^{(\varepsilon, \delta)}$}, and let $w^* = [w_1^*, \ldots , w_n^*] \in \mathbb{X}$ be the solution to \hyperlink{(AC)$^{(\varepsilon, \delta)}$}{(AC)$^{(\varepsilon, \delta)}$} for the forcing term $u^*$ and initial data $w_0$. Then, it holds that:
      \begin{equation*}%\label{Thm.5-00}
        \ds M_u \left(p^*_i +  u^*_i \right)= 0 \mbox{ in } X, \mbox{ for any } i = 1, 2, 3, \ldots, n.    
       \end{equation*}
    In this context,  $p^{*} = [p_1^*, \ldots , p_n^*] \in \mathbb{X}$ is a unique solution to the following variational system:
\begin{align*}%\label{Thm.5-01}
  \ds  \frac{1}{\tau} ( p^{*}_i &\, - p^*_{i + 1}, \varphi )_X + \bigl( (f^\varepsilon)''(\partial_x w_i^*)\partial_x p_i^*, \partial_x \varphi \bigr)_X + \nu^2(\partial_x p^*_i, \partial_x \varphi )_X 
    \nonumber
    \\[0.5ex]
    &\, +\bigl( g'(w_i^*) p^{*}_i, \varphi \bigr)_{X} +\bigl( (K^\delta)'(w_i^*)p_i^* , \varphi \bigr)_{X} = \bigl( M_w (w^{*}_i - w_i^{\mbox{\scriptsize ad}}), \varphi \bigr)_{X},
    \\[0.5ex]
     &\, \mbox{ for any $ \varphi \in Y $, and $ i = n, \ldots, 3, 2, 1 $,}
    \nonumber
\end{align*}
subject to the terminal condition:
\begin{equation*}%\label{Thm.5-02}
p^{*}_{n+1} = 0 \mbox{ in } X.
\end{equation*}
\end{mTh}

    \begin{mTh}\label{mainTh04}
        Let us assume (\hyperlink{A1l}{A1})--(\hyperlink{A6l}{A6}), and let us fix an initial data $w_0 \in Y$ satisfying $|w_0| \leq 1$ on $\overline{\Omega}$, i.e. $K^0(w_0) \in L^1(\Omega)$. Also, Let us define a duality map $F : Y \longrightarrow Y^*$, as follows:
        \begin{equation*}
            \langle F\varphi, \psi \rangle_Y :=  (\varphi, \psi)_Y, \mbox{ for any } \varphi, \psi
             \in Y. 
        \end{equation*}
        Then, there exists an optimal control $ u^{\circ} = [u_1^{\circ}, \ldots , u_n^{\circ}] \in \mathbb{X} $ of the problem (OP)$^{(0, 0)} $, together with the solution $ w^\circ = [w_1^{\circ}, \ldots, w_n^{\circ}]$ to the state-system (AC)$^{(0, 0)}$ for the forcing term $ u^\circ$ and initial data $ w_0 $, and moreover, there exist pairs of functions $ p^{\circ} = [p_1^\circ, \ldots , p_n^\circ] \in \mathbb{X} $, and $ \zeta^\circ = [\zeta_1^\circ, \ldots , \zeta_n^\circ] \in \mathbb{Y}^* $, such that:
      \begin{equation}\label{Thm.5-10}
        M_u \left(p^\circ_i + u^\circ_i \right)= 0, \mbox{ in } X, \mbox{ for any } i = 1, 2, 3, \ldots, n,
      \end{equation}
    \begin{equation}\label{Thm.5-10-1}
        p^\circ \in \mathbb{Y} \subset [C(\overline{\Omega})]^n,  
    \end{equation}
    and
    \begin{align}\label{Thm.5-11}
      \ds  \frac{1}{\tau}( p^\circ_i &\, - p_{i+1}^\circ  , \varphi )_{X} + \nu^2 \langle Fp_i^\circ, \varphi \rangle_{Y} - \nu^2(p_i^\circ , \varphi)_X + \bigl( g'(w^\circ_i)p^\circ_i , \varphi)_X + \langle \zeta^\circ_i, \varphi \rangle_{Y} \nonumber
        \\
        &\, = \bigl( M_w (w^\circ_i -w^\mathrm{ad}_i), \varphi \bigr)_{X}, \mbox{ for any $ \varphi \in Y $, and $i = n, \ldots , 3, 2, 1$,}
        \\
        &\, \quad \mbox{ subject to $ p^\circ_{n+1} = 0 $ in $ X $.}
        \nonumber
    \end{align}
Moreover, for any $\gamma_0 \in C^1(\mathbb{R}) \cap W^{1, \infty}(\mathbb{R})$ satisfying $\gamma_0(0) = \gamma'_0(0) = 0$, it follows that: 
\begin{align}\label{Thm.5-13}
 \ds \gamma_0(\partial_x w_i^\circ) &\,\left( \frac{1}{\tau} (p_i^\circ - p_{i+1}^\circ) - \nu^2 \partial_x^2 p_i^\circ + g'(w_i^\circ)p_i^\circ - M_w(w_i^\circ - w_i^{\mathrm{ad}}) \right)  = 0 \mbox{ in } X, \nonumber
 \\
  &\, \mbox{ for any } i = n, \ldots, 3, 2, 1,
\end{align}
and therefore, 
\begin{align}\label{Thm.5-14}
 \ds \frac{1}{\tau} (p_i^\circ - p_{i+1}^\circ) &\, - \nu^2 \partial_x^2 p_i^\circ + g'(w_i^\circ)p_i^\circ = M_w(w_i^\circ - w_i^{\mathrm{ad}}) \mbox{ a.e. in } \{\partial_x w_i^\circ \ne 0 \}, \nonumber
 \\
 &\, \quad \mbox{ for any } i = n, \ldots, 3, 2, 1.
\end{align}
\end{mTh}

\begin{rem}
From \eqref{Thm.5-13} and \eqref{Thm.5-14}, it is observed that: 
\begin{align*}
 \mathrm{spt} (\zeta_i^\circ) \subset \{ \partial_x w_i^\circ = 0\}, \mbox{ for any } i = 1, 2, 3, \ldots , n.
\end{align*}
and 
\begin{align*}
 \ds \zeta_i^\circ &\, = M_w(w_i^\circ - w_i^{\mathrm{ad}}) - \frac{1}{\tau} (p_i^\circ - p_{i+1}^\circ) + \nu^2 \partial_x^2 p_i^\circ - g'(w_i^\circ)p_i^\circ, 
 \\
 &\, \mbox{ a.e. in } \{ \partial_x w_i^\circ \ne 0\} ~ \bigl(\supset \Omega \setminus \mathrm{spt} (\zeta_i^\circ)\bigr), \mbox{ for any } i = n, \ldots, 3, 2, 1.
\end{align*}  
These are to answer the functional expression of the distribution $\zeta^\circ \in \mathbb{Y}^*$, which somehow link to the second unfinished issue, mentioned in the Introduction.   
  
\end{rem}
\medskip

\section{Proof of Key-Theorem} 

In this Section, we give the proof of Key-Theorem. Before the proof, we prepare a Key-Lemma.

\paragraph{\boldmath Key-Lemma.}{
Let us fix $\varepsilon_0 \in [0, 1]$ and $\delta_0 \in [0, 1]$. 
Besides, let $\{\varepsilon_m \}_{m \in \mathbb{N}} \subset (0, 1]$ and $\{\delta_m \}_{m \in \mathbb{N}} \subset (0, 1]$ be given sequences such that $\varepsilon_m \to \varepsilon_0$ and $\delta_m \to \delta_0$, as $m \to \infty$, respectively. Then, for the sequence $\{ \Psi^{(\varepsilon_m, \delta_m)}\}_{m \in \mathbb{N}}$, it holds that: 
\begin{align*}
 \Psi^{(\varepsilon_m, \delta_m)} \to \Psi^{(\varepsilon_0, \delta_0)} \mbox{ on } X, \mbox{ in the sense of Mosco}, \mbox{ as } m \to \infty.
\end{align*}
}

\paragraph{Proof.}{
First, we show the lower-bound condition (\hyperlink{M_lb}{M1}) in Definition \ref{Def.Mosco}. 
Let $z \in X$ and $\{ z^m\}_{m \in \mathbb{N}} \subset X$ such that: 
\begin{align}\label{keylem01}
 z^m \to z \mbox{ weakly in } X, \mbox{ as } m \to \infty.
\end{align}
Then, we may suppose $\varliminf_{m \to \infty} \Psi^{(\varepsilon_m, \delta_m)}(z^m) < \infty $, since the other case is trivial. 
So, by taking a subsequence $\{ m_k\}_{k \in \mathbb{N}} \subset \{m\}$, one can also say that: 
\begin{align}\label{keylem02}
 \varliminf_{m \to \infty}\Psi^{(\varepsilon_m, \delta_m)}(z^m) = \lim_{k \to \infty}\Psi^{(\varepsilon_{m_k}, \delta_{m_k})}(z^{m_k}) < \infty.
\end{align} 
Here, with (\hyperlink{A2l}{A2}), (\hyperlink{A4l}{A4}), \eqref{Phi.eps}, \eqref{K.del}, \eqref{Psi.eps}, \eqref{keylem01}, and \eqref{keylem02} in mind, we further see that: 
\begin{align}\label{keylem03}
 \partial_x z^{m_k} \to \partial_x z \mbox{ weakly in } X, \mbox{ as } k \to \infty,
\end{align}
by taking a subsequence if necessary.
In the light of Remark \ref{f.sub}, (\hyperlink{A2l}{A2}), (\hyperlink{A4l}{A4}), (\hyperlink{Fact2}{Fact\,2}), \eqref{keylem01}--\eqref{keylem03}, and weakly lower semi-continuity of the norms, the lower-bound condition (\hyperlink{M_lb}{M1}) can be verified as follows: 
\begin{align*}
 \ds &\, \varliminf_{m \to \infty}\Psi^{(\varepsilon_{m}, \delta_{m})}(z^{m}) = \lim_{k \to \infty}\Psi^{(\varepsilon_{m_k}, \delta_{m_k})}(z^{m_k}) 
 \\
 &\, \geq \varliminf_{k \to \infty}\int_{\Omega}f^{\varepsilon_{m_k}}(\partial_x z^{m_k}) \, dx + \frac{\nu^2}{2}\varliminf_{k\to \infty}\int_{\Omega} |\partial_x z^{m_k}|^2 \, dx + \varliminf_{k \to \infty} \int_\Omega \hat{K}^{\delta_{m_k}}(z^{m_k})\, dx\\
  &\, \geq \Psi^{(\varepsilon_0, \delta_0)}(z).
\end{align*}

Next, we show the optimality condition (\hyperlink{M_opt}{M2}) in Definition \ref{Def.Mosco}. 
Let us fix any $z \in D(\Psi^{(\varepsilon_0, \delta_0)})$. 
In the light of (\hyperlink{A2l}{A2}), we can say that:
\begin{align}\label{keylem05}
0 & \leq f^{\varepsilon_m}(r) \leq (f^{\varepsilon_m})'(r)r \nonumber
\\
  & \leq C_0(1 + |r|)|r|, \mbox{ for any } r \in \mathbb{R}, \mbox{ and } m \in \mathbb{N}.
\end{align}
With the assumption (\hyperlink{A2l}{A2}) and \eqref{keylem05} in mind, we can infer that: 
\begin{subequations}\label{cor08-5}
\begin{align}\label{cor08-5-1}
 f^{\varepsilon_m}(\partial_x z) \to f^{\varepsilon_0}(\partial_x z) \mbox{ in the pointwise sense a.e. in } \Omega, \mbox{ as } m \to \infty,
\end{align}
\begin{align}\label{cor08-5-2}
 f^{\varepsilon_m}(\partial_x z) \leq C_0 (1 + |\partial_x z|)|\partial_x z|, \mbox{ a.e. in } \Omega, \mbox{ for any } m\in \mathbb{N}.
\end{align}
\end{subequations}
By \eqref{cor08-5} and Lebesgue's dominated convergence theorem, it is verified that: 
\begin{align}\label{cor08-6}
 f^{\varepsilon_m}(\partial_x z) \to f^{\varepsilon_0}(\partial_x z) \mbox{ in } L^1(\Omega), \mbox{ as } m \to \infty.
\end{align} 

Furthermore, we can show: 
\begin{align}\label{cor08-7}
 \hat{K}^{\delta_m}(z) \to \hat{K}^{\delta_0}(z) \mbox{ in } C(\overline{\Omega}), \mbox{ as } m \to \infty.
\end{align}
In fact, the case when $\delta_0 = 0$ is trivial since $\hat{K}^{\delta_m}(z) = \hat{K}^{0}(z) ( = I_{[-1, 1]}(z)) = 0$ a.e. in $\Omega$. 
Meanwhile, when $\delta_0 \in (0, 1]$, the uniform convergence \eqref{cor08-7} is obtained as a consequence of (\hyperlink{A4l}{A4}) and the embedding $Y \subset C(\overline{\Omega})$.

Based on these, let us define $z^m := z$, for any $m \in \mathbb{N}$. 
Taking into account \eqref{cor08-6} and \eqref{cor08-7}, we compute that: 
\begin{align*}%\label{cor08-9}
 \ds \bigl| \Psi^{(\varepsilon_m, \delta_m)}(z^m) - \Psi^{(\varepsilon_0, \delta_0)}(z)\bigr| &\, \leq \int_\Omega \bigl| f^{\varepsilon_m}(\partial_x z) - f^{\varepsilon_0}(\partial_x z)\bigr| \, dx + \int_\Omega \bigl| \hat{K}^{\delta_m}(z) - \hat{K}^{\delta_0}(z) \bigr| \, dx \\
           &\, \leq \int_\Omega \bigl| f^{\varepsilon_m}(\partial_x z) - f^{\varepsilon_0}(\partial_x z)\bigr| dx + 2L\bigl| \hat{K}^{\delta_m}(z) - \hat{K}^{\delta_0}(z) \bigr|_{C(\overline{\Omega})} \\
           &\, \quad  \to 0, \mbox{ as } m \to \infty,
\end{align*}
and therefore,
\begin{align*}
 \Psi^{(\varepsilon_m, \delta_m)}(z^m) \to \Psi^{(\varepsilon_0, \delta_0)}(z), \mbox{ as } m \to \infty.
\end{align*}
\qed
}
\medskip

For efficiency of explanation, we prove the Key-Theorem in accordance with the following two Steps.
    \begin{description}
        \item[\textbf{\boldmath{Step\,$1$:}}]For arbitrary $\varepsilon \in (0, 1]$ and $\delta \in (0, 1]$, the subdifferential $\partial \Psi^{(\varepsilon, \delta)} \subset X \times X$ of the convex function $\Phi^{(\varepsilon, \delta)}$ is decomposed as follows:  
        \begin{align*}
         \partial \Psi^{(\varepsilon, \delta)} = \partial V^\varepsilon +  \partial V^D + K^\delta \mbox{ in } X \times X,  
        \end{align*}
        where $\Psi^{(\varepsilon, \delta)}$, $\Phi^\varepsilon$, $V^\varepsilon$, and $V^D$ are convex functions given in \eqref{Main.Vd}, \eqref{Phi.eps}, \eqref{Psi.eps}, and \eqref{V.eps}, respectively, and $K^\delta$ is the function as in (\hyperlink{A4l}{A4}).
        \item[\textbf{\boldmath{Step\,$2$:}}]For arbitrary $\varepsilon \in [0, 1]$ and $\delta \in [0, 1]$, the subdifferential $\partial \Psi^{(\varepsilon, \delta)} \subset X \times X$ of the convex function $\Phi^{(\varepsilon, \delta)}$ is decomposed as follows: 
        \begin{align*}
         \partial \Psi^{(\varepsilon, \delta)} = \partial V^\varepsilon +  \partial V^D + \mathcal{K}^{\delta} \mbox{ in } X \times X, 
        \end{align*}
        where $\mathcal{K}^{\delta}$ is the operator given in \eqref{K.del1}.
    \end{description}  
\medskip

\noindent
\textbf{\boldmath\underline{Verification of Step\,$1$.}}~

Let us fix $\varepsilon \in (0, 1]$ and $\delta \in (0, 1]$. 
The assumption (\hyperlink{A4l}{A4}) guarantees $D( \hat{\mathcal{K}}^\delta) = X$,
i.e. 
\begin{align}\label{step2-01}
 \ds \mathrm{int}D( \hat{\mathcal{K}}^\delta) = X.
\end{align}
Due to \eqref{step2-01}, the decomposition of the subdifferential $\partial \Psi^{(\varepsilon, \delta)} \subset X \times X$ of the convex function $\Psi^{(\varepsilon, \delta)} = \Phi^\varepsilon + \hat{\mathcal{K}}^\delta$ will be a straightforward consequence of \cite[Theorem 2.10]{MR2582280}, and \cite[Corollary 2.11]{MR0348562}, i.e.
\begin{align}
 \partial \Psi^{(\varepsilon, \delta)} = \partial \Phi^\varepsilon + K^\delta \mbox{ in } X \times X.
\end{align} 

Hence, it is sufficient to prove that the subdifferential $\partial \Phi^\varepsilon \subset X \times X $ of the convex function $\Phi^\varepsilon$ is decomposed as follows: 
\begin{align}
\partial \Phi^\varepsilon = \partial V^\varepsilon + \partial V^D \mbox{ in } X \times X.
\end{align}

First, we verify the following inclusion:
\begin{align}\label{step1-00}
 \partial\Phi^\varepsilon \subset \partial V^\varepsilon + \partial V^D \mbox{ in } X \times X.
\end{align}

Let us take any $[z, z^*] \in \partial \Phi^\varepsilon$ in $X \times X$, and apply Remark \ref{lem.op.ap} to the case when $\alpha \equiv 1$ and $\beta \equiv \nu^2$.
Then, we have: 
\begin{align}\label{step1-01}
z^* & = - \partial_x \left( (f^\varepsilon)'(\partial_x z) + \nu^2 \partial_x z\right) \mbox{ in } X, \mbox{ and } \nonumber
\\
& (f^\varepsilon)'(\partial_x z) + \nu^2 \partial_x z \in H_0^1(\Omega). 
\end{align}
Now, we define a function $\mu_\varepsilon^\nu : \mathbb{R} \longrightarrow \mathbb{R}$ as follows: 
\begin{align}\label{step1-02}
 \mu_\varepsilon^\nu(r) := (f^\varepsilon)'(r) + \nu^2 r, \mbox{ for any } r \in \mathbb{R}. 
\end{align}
On this basis, we infer that: 
\begin{align}\label{step1-03}
 \mu_\varepsilon^\nu(\partial_x z) & = (f^\varepsilon)'(\partial_x z) + \nu^2 \partial_x z = \left((f^\varepsilon)' + \nu^2 I_d \right)(\partial_x z) \mbox{ in } Y.
\end{align}
Based on Remark \ref{rem.f.lip} and \eqref{step1-03}, we can apply the generalized chain rule in BV-theory \cite[Theorem 3.99]{MR1857292}, and can infer that:
\begin{align}\label{step1-05}
 \partial_x z = (\mu_\varepsilon^\nu)^{-1}(\mu_\varepsilon^\nu(\partial_x z)) \in Y, \mbox{ i.e. } z \in H^2(\Omega).
\end{align}
Furthermore, having in mind \eqref{step1-01}, \eqref{step1-03}, \eqref{step1-05}, and applying Remark \ref{Rem.W_bt} to the case when $\beta \equiv \nu^2$, one can see that: 
\begin{align}\label{step1-06}
 \partial_x z \in H_0^1(\Omega) \mbox{ and } [z, -\nu^2\partial_x^2 z] \in \partial V^D \mbox{ in } X \times X.
\end{align}
%where $\beta \equiv \nu^2$. in Remark \ref{Rem.W_bt}. 
By using (\hyperlink{A2l}{A2}), \eqref{step1-01}, \eqref{step1-06}, and the integration by part, we observe that: 
\begin{align}\label{step1-07}
 (f^\varepsilon)'(\partial_x z) \in H_0^1(\Omega),
\end{align}
and
\begin{align*}%\label{step1-08}
 \ds \int_\Omega & -\partial_x \left( (f^\varepsilon)'(\partial_x z)\right)(\varphi - z)\, dx = \int_\Omega (f^\varepsilon)'(\partial_x z)\partial_x (\varphi - z)\, dx \\
 & \leq \int_\Omega f^\varepsilon(\partial_x \varphi)\, dx - \int_\Omega f^\varepsilon(\partial_x z)\, dx, \mbox{ for any } \varphi \in Y, 
\end{align*}
i.e.  
\begin{align}\label{step1-09}
[z, -\partial_x \left( (f^\varepsilon)'(\partial_x z)\right)] \in \partial V^\varepsilon \mbox{ in } X \times X.
\end{align}
\eqref{step1-05}--\eqref{step1-09} enable us to verify the inclusion \eqref{step1-00}. 

Now, from the maximality of $\partial \Phi^\varepsilon$ in $X \times X$, we can see the coincidence $\partial \Phi^\varepsilon = \partial V^\varepsilon + \partial V^D$ in $X \times X$. 
\medskip

\noindent
\textbf{\boldmath\underline{Verification of Step\,$2$.}}~ 

Let us fix $\varepsilon \in [0, 1]$ and $\delta \in [0, 1]$. 
By the general theory of the convex analysis \cite[Chapter 1]{MR1727362}, we immediately obtain that $\partial \Psi^{(\varepsilon, \delta)} \supset \partial V^\varepsilon + \partial V^D + \mathcal{K}^\delta$ in $X \times X$. 
This inclusion implies the monotonicity of $\partial V^\varepsilon + \partial V^D + \mathcal{K}^\delta$ in $X \times X$. 
So, we next see the maximality of $\partial V^\varepsilon + \partial V^D + \mathcal{K}^\delta$, by verifying:
\begin{align*}
 (\partial V^\varepsilon + \partial V^D + \mathcal{K}^\delta + I_d)X = X,
\end{align*}
and by applying Minty's theorem (cf. \cite[Theorem 2.2]{MR2582280}). 

Since, the inclusion $(\partial V^\varepsilon + \partial V^D + \mathcal{K}^\delta + I_d)X \subset X$ is trivial, it is sufficient to prove the converse inclusion. 
Let us fix $h \in X$. 
Then, by Step\,$1$, we can configure a class of functions $\{ z^{(\tilde{\varepsilon}, \tilde{\delta})} \}_{\tilde{\varepsilon}, \tilde{\delta} \in (0, 1]} \subset Y$, by setting: 
\begin{align*}
 \{ z^{(\tilde{\varepsilon}, \tilde{\delta})} := (\partial V^{\tilde{\varepsilon}} + \partial V^D + K^{\tilde{\delta}} + I_d)^{-1} (h)\}_{\tilde{\varepsilon}, \tilde{\delta} \in (0, 1]} \mbox{ in } X,
\end{align*}
i.e. 
\begin{align}\label{claim7-01}
 h - z^{(\tilde{\varepsilon}, \tilde{\delta})} = (\partial V^{\tilde{\varepsilon}} + \partial V^D + K^{\tilde{\delta}})(z^{(\tilde{\varepsilon}, \tilde{\delta})}) = \partial \Psi^{(\tilde{\varepsilon}, \tilde{\delta})}(z^{(\tilde{\varepsilon}, \tilde{\delta})}) \mbox{ in } X, \mbox{ for any } \tilde{\varepsilon}, \tilde{\delta} \in (0, 1].
\end{align}
In the light of \eqref{step1-06}, \eqref{step1-09}, and \eqref{claim7-01}, there exist $\varpi^{\tilde{\varepsilon}} \in Y \cap L^\infty(\Omega)$ and $\xi^{\tilde{\delta}} \in X$ such that
\begin{align}\label{s2-02}
  -\partial_x \varpi^{\tilde{\varepsilon}} - \nu^2 \partial_x^2 z^{(\tilde{\varepsilon}, \tilde{\delta})} &\, + \xi^{\tilde{\delta}} = h - z^{(\tilde{\varepsilon}, \tilde{\delta})} \mbox{ in } X, \nonumber
 \\
 &\, \mbox{ where } \varpi^{\tilde{\varepsilon}} = \partial f^{\tilde{\varepsilon}}(\partial_x z^{(\tilde{\varepsilon}, \tilde{\delta})}) \mbox{ and } \xi^{\tilde{\delta}} = K^{\tilde{\delta}}(z^{(\tilde{\varepsilon}, \tilde{\delta})}) \mbox{ a.e. in } \Omega.
\end{align}
Using \eqref{step1-06}, \eqref{step1-07}, and the integration by part, we can see that: 
\begin{align}\label{claim7-02}
 \int_\Omega \varpi^{\tilde{\varepsilon}} \partial_x \varphi\, dx  &\, + \int_\Omega \nu^2 \partial_x z^{(\tilde{\varepsilon}, \tilde{\delta})}\partial_x \varphi\, dx + \int_\Omega \xi^{\tilde{\delta}}\varphi + \int_\Omega z^{(\tilde{\varepsilon}, \tilde{\delta})}\varphi \, dx \nonumber
 \\
 &\, = \int_\Omega h\varphi\, dx, \mbox{ for any }\varphi \in Y, \mbox{ and  any } \tilde{\varepsilon}, \tilde{\delta} \in (0, 1].
\end{align}
In the variational form \eqref{claim7-02}, let us put $\varphi = z^{(\tilde{\varepsilon}, \tilde{\delta})}$. 
Then, with Remark \ref{K.rem}, \eqref{keylem05}, and Young's inequality in mind, we deduce that: 
\begin{align}\label{claim7-03}
 \frac{1}{2}|z^{(\tilde{\varepsilon}, \tilde{\delta})}|_X^2 + \nu^2|\partial_x z^{(\tilde{\varepsilon}, \tilde{\delta})}|_X^2 \leq \frac{1}{2}|h|_X^2, \mbox{ for any } \tilde{\varepsilon}, \tilde{\delta} \in (0, 1],
\end{align}
so that
\begin{align}\label{claim7-04}
 \ds |z^{(\tilde{\varepsilon}, \tilde{\delta})}|_Y^2 \leq \frac{1}{1 \wedge (2\nu^2)}|h|_X^2, \mbox{ for any } \tilde{\varepsilon}, \tilde{\delta} \in (0, 1].
\end{align}  
On account of \eqref{claim7-04}, we find a function $z \in Y$ and sequences $\{\tilde{\varepsilon}_m \}_{m \in \mathbb{N}} \subset \{ \tilde{\varepsilon}\}$ with $\tilde{\varepsilon}_m \to \varepsilon$ as $m \to \infty$, and $\{\tilde{\delta}_m \}_{m \in \mathbb{N}} \subset \{ \tilde{\delta}\}$ with $\tilde{\delta}_m \to \delta$ as $m \to \infty$, such that:
\begin{align}\label{claim7-05}
 z^{(\tilde{\varepsilon}_m, \tilde{\delta}_m)} \to z \mbox{ in } X, \mbox{ weakly in } Y, \mbox{ as } m \to \infty.
\end{align}
\begin{comment}
In the light of Key-Lemma, \eqref{claim7-01}, and \eqref{claim7-05}, we can apply (\hyperlink{Fact1}{Fact\,1}) to see that: 
\begin{align*}
 h - z \in \partial \Psi^{(\varepsilon, \delta)}(z) \mbox{ in } X,
\end{align*}
and 
\begin{align}\label{claim7-06}
 \Psi^{(\tilde{\varepsilon}_m, \tilde{\delta}_m)}(z^{(\tilde{\varepsilon}_m, \tilde{\delta}_m)}) \to \Psi^{(\varepsilon, \delta)}(z), \mbox{ as } m \to \infty.
\end{align}
\end{comment}

Furthermore, let us multiply $\xi^{\tilde{\delta}_m}$ the both sides of \eqref{claim7-01}. 
Then, having in mind (\hyperlink{A4l}{A4}), \eqref{keylem05}, \eqref{claim7-03}, and Young's inequality, we deduce that: 
\begin{align}\label{claim7-10}
 \frac{1}{2}|\xi^{\tilde{\delta}_m}|_X^2 \leq |h|_X^2 + |z^{(\tilde{\varepsilon}_m, \tilde{\delta}_m)}|_X^2 \leq 2|h|_X^2, \mbox{ for any } m \in \mathbb{N}, 
\end{align}
via the following calculation: 
\begin{align}\label{claim7-11}
\ds &\, \int_\Omega \left( -\partial_x \varpi^{\tilde{\varepsilon}_m} - \nu^2 \partial_x^2 z^{(\tilde{\varepsilon}_m, \tilde{\delta}_m)}\right)\xi^{\tilde{\delta}_m}\, dx \nonumber
\\
 &\, = \int_\Omega \left( -\partial_x \bigl( \varpi^{\tilde{\varepsilon}_m} + \nu^2 \partial_x z^{(\tilde{\varepsilon}_m, \tilde{\delta}_m)}\bigr)\right)\xi^{\tilde{\delta}_m}\, dx \nonumber
\\
 &\, = \int_\Omega \left( -\partial_x \bigl( (f^{\tilde{\varepsilon}_m})'(\partial_x z^{(\tilde{\varepsilon}_m, \tilde{\delta}_m)}) + \nu^2 \partial_x z^{(\tilde{\varepsilon}_m, \tilde{\delta}_m)}\bigr)\right)K^{\tilde{\delta}_m}(z^{(\tilde{\varepsilon}_m, \tilde{\delta}_m)})\, dx \nonumber
\\
  &\, = \int_\Omega \left( (f^{\tilde{\varepsilon}_m})'(\partial_x z^{(\tilde{\varepsilon}_m, \tilde{\delta}_m)}) + \nu^2 \partial_x z^{(\tilde{\varepsilon}_m, \tilde{\delta}_m)}\right)\partial_x K^{\tilde{\delta}_m}(z^{(\tilde{\varepsilon}_m, \tilde{\delta}_m)})\, dx \nonumber
  \\
  &\, = \int_\Omega (K^{\tilde{\delta}_m})'(z^{(\tilde{\varepsilon}_m, \tilde{\delta}_m)})(f^{\tilde{\varepsilon}_m})'(\partial_x z^{(\tilde{\varepsilon}_m, \tilde{\delta}_m)})\partial_x z^{(\tilde{\varepsilon}_m, \tilde{\delta}_m)} \, dx \nonumber
  \\
  &\, \qquad + \int_\Omega \nu^2 (K^{\tilde{\delta}_m})'(\partial_x z^{(\tilde{\varepsilon}_m, \tilde{\delta}_m)})|\partial_x z^{(\tilde{\varepsilon}_m, \tilde{\delta}_m)}|^2 \, dx \nonumber
  \\
  &\, \geq 0, \mbox{ for any } m \in \mathbb{N}.
\end{align}
From the estimate \eqref{claim7-10}, it is observed that:
\begin{align}\label{claim7-12}
\xi^{\tilde{\delta}_m} \to \xi \mbox{ weakly in } X, \mbox{ as } m \to \infty, \mbox{ for some } \xi \in X, 
\end{align}
by taking a subsequence if necessary. 

Additionally, by virtue of \eqref{s2-02}, \eqref{claim7-03}, and \eqref{claim7-10}, one can observe that :
\begin{description}
    \item[\textmd{$(\hypertarget{lozenge1}{\lozenge\,1})$}] the sequence $\{ \partial_x \varpi^{\tilde{\varepsilon}_m} + \nu^2 \partial_x^2 z^{(\tilde{\varepsilon}_m, \tilde{\delta}_m)}\}_{m \in \mathbb{N}} =\{ \partial_x(\varpi^{\tilde{\varepsilon}_m} + \nu^2 \partial_x z^{(\tilde{\varepsilon}_m, \tilde{\delta}_m)})\}_{m \in \mathbb{N}} $ is bounded in $ X $. 
\end{description}

Note that \eqref{step1-02} leads to: 
\begin{align}\label{s2-03}
 \partial_x z^{(\tilde{\varepsilon}_m, \tilde{\delta}_m)} &\, = (\mu_{\tilde{\varepsilon}_m}^\nu)^{-1}\bigl(\mu_{\tilde{\varepsilon}_m}^\nu(\partial_x z^{(\tilde{\varepsilon}_m, \tilde{\delta}_m)}) \bigr) \\ \nonumber
 &\, = (\mu_{\tilde{\varepsilon}_m}^\nu)^{-1}(\varpi^{\tilde{\varepsilon}_m} + \nu^2 \partial_x z^{(\tilde{\varepsilon}_m, \tilde{\delta}_m)}), \mbox{ for any } m \in \mathbb{N}.
\end{align}
With Remark \ref{rem.f.lip} and \eqref{s2-03} in mind, we apply the generalized chain rule in BV-theory \cite[Theorem 3.99]{MR1857292}, and we infer that:
\begin{align}\label{s2-04}
 |\partial_x^2 z^{(\tilde{\varepsilon}_m, \tilde{\delta}_m)}|_X &\, = \bigl|\bigl((\mu_{\tilde{\varepsilon}_m}^\nu)^{-1}\bigr)'(\varpi^{\tilde{\varepsilon}_m} + \nu^2 \partial_x z^{(\tilde{\varepsilon}_m, \tilde{\delta}_m)})\partial_x(\varpi^{\tilde{\varepsilon}_m} + \nu^2 \partial_x z^{(\tilde{\varepsilon}_m, \tilde{\delta}_m)}) \bigr|_X \nonumber
 \\
 &\, \ds \leq \frac{1}{\nu^2}|\partial_x(\varpi^{\tilde{\varepsilon}_m} + \nu^2 \partial_x z^{(\tilde{\varepsilon}_m, \tilde{\delta}_m)})|_X, \mbox{ for any } m \in \mathbb{N}.
\end{align}
From $(\hyperlink{lozenge1}{\lozenge\,1})$ \eqref{claim7-04}, and \eqref{s2-04}, it is deduced that: 
\begin{description}
    \item[\textmd{$(\hypertarget{lozenge2}{\lozenge\,2})$}]the sequence $\{ \partial_x^2 z^{(\tilde{\varepsilon}_m, \tilde{\delta}_m)}\}_{m \in \mathbb{N}} $ is bounded in $ X $, and hence $\{ z^{(\tilde{\varepsilon}_m, \tilde{\delta}_m)}\}_{m \in \mathbb{N}}$ is bounded in $ H^2(\Omega)$.
\end{description}
As a consequence of the one-dimensional compact embedding $H^2(\Omega) \subset C^1(\overline{\Omega})$, \eqref{claim7-05}, and $(\hyperlink{lozenge2}{\lozenge\,2})$, it is observed that: 
\begin{align}\label{s2-05}
 z^{(\tilde{\varepsilon}_m, \tilde{\delta}_m)} \to &\, z \mbox{ in } Y, \mbox{ in } C^1(\overline{\Omega}), \mbox{ and weakly in } H^2(\Omega), \mbox{ as } m \to \infty,
\end{align}
by taking a subsequence if necessary. 

Note that by the assumption (\hyperlink{A2l}{A2}) and \eqref{s2-05}, we can compute that: 
\begin{align*}%\label{claim4-09}
 \ds |\varpi^{\tilde{\varepsilon}_m} |_X^2 & = |(f^{\tilde{\varepsilon}_m})'(\partial_x z^{(\tilde{\varepsilon}_m, \delta_m)})|_X^2 \nonumber
 \\
 & \leq \int_\Omega \bigl(C_0 (1 + |\partial_x z^{(\tilde{\varepsilon}_m, \tilde{\delta}_m)}|) \bigr)^2 dx \nonumber
 \\
 & \leq 2C_0^2 \int_\Omega (1 + |\partial_x z^{(\tilde{\varepsilon}_m, \tilde{\delta}_m)}|^2) dx\nonumber
 \\
 & \leq 2C_0^2\bigl( 2L + \sup_{m \in \mathbb{N}}|\partial_x z^{(\tilde{\varepsilon}_m, \tilde{\delta}_m)}|_X^2\bigr) < \infty, \mbox{ for any } m \in \mathbb{N}.
\end{align*}
Based on these, it enables us to say 
\begin{align}\label{claim7-09}
 \varpi^{\tilde{\varepsilon}_m} \to \varpi \mbox{ weakly in } X, \mbox{ as } m \to \infty, \mbox{ for some } \varpi \in X, 
\end{align}
by taking a subsequence if necessary. 
Besides, from $(\hyperlink{lozenge1}{\lozenge\,1})$ and $(\hyperlink{lozenge2}{\lozenge\,2})$, it is deduced that: 
\begin{description}
    \item[\textmd{$(\hypertarget{lozenge3}{\lozenge\,3})$}] the sequence $\{ -\partial_x \varpi^{\tilde{\varepsilon}_m}\}_{m \in \mathbb{N}} = \{ \omega^{\tilde{\varepsilon}_m} \}_{m \in \mathbb{N}} $ is bounded in $ X $; 
\end{description}
and we obtain that 
\begin{align}\label{s2-06}
\omega^{\tilde{\varepsilon}_m} \to \omega \mbox{ weakly in } X, \mbox{ as } m \to \infty, \mbox{ for some } \omega \in X,
\end{align}
by taking a subsequence if necessary. 

In the light of (\hyperlink{A2l}{A2}), (\hyperlink{A4l}{A4}), (\hyperlink{Fact1}{Fact\,1}), (\hyperlink{Fact2}{Fact\,2}), \eqref{Sgn^d}, \eqref{claim7-12}, and \eqref{s2-05}--\eqref{s2-06}, it is inferred that: 
\begin{subequations}\label{claim7-13}
\begin{align}\label{claim7-13a}
 \xi \in \mathcal{K}^\delta(z),\, \omega \in \partial V^\varepsilon (z), \mbox{ and } -\nu^2 \partial_x^2 z \in \partial V^D(z) \mbox{ in } X. 
\end{align}
Additionally, taking into account Remark \ref{f.sub}, (\hyperlink{Fact1}{Fact\,1}), \eqref{s2-05}, \eqref{claim7-09}, \cite[Theorem 2.10]{MR2582280}, and \cite[Corollary 2.11]{MR0348562}, one can obtain that:
\begin{equation}\label{claim7-13b}
    \varpi \in \partial f^\varepsilon(\partial_x z) \mbox{ a.e. in $ \Omega $.}
\end{equation}
\end{subequations}

In the meantime, letting $m \to \infty$ in \eqref{claim7-02} yields that: 
\begin{align*}
 \int_\Omega \varpi \partial_x \varphi\, dx  &\, + \int_\Omega \nu^2 \partial_x z\partial_x \varphi\, dx + \int_\Omega \xi \varphi + \int_\Omega z\varphi \, dx = \int_\Omega h\varphi\, dx, \mbox{ for any }\varphi \in Y.
\end{align*}
This equation and \eqref{claim7-13} imply: 
\begin{align*}
-\partial_x \varpi - \nu^2 \partial_x^2 z + \xi + z = h \mbox{ in } X,
\end{align*}
i.e.
\begin{align*}
 (\partial V^\varepsilon + \partial V^D + \mathcal{K}^\delta + I_d)(z) \ni h. 
\end{align*}
By applying Minty's theorem (cf. \cite[Theorem 2.2]{MR2582280}), $\partial V^\varepsilon + \partial V^D + \mathcal{K}^\delta$ is a maximal monotone. 
Now, from this maximality, we can see the coincidence $ \partial \Psi^{(\varepsilon, \delta)} = \partial V^\varepsilon + \partial V^D + \mathcal{K}^\delta $ in $X \times X$.

Based on Step\,$1$--Step\,$2$, we conclude that for any $\varepsilon \in [0, 1]$ and $\delta \in [0, 1]$, $[z, z^*] \in \Psi^{(\varepsilon, \delta)}$ in $X \times X$ if and only if there exist $\varpi^* \in Y$ and $\xi \in X$ such that:
\begin{itemize}
\item $z \in H^2(\Omega)$, with $(\varpi^* + \nu^2 \partial_x z)(\pm L) = 0$; 
\item $z^* = -\partial_x \varpi^* - \nu^2 \partial_x^2 z + \xi$ in $X$, where $\varpi^* \in \partial f^\varepsilon(\partial_x z)$ and $\xi \in K^\delta(z)$ a.e. in $\Omega$.
\end{itemize}
Hence, we finish the proof of Key-Theorem.
\qed

\begin{rem}
The Key-Theorem is obtained on the basis of some previous works \cite{MR4279540,MR2101878}. 
In fact, the proof of Step\,$1$ is referred to the proving method of \cite[Main Theorem]{MR4279540}, and the case when $\varepsilon = 0$ and $\delta \in (0, 1]$ is verified as a consequence of \cite[Main Theorem]{MR4279540} and \cite[Theorem 3.1]{MR2101878}. 
However, the Key-Theorem covers various approximating sequences $\{ f^\varepsilon\}_{\varepsilon \in (0, 1]}$ and $\{ K^\delta\}_{\delta \in (0, 1]}$ in the range of the assumptions (\hyperlink{A2l}{A2}) and (\hyperlink{A4l}{A4}) which includes the previous setting adopted in \cite{MR2836557,MR2459669,MR2509574,MR2101878,MR3661429}. 

In this light, it can be said that our Key-Theorem provides a general theory of approximating method for $-\partial_x (\frac{\partial_x w}{|\partial_x w|} + \nu^2 \partial_x w) + \partial I_{[-1, 1]}(w)$, and also brings a key-answer for unfinished issue mentioned in Introduction. 
\end{rem}
\medskip

\section{Proof of Main Theorem \ref{mainTh01}} 

In this Section, we give the proof of Main Theorem \ref{mainTh01}. 
Let us fix $\varepsilon \in [0, 1]$ and $\delta \in [0, 1]$. 
Let us fix a forcing term $u = [u_1, \ldots , u_n] \in \mathbb{X}$, and an initial data $w_0 \in Y$ satisfying $\hat{K}^\delta(w_0) \in L^1(\Omega)$.  
Let us fix $i \in \{ 1, 2, 3, \ldots, n\}$. 

On this basis, we define a functional $\mathscr{G} : X \longrightarrow (-\infty, \infty]$, by letting: 
\begin{align}\label{maTh01-01}
    w & \in X \mapsto \mathscr{G}(w) 
    \nonumber
    \\
    &:= \left\{
        \begin{array}{lll}
            \multicolumn{2}{l}{\displaystyle \frac{1}{2\tau}\int_\Omega|w-w_{i-1}|^2\, dx + \Psi^{(\varepsilon, \delta)}(w) + \int_\Omega G(w)\, dx - (M_u u_i, w)_X,}
            \\[2ex]
            & \mbox{if $ w \in Y$,}
            \\[2ex]
            \infty, & \mbox{otherwise.}
        \end{array}
    \right. 
\end{align}
Since the assumptions (\hyperlink{A3l}{A3}), (\hyperlink{A5l}{A5}), and Remark \ref{g.prime} imply 
\begin{align*}
 \ds &\, \frac{d^2}{d w^2}\left( \frac{1}{4\tau}|w - \tilde{w}_0|^2 + G(w) \right) = \frac{1}{2\tau} + g'(w) > \frac{1}{2\tau^*} - C_g > \frac{C_g}{\tau^*}\left( \frac{1}{8(C_g + 1)} - \tau^*\right) = 0, 
 \\
 \ds &\, \quad \mbox{ for any fixed } \tilde{w}_0 \in \mathbb{R},
\end{align*}
and the functional:  
\begin{align*}%\label{maTh01-02}
 w &\, \in D(\mathscr{G}) = Y \subset X \mapsto \mathscr{G}(w) \nonumber
 \\
 \ds &\, = \frac{1}{4\tau}|w - w_{i-1}|_X^2 + \int_\Omega \left(\frac{1}{4\tau}|w - w_{i-1}|^2 + G(w) \right)\, dx + \Psi^{(\varepsilon, \delta)}(w) - (M_u u_i, w)_X,
\end{align*}
is proper, l.s.c., strictly convex, and coercive on X.

Based on these, we find a unique minimizer of $\mathscr{G}$, denoted by $w^* \in X$, by applying \cite[Proposition 1.2, Chapter II]{MR1727362}. 
Since $w^*$ is the minimizer of $\mathscr{G}$ and $\Psi^{(\varepsilon, \delta)} $ is a convex on $X$, we can compute that: 
\begin{align}\label{maTh01-05}
 \ds 0 &\, \leq \frac{1}{\lambda}\bigl( \mathscr{G}(w^* + \lambda (\varphi- w^*)) - \mathscr{G}(w^*) \bigr) \nonumber 
 \\
       &\, \leq \frac{1}{\tau}\int_\Omega (\varphi - w^*)(w^* - w_{i-1})\, dx + \frac{\lambda}{2\tau} \int_\Omega |\varphi - w^*|^2\, dx \nonumber
       \\
       &\, \quad + \frac{1}{\lambda}\int_\Omega \bigl(G(w^* + \lambda (\varphi -w^* )) - G(w^*)  \bigr)\, dx \nonumber 
       \\
       &\, \quad  + \Psi^{(\varepsilon, \delta)}(\varphi) - \Psi^{(\varepsilon, \delta)}(w^*) - (M_u u_i, \varphi - w^*)_X, \mbox{ for any } \varphi \in Y, \mbox{ and } \lambda \in (0, 1).       
\end{align}
Letting $\lambda \downarrow 0$ in \eqref{maTh01-05}, it is inferred that: 
\begin{align*}%\label{maTh01-06}
\ds  \Psi^{(\varepsilon, \delta)}(\varphi) - \Psi^{(\varepsilon, \delta)}(w^*) &\, \leq \left(-\left(\frac{1}{\tau}(w^* - w_{i-1}) + g(w^*)-  M_u u_i\right), \varphi - w^* \right)_X, 
\\
&\, \mbox{ for any } \varphi \in Y,  
\end{align*}
i.e.
\begin{align}\label{maTh01-07}
 \ds -\left(\frac{1}{\tau}(w^* - w_{i-1}) + g(w^*)- M_u u_i \right) \in \partial \Psi^{(\varepsilon, \delta)}(w^*) \mbox{ in } X.
\end{align}
Here, on account of \eqref{maTh01-07} and Key-Theorem, $w_i \in H^2(\Omega)$, and there exist $\varpi_i^* \in Y \cap L^\infty(\Omega)$ and $\xi_i \in X$ such that: 
\begin{align}\label{maTh01-07-1}
 \ds \frac{1}{\tau}(w^* - w_{i-1})_X - \partial_x \varpi_i^* - \nu^2\partial_x^2 w^* + g(w^*) + \xi_i = M_u u_i, \mbox{ in } X,
\end{align}
with
        \begin{equation*}
\left\{
\begin{array}{l}
\bullet \, \varpi^{*}_i \in \partial f^\varepsilon(\partial_x w^*) \mbox{ in } \mathbb{R}, \mbox{ a.e. in } \Omega, \\[3mm]
\bullet \, \xi_i \in K^\delta(w^*) \mbox{ in } \mathbb{R}, \mbox{ a.e. in } \Omega,
\end{array}
\right.
i = 1, 2, 3, \ldots , n.
\end{equation*}
The equation \eqref{maTh01-07-1} implies that the state-system \hyperlink{(AC)$^{(\varepsilon, \delta)}$}{(AC)$^{(\varepsilon, \delta)}$} has a solution $w = [w_1, \ldots , w_n] \in [H^2(\Omega)]^n$. 
\medskip

Next, we prove the energy inequality as in (III). 
Let us fix $i \in \{ 1, 2, 3, \ldots, n\}$ and multiply $w_i - w_{i-1}$ the both sides of \eqref{MaTh01-01}, we arrive at:
\begin{align}\label{maTh01-14}
 \ds &\, \frac{1}{\tau}|w_i - w_{i-1}|_X^2 + (\varpi_i^{*}, \partial_x(w_i - w_{i-1}))_X + \nu^2(\partial_x w_i, \partial_x(w_i - w_{i-1}))_X \nonumber
 \\
 &\, \quad + (\xi_i, w_i - w_{i-1})_X + \bigl(g(w_i) , w_i -w_{i-1} \bigr)_X = M_u (u_i, w_i - w_{i-1})_X,
\end{align} 
        with
        \begin{equation*}
\left\{
\begin{array}{l}
\bullet \, \varpi_i^{*} \in \partial f^\varepsilon(\partial_x w_i) \mbox{ in } \mathbb{R}, \mbox{ a.e. in } \Omega, \\[3mm]
\bullet \, \xi_i \in K^\delta(w_i) \mbox{ in } \mathbb{R}, \mbox{ a.e. in } \Omega.
\end{array}
 \right.
\end{equation*}
Also, by applying the assumption (\hyperlink{A3l}{A3}), Remark \ref{g.prime}, and Taylor's theorem, we have: 
\begin{align}\label{maTh01-02}
 G(w_{i-1}) &\, \geq G(w_i) + g(w_i)(w_{i-1} - w_i) - \frac{1}{2}|[g']^-|_{L^\infty(\mathbb{R})}(w_{i-1} - w_i)^2 \nonumber
 \\
 &\, \geq G(w_i) -g(w_i)(w_i - w_{i-1}) - \frac{C_g}{2}|w_i - w_{i-1}|^2, \mbox{ a.e. in } \Omega. 
\end{align}

By using \eqref{Sgn^d}, \eqref{maTh01-02}, the assumptions (\hyperlink{A1l}{A1})--(\hyperlink{A4l}{A4}), and Young's inequality, we can deduce from \eqref{maTh01-14} that: 
\begin{align}\label{maTh01-15}
 \ds &\, \frac{3}{4\tau}|w_i - w_{i-1}|_X^2 + \frac{\nu^2}{2}|\partial_x w_i|_X^2 - \frac{\nu^2}{2}|\partial_x w_{i-1}|_X^2 \nonumber 
 \\
 &\, \quad + \int_\Omega f^\varepsilon(\partial_x w_i)\, dx - \int_\Omega f^\varepsilon (\partial_x w_{i-1})\, dx \nonumber 
 \\
 &\, \quad + \int_\Omega G(w_i)\, dx - \int_\Omega G(w_{i-1})\, dx - \frac{C_g}{2}|w_i - w_{i-1}|_X^2 \nonumber 
 \\
 &\, \quad + \int_\Omega \hat{K}^\delta(w_i)\, dx - \int_\Omega \hat{K}^\delta(w_{i-1})\, dx \leq \tau M_u^2 |u_i|_X^2,
\end{align}
via
\begin{align*}
 \ds (\varpi_i^*, \partial_x (w_i - w_{i-1}))_X \geq \int_\Omega f^\varepsilon(\partial_x w_i)\, dx - \int_\Omega f^\varepsilon (\partial_x w_{i-1})\, dx,
 \end{align*}
 \begin{align*}
  \ds \nu^2\bigl(\partial_x w_i, \partial_x (w_i - w_{i-1}) \bigr)_X \geq \frac{\nu^2}{2}|\partial_x w_i|_X^2 - \frac{\nu^2}{2}|\partial_x w_{i-1}|_X^2,
 \end{align*}
 \begin{align*}
 \ds (\xi_i, w_i - w_{i-1})_X \geq \int_\Omega \hat{K}^\delta(w_i)\, dx - \int_\Omega \hat{K}^\delta(w_{i-1})\, dx,
 \end{align*}
 \begin{align*}%\label{maTh01-021}
  \ds \bigl(g(w_i) , w_i -w_{i-1} \bigr)_X \geq \int_\Omega G(w_i)\, dx - \int_\Omega G(w_{i-1})\, dx - \frac{C_g}{2}|w_i - w_{i-1}|_X^2,
 \end{align*}
 and 
 \begin{align*}
 \ds M_u (u_i, w_i - w_{i-1})_X \leq \frac{1}{4\tau}|w_i - w_{i-1}|_X^2 + \tau M_u^2|u_i|_X^2.
 \end{align*}
In view of (\hyperlink{A5l}{A5}) and \eqref{maTh01-15}, we can see that: 
         \begin{align*}
          \displaystyle \frac{1}{2\tau}|w_i - w_{i-1}|_X^2 + \mathscr{F}^{(\varepsilon, \delta)}(w_i) - \mathscr{F}^{(\varepsilon, \delta)}(w_{i-1}) \leq \tau M_u^2|u_i|_X^2,\, i = 1, 2, 3, \ldots, n.
         \end{align*}
         
Thus we conclude Main Theorem \ref{mainTh01}.
\qed

\begin{rem}\label{rem10}
For every $\varepsilon, \delta \in (0, 1]$ and $i \in \{ 1, 2, 3, \ldots , n\}$, let us multiply the both sides of \eqref{MaTh01-01} by $\xi_i$. 
Then, by using the assumption (\hyperlink{A3l}{A3}) and \eqref{claim7-11}, one can observe that: 
\begin{align}\label{rem10-2}
 \ds \frac{1}{\tau}(w_i - w_{i-1}, \xi_i)_X + |\xi_i|^2_X - C_g(w_i, \xi_i)_X + (g(0), \xi_i)_X \leq M_u(u_i, \xi_i)_X, 
\end{align}
via
\begin{align*}
  (g(w_i), \xi_i)_X &\, = \bigl((g(w_i) + C_g w_i) - g(0), \xi_i \bigr)_X - (C_g w_i, \xi_i)_X + (g(0), \xi_i)_X 
 \\
 &\, \geq -C_g(w_i, \xi_i)_X + (g(0), \xi_i)_X. 
\end{align*} 

By using \eqref{rem10-2}, H\"{o}lder's and Young's inequalities, we obtain that: 
\begin{align}\label{rem10-3}
 \ds \frac{1}{4}|\xi_i|_X^2 \leq \frac{1}{\tau^2}|w_i - w_{i-1}|_X^2 + C_g^2|w_i|_X^2 + 2|g(0)|_X^2 + 2M_u|u_i|_X^2.
\end{align}
\end{rem}
\medskip

The above estimate will be used later.

\section{Proof of Main Theorem \ref{MainTh02}}
Let us fix $m \in \mathbb{N}$. 
On account of the energy inequality, we can see that: 
\begin{align*}
 \ds \frac{1}{2\tau}|w_i^m - w_{i-1}^m|_X^2 + \mathscr{F}^{(\varepsilon_m, \delta_m)}(w_i^m) - \mathscr{F}^{(\varepsilon_m, \delta_m)}(w_{i-1}^m) \leq \tau M_u^2|u_i^m|_X^2,\, i = 1, 2, 3, \ldots, n.
\end{align*}
Taking the sum of the above inequalities, for $i = 1, 2, 3, \ldots , n$, we have: 
\begin{align*}
 \ds \frac{1}{2\tau}&\, \sum_{i = 1}^l |w_i^m - w_{i - 1}^m|_X^2 + \mathscr{F}^{(\varepsilon_m, \delta_m)}(w_i^m) 
 \\ 
 &\, \leq \mathscr{F}^{(\varepsilon_m, \delta_m)}(w_0^m) + \tau M_u^2 \sum_{i = 1}^l |u_i^m|_X^2, l = 1, 2, 3, \ldots , n, \mbox{ and } m = 1, 2, 3, \ldots .
\end{align*}
Here, from (\hyperlink{A3l}{A3}) and \eqref{w.i}, we will find a positive constant $R_0$, independent of $m \in \mathbb{N}$, such that: 
\begin{align}\label{ken-10}
 |u^m|_{\mathbb{X}} \leq R_0,&\, |w_0^m|_{C(\overline{\Omega})} \leq R_0, |w_0^m|_Y \leq R_0, \nonumber 
 \\
 &\, \mbox{ and } |G(w_0^m)|_{C(\overline{\Omega})} \leq R_0, \mbox{ for } m = 1, 2, 3, \ldots .
\end{align}
Also, by using (\hyperlink{A2l}{A2}), it is estimated that: 
\begin{align}\label{ken-11}
 \ds |f^{\varepsilon_m}(\partial_x w_0^m)|_{L^1(\Omega)} &\, \leq \int_\Omega \left( \int_0^1 |(f^{\varepsilon_m})'(\varsigma \partial_x w_0^m)|\, d\varsigma \right)\, dx \nonumber
 \\
 &\,\leq  C_0 \left(2L + \frac{1}{2}\int_\Omega |\partial_x w_0^m|\, dx \right) \nonumber
 \\
 &\, \leq  C_0 \left(2L +  \sqrt{\frac{L}{2}} R_0 \right), \mbox{ for } m = 1, 2, 3, \ldots .
\end{align}
On account of \eqref{w.i}, \eqref{ken-10}, and \eqref{ken-11}, we compute that: 
\begin{align}\label{Main02-01}
 \ds \frac{1}{4T}|w_i^m|_X^2 &\, + \mathscr{F}^{(\varepsilon_m, \delta_m)}(w_i^m) \nonumber
 \\
 &\, \leq \frac{1}{2T}|w_0^m|_X^2 + \frac{1}{2\tau}\sum_{i=1}^l|w_i^m - w_{i-1}|_X^2 + \mathscr{F}^{(\varepsilon_m, \delta_m)}(w_i^m) \nonumber 
 \\
 &\, \leq \frac{R_0^2}{2T} +  \mathscr{F}^{(\varepsilon_m, \delta_m)}(w_0^m) +\tau M_u^2|u_i^m|_X^2 \nonumber
 \\
 &\, \leq \frac{R_0^2}{2T} +  C_0 \left(2L + \sqrt{\frac{L}{2}} R_0 \right) + \frac{\nu^2}{2}R_0^2 + \hat{K}_* + 2L R_0 + \tau M_u^2 R_0^2, \nonumber
 \\
 &\, \mbox{ for } i = 1, 2, 3, \ldots , n, \mbox{ and }  m = 1, 2, 3, \ldots .
\end{align}

As a consequence of \eqref{free.en}, \eqref{w.i} and \eqref{Main02-01}, one can observe that: 
\begin{description}
    \item[\textmd{$(\hypertarget{blacklozenge1}{\blacklozenge\,1})$}]the sequence $\{ w^m\}_{m \in \mathbb{N}} = \{ [w_1^m, \ldots , w_n^m]\}_{m \in \mathbb{N}}$ is bounded in $ \mathbb{Y} $.
\end{description}
Furthermore, by virtue of (\hyperlink{A3l}{A3}), \eqref{MaTh01-01}, \eqref{w.i}, \eqref{rem10-3}, and $(\hyperlink{blacklozenge1}{\blacklozenge\,1})$, it is derived that: 
\begin{description}
    \item[\textmd{$(\hypertarget{blacklozenge2}{\blacklozenge\,2})$}]the sequence $\{ \xi^m\}_{m \in \mathbb{N}} = \{ [\xi_1^m, \ldots , \xi_n^m]\}_{m \in \mathbb{N}} $ is bounded in $ \mathbb{X} $,
\end{description} 
\begin{description}
    \item[\textmd{$(\hypertarget{blacklozenge3}{\blacklozenge\,3})$}]the sequence $\{ \partial_x(\varpi^{(*, m)} + \nu^2 \partial_x w^m)\}_{m \in \mathbb{N}} = \{ [\partial_x(\varpi_1^{(*, m)} + \nu^2 \partial_x w_1^m), \ldots , \partial_x(\varpi_n^{(*, m)} + \nu^2 \partial_x w_n^m) ]\}_{m \in \mathbb{N}} $ is bounded in $ \mathbb{X} $, with a sequence $\{ \varpi_i^{(*, m)}\}_{m\in \mathbb{N}} \subset X$ satisfying $\varpi_i^{(*, m)} \in \partial f^{\varepsilon_m}(\partial_x w_i^m)$ in $\mathbb{R}$, a.e. in $\Omega$. 
\end{description}

Note that \eqref{step1-02} leads to: 
\begin{align}\label{Main02-1}
 \partial_x w_i^m &\, = (\mu_{\varepsilon_m}^\nu)^{-1}\bigl(\mu_{\varepsilon_m}^\nu(\partial_x w_i^m) \bigr) 
 \\ \nonumber
 &\, = (\mu_{\varepsilon_m}^\nu)^{-1}(\varpi_i^{(*, m)} + \nu^2 \partial_x w_i^m), \mbox{ for any } m \in \mathbb{N} \mbox{ and } i = 1, 2, 3, \ldots , n.
\end{align}
With Remark \ref{rem.f.lip} and \eqref{Main02-1} in mind, we apply the generalized chain rule in BV-theory \cite[Theorem 3.99]{MR1857292}, and we infer that:
\begin{align}\label{Main02-2}
 |\partial_x^2 w_i^m|_X &\, = \bigl|\bigl((\mu_{\varepsilon_m}^\nu)^{-1}\bigr)'(\varpi_i^{(*, m)} + \nu^2 \partial_x w_i^m)\partial_x(\varpi_i^{(*, m)} + \nu^2 \partial_x w_i^m) \bigr|_X \nonumber
 \\
 &\, \ds \leq \frac{1}{\nu^2}|\partial_x(\varpi_i^{(*, m)} + \nu^2 \partial_x w_i^m)|_X, \mbox{ for any } m \in \mathbb{N} \mbox{ and } i = 1, 2, 3, \ldots, n.
\end{align}
From $(\hyperlink{blacklozenge3}{\blacklozenge\,3})$ and \eqref{Main02-2}, it is deduced that: 
\begin{description}
    \item[\textmd{$(\hypertarget{blacklozenge4}{\blacklozenge\,4})$}]the sequence $\{ \partial_x^2 w^m\}_{m \in \mathbb{N}} = \{ [\partial_x^2 w_1^m, \ldots , \partial_x^2 w_n^m]\}_{m \in \mathbb{N}} $ is bounded in $ \mathbb{X} $, and hence $\{ w^m\}_{m \in \mathbb{N}}$ is bounded in $[H^2(\Omega)]^n$.
\end{description}
As a consequence of the one-dimensional compact embedding $H^2(\Omega) \subset C^1(\overline{\Omega})$ and $(\hyperlink{blacklozenge4}{\blacklozenge\,4})$, there exist $\{ m_k\}_{k \in \mathbb{N}}\subset \{ m \}$, $\tilde{w} = [\tilde{w}_1, \ldots, \tilde{w}_n] \in [H^2(\Omega)]^n$ such that: 
\begin{align}\label{Main02-03}
 w_i^{m_k} \to &\, \tilde{w}_i \mbox{ in } Y, \mbox{ in } C^1(\overline{\Omega}), \mbox{ and weakly in } H^2(\Omega), \nonumber
 \\
 &\, \mbox{ for any } i = 1, 2, 3, \ldots, n, \mbox{ as } k \to \infty.
\end{align}

Next, we verify that the limit $\tilde{w}$ is the solution to the state-system \hyperlink{(AC)$^{(\varepsilon, \delta)}$}{(AC)$^{(\varepsilon, \delta)}$}. 
Let us fix $k \in \mathbb{N}$.
From Remark \ref{MaTh.rem}, the solution $w^{m_k}$ admits the following variational inequality: 
\begin{align}\label{Main02-08}
 \displaystyle &\, \frac{1}{\tau}(w_i^{m_k} - w_{i-1}^{m_k}, w_i^{m_k} - z)_X + (g(w_i^{m_k}) - M_u u_i^{m_k}, w_i^{m_k} - z)_X \nonumber 
 \\
 &\,\quad  + \Psi^{(\varepsilon_{m_k}, \delta_{m_k})}(w_i^{m_k}) - \Psi^{(\varepsilon_{m_k}, \delta_{m_k})}(z) \leq 0, \nonumber
 \\
 &\, \qquad \mbox{ for any } z \in D(\Psi^{(\varepsilon_{m_k}, \delta_{m_k})}) \mbox{ and } i = 1, 2, 3, \ldots, n.
\end{align}
Here, we compute the limit of the both sides of \eqref{Main02-08}, as $k \to \infty$. 
Then, with (\hyperlink{A3l}{A3}), Key-Lemma and \eqref{Main02-03} in mind, one can see that: 
\begin{align*}%\label{Main02-09}
\displaystyle &\, \frac{1}{\tau}(\tilde{w}_i - \tilde{w}_{i-1}, \tilde{w}_i - z)_X + (g(\tilde{w}_i), \tilde{w}_i - z)_X + \Psi^{(\varepsilon, \delta)}(\tilde{w}_i)  \nonumber 
 \\
 \displaystyle &\, \leq \frac{1}{\tau}\lim_{k \to \infty}(w^{m_k}_i - w^{m_k}_{i-1}, w^{m_k}_i - z)_X + \lim_{k \to \infty}(g(w^{m_k}_i), w^{m_k}_i - z)_X + \varliminf_{k \to \infty}\Psi^{(\varepsilon_{m_k}, \delta_{m_k})}(w_i^{m_k}) \nonumber 
 \\
 &\, \leq  \Psi^{(\varepsilon, \delta)}(z) + (M_u u_i, \tilde{w}_i - z)_X , \mbox{ for any } z \in D(\Psi^{(\varepsilon, \delta)}) \mbox{ and } i = 1, 2, 3, \ldots, n.
\end{align*}
This implies that $\tilde{w}$ coincides with the solution $w$ to \hyperlink{(AC)$^{(\varepsilon, \delta)}$}{(AC)$^{(\varepsilon, \delta)}$}. 

Finally, we prove the convergence \eqref{mThConv00}. 
Let us fix $i \in \{ 1, 2, 3, \ldots , n\}$. 
Then, one can see from \eqref{Main02-03} that: 
\begin{align}\label{Sho1}
 \partial_x w_i^{m_k} \to \partial_x w_i \mbox{ in } C(\overline{\Omega}), \mbox{ as } k \to \infty.
\end{align}
In the light of (\hyperlink{A2l}{A2}) and \eqref{Sho1}, we compute that: 
\begin{align}\label{Sho2}
&\, |f^{\varepsilon_{m_k}}(\partial_x w_i^{m_k}) - f^\varepsilon(\partial_x w_i)| \nonumber
 \\
 &\,~ \leq |f^{\varepsilon_{m_k}}(\partial_x w_i^{m_k}) - f^{\varepsilon_{m_k}}(\partial_x w_i)| + |f^{\varepsilon_{m_k}}(\partial_x w_i)  - f^\varepsilon(\partial_x w_i)| \nonumber
 \\
 &\,~  \leq \ds \left| \left(\int_0^1 (f^{\varepsilon_{m_k}})'(\partial_x w_i + \varsigma(\partial_x w_i^{m_k} - \partial_x w_i))\, d\varsigma \right)(\partial_x w_i^{m_k} - \partial_x w_i)\right| \nonumber 
 \\
 &\, \qquad  + |f^{\varepsilon_{m_k}}(\partial_x w_i)  - f^\varepsilon(\partial_x w_i)| \nonumber 
 \\
  &\, ~ \leq \left( \int_0^1 C_0(1 + |\partial_x w_i + \varsigma(\partial_x w_i^{m_k} - \partial_x w_i)|)\, d\varsigma \right)|\partial_x w_i^{m_k} - \partial_x w_i| \nonumber 
  \\
  &\, \qquad  + |f^{\varepsilon_{m_k}}(\partial_x w_i)  - f^\varepsilon(\partial_x w_i)| \nonumber 
  \\
  &\, ~ \leq \ds C_0|\partial_x w_i^{m_k} - \partial_x w_i|_{C(\overline{\Omega})} \int_0^1 (1 + \varsigma|\partial_x w_i| + (1 - \varsigma)|\partial_x w_i^{m_k}|)\, d\varsigma \nonumber 
  \\
  &\, \qquad  + |f^{\varepsilon_{m_k}}(\partial_x w_i)  - f^\varepsilon(\partial_x w_i)| \nonumber 
  \\
  &\, ~ = \ds C_0|\partial_x w_i^{m_k} - \partial_x w_i|_{C(\overline{\Omega})}\left(1 + \frac{1}{2}|\partial_x w_i| + \frac{1}{2}|\partial_x w_i^{m_k}|\right) \nonumber 
  \\
  &\, \qquad  + |f^{\varepsilon_{m_k}}(\partial_x w_i)  - f^\varepsilon(\partial_x w_i)| \nonumber 
  \\
  &\, ~ \to 0 \mbox{ in the pointwise sense on } \overline{\Omega}, \mbox{ as } k \to \infty.
\end{align}
On account of (\hyperlink{A2l}{A2}) and \eqref{Sho1}, we can say that: 
\begin{align}\label{Sho3}
 \ds |f^{\varepsilon_{m_k}}(\partial_x w_i^{m_k})| &\, \leq \left|\int_0^1 (f^{\varepsilon_{m_k}})(\varsigma \partial_x w_i^{m_k})\, d\varsigma \right| \leq \int_0^1 C_0(1 + \varsigma|\partial_x w_i^{m_k}|)\, d\varsigma \nonumber 
 \\
 &\, \leq C_0 \left( 1 + \frac{1}{2} \sup_{k \in \mathbb{N}}|\partial_x w_i^{m_k}| \right) \mbox{ on } \overline{\Omega}, \mbox{ for any } k \in \mathbb{N}.
\end{align}
By \eqref{Sho2}, \eqref{Sho3}, and Lebesgue's dominated convergence theorem, we can infer that: 
\begin{align}\label{Sho4}
 f^{\varepsilon_{m_k}}(\partial_x w_i^{m_k}) \to f^\varepsilon(\partial_x w_i) \mbox{ in } X, \mbox{ as } k \to \infty. 
\end{align}

Now, taking into account \eqref{Main02-03}, \eqref{Sho2}, \eqref{Sho4}, and the uniqueness of the solution $w$, we conclude the convergences \eqref{mThConv} and \eqref{mThConv00}, with no use of subsequence. 
\qed

\section{Proof of Main Theorem \ref{mainTh02}}

In this section, we prove the third Main Theorem \ref{mainTh02}. 
Let us fix the constants $\varepsilon, \delta \in [0, 1]$, and fix the initial data $w_0 \in Y$ satisfying $\hat{K}^\delta(w_0) \in L^1(\Omega)$. 
Let us fix any forcing term $\bar{u} = [\bar{u}_1, \ldots , \bar{u}_n] \in \mathbb{X}$.
Then, invoking \eqref{J}, it is estimated that:
\begin{align}\label{mTh02-00}
    0 &~ \leq \underline{J}^{(\varepsilon, \delta)}:= \inf \mathcal{J}^{(\varepsilon, \delta)}(\mathbb{X}) \leq \overline{J}^{(\varepsilon, \delta)}:= \mathcal{J}^{(\varepsilon, \delta)}(\bar{u}) < \infty, \nonumber
    \\
    &\, \qquad \mbox{ for all } \varepsilon \in [0, 1] \mbox{ and } \delta \in [0, 1].
\end{align}
Also, for any $ \varepsilon \in [0, 1] $ and $\delta \in [0, 1]$, we denote by $\bar{w} = [\bar{w}_1, \ldots , \bar{w}_n]$ the solution to \hyperlink{(AC)$^{(\varepsilon, \delta)}$}{(AC)$^{(\varepsilon, \delta)}$} for the forcing term $ \bar{u}$ and initial data $w_0$.
\bigskip

Based on these, the Main Theorem \ref{mainTh02} is proved  as follows.

\paragraph{Proof of Main Theorem \ref{mainTh02} (III-A).}
From the estimate \eqref{mTh02-00}, we immediately find a sequence of forcing functions $\{u^{m} \}_{m \in \mathbb{N}} = \{[u_1^{m}, \ldots , u_n^m] \}_{m \in \mathbb{N}} \subset \mathbb{X}$, such that:
\begin{subequations}\label{mTh02-0102}
\begin{equation}\label{mTh02-01}
\mathcal{J}^{(\varepsilon, \delta)}(u^{m}) \downarrow \underline{J}^{(\varepsilon, \delta)},\ \mbox{as}\ m \to \infty,
\end{equation}
and
\begin{align}\label{mTh02-02}
    \sup_{m \in \mathbb{N}}& \left| \sqrt{\frac{M_u}{2}} u^{m} \right|_{\mathbb{X}}^{2} \leq \mathcal{J}^{(\varepsilon, \delta)}(\bar{u}) < \infty.
\end{align}
\end{subequations}
Also, the estimate \eqref{mTh02-02} enables us to take a subsequence of $\{u^{m} \}_{m \in \mathbb{N}} \subset \mathbb{X}$ (not relabeled), and to find a pair of functions $u^{*} = [u_1^*, \ldots , u_n^*] \in \mathbb{X}$, such that:
\begin{equation}\label{mTh02-03}
\begin{array}{c}
\sqrt{M_u} u^m \to \sqrt{M_u} u^{*} \mbox{ weakly in } \mathbb{X}, \mbox{ as } m \to \infty. 
\end{array}
\end{equation}

Let $w^* = [w_1^*, \ldots , w_n^*] \in \mathbb{X}$ be the solution to \hyperlink{(AC)$^{(\varepsilon, \delta)}$}{(AC)$^{(\varepsilon, \delta)}$} for the forcing term $u^{*}$ and initial data $w_0$. 
Also, for any $m \in \mathbb{N}$, let $w^m = [w_1^m, \ldots , w_n^m] \in \mathbb{X}$ be the solution to (AC)$^{(\varepsilon_m, \delta_m)}$ for the forcing term $u^{m}$ and the initial data $w_0$. Then, having in mind \eqref{w.i}, \eqref{mTh02-03}, we can apply Main Theorem \ref{MainTh02}, to see that:
\begin{equation}\label{mTh02-04}
w^m \to w^{*} \mbox{ in } \mathbb{Y}, \mbox{ in } [C^1(\overline{\Omega})]^n, \mbox{ and weakly in } [H^2(\Omega)]^n, \mbox{ as } m \to \infty.
\end{equation}
On account of \eqref{mTh02-01}, \eqref{mTh02-03}, and \eqref{mTh02-04}, it is computed that:
\begin{align*}
    \mathcal{J}^{(\varepsilon, \delta)}(u^{*}) & = \frac{1}{2} \bigl| [\sqrt{M_w} (w^{*}-w^\mathrm{ad}) \bigr|_{\mathbb{X}}^{2}  + \frac{1}{2} \bigl| \sqrt{M_u} u^{*} \bigr|_{\mathbb{X}}^{2}
    \\
                                            & \leq \frac{1}{2}\lim_{m \to \infty} \bigl| [\sqrt{M_w} (w^{m}-w^\mathrm{ad})] \bigr|_{\mathbb{X}}^{2} + \frac{1}{2} \varliminf_{m \to \infty} \bigl| \sqrt{M_u} u^{m} \bigr|_{\mathbb{X}}^{2}
    \\
    & = \lim_{m \to \infty}\mathcal{J}^{(\varepsilon, \delta)}(u^m) = \underline{J}^{(\varepsilon, \delta)} ~ (\leq \mathcal{J}^{(\varepsilon, \delta)}(u^{*})),
\end{align*}
and this leads to:    
\begin{equation*}
    \mathcal{J}^{(\varepsilon, \delta)}(u^{*}) = \min_{u \in \mathbb{X}} \mathcal{J}^{(\varepsilon, \delta)}(u).
\end{equation*}

Thus, we conclude the item (\hyperlink{III-A}{III-A}). 
\qed

\paragraph{Proof of Main Theorem \ref{mainTh02} (III-B).} 
    Let $\varepsilon, \delta \in [0, 1]$, $ \{ \varepsilon_m \}_{m \in \mathbb{N}} \subset (0, 1] $, $\{ \delta_m\}_{m \in \mathbb{N}} \subset (0, 1]$, and $ \{ w_0^m \}_{m \in \mathbb{N}} \subset Y $ be the sequences as in \eqref{w.i}. 
    Let $ \bar{w} = [\bar{w}_1, \ldots , \bar{w}_n] \in \mathbb{X} $ be the solution to \hyperlink{(AC)$^{(\varepsilon, \delta)}$}{(AC)$^{(\varepsilon, \delta)}$} for the forcing term $ \bar{u} $ and initial data $w_0$, and for any $m \in \mathbb{N}$, let $\bar{w}^m = [\bar{w}_1^m, \ldots , \bar{w}_n^m] \in \mathbb{X} $ be the solution to (AC)$^{(\varepsilon_m, \delta_m)}$ for the forcing term $\bar{u}$ and initial data $w_0^m$.
Then, applying (\hyperlink{Fact2}{Fact\,2}) and Main Theorem \ref{MainTh02} to these solutions, it is observed that: 
    \begin{subequations}\label{mTh02-09}
        \begin{equation}\label{mTh02-09a}
            \bar{w}^m \to \bar{w} \mbox{ in $ \mathbb{Y} $,} \mbox{ in } [C^1(\overline{\Omega})]^n, \mbox{ and weakly in } [H^2(\Omega)]^n, 
        \end{equation}
        and 
        \begin{align}\label{mTh02-09b}
            \ds \sup_{m \in \mathbb{N}}|w_0^m|_Y < \infty,\, \hat{K}^\delta(w_0) \leq \varliminf_{m \to \infty}\int_\Omega \hat{K}^{\delta_m}(w_0^m) \leq \hat{K}_* < \infty.
        \end{align}
    \end{subequations}
    and hence,
    \begin{equation}\label{mTh02-10}
        \overline{J}_\mathrm{sup} := \sup_{m \in \N} J^{(\varepsilon_m, \delta_m)}(\bar{u}) < \infty.
    \end{equation}

    Next, for any $ m \in \N $, let us denote by $ w^{(*, m)} = [w_1^{(*, m)}, \ldots , w_n^{(*, m)}] \in \mathbb{X}$ the solution to (AC)$^{(\varepsilon_m, \delta_m)}$ for the forcing term $ u^{(*, m)} = [u_1^{(*, m)}, \ldots , u_n^{(*, m)}] $ of the optimal control of (OP)$^{(\varepsilon_m, \delta_m)}$ and initial data $ w_{0}^m$. Then, in the light of \eqref{mTh02-00} and \eqref{mTh02-10}, it is observed that:
    \begin{equation*}
        0 \leq \frac{1}{2} \bigl| \sqrt{M_u} u^{(*, m)} \bigr|_{\mathbb{X}}^2 \leq \underline{J}^{(\varepsilon_m, \delta_m)} \leq \overline{J}_\mathrm{sup} < \infty, \mbox{ for any } m \in \mathbb{N}. 
    \end{equation*}
    Therefore, one can find a subsequence $ \{ m_k \}_{k \in \mathbb{N}} \subset \{m\} $, together with a limiting pair of functions $ u^{**} = [u_1^{**}, \ldots , u_n^{**}] \in \mathbb{X} $, such that:
\begin{equation}\label{mTh02-11}
\begin{array}{c}
\sqrt{M_u} u^{(*, m_k)} \to \sqrt{M_u} u^{**}\ \mbox{weakly in}\ \mathbb{X},\ \mbox{as}\ k \to \infty, 
\\[0.5ex]
\mbox{and as well as}\ M_u u^{(*, m_k)} \to M_u u^{**}\ \mbox{weakly in}\ \mathbb{X}, \ \mbox{as}\ k \to \infty.
\end{array}
\end{equation}

Now, let us denote by $w^{**} = [w_1^{**}, \ldots , w_n^{**}] \in \mathbb{X} $ the solution to \hyperlink{(AC)$^{(\varepsilon, \delta)}$}{(AC)$^{(\varepsilon, \delta)}$} for the forcing term $u^{**}$ and initial data $ w_0 $.
Then, applying Main Theorem \ref{MainTh02}, again, to the solutions $w^{**}$ and $w^{(*, m_k)}$, $k = 1, 2, 3, \dots$, one can observe that:
\begin{align}\label{mTh02-12}
    w^{(*, m_k)} & \to w^{**} \mbox{ in $ \mathbb{Y} $,} \mbox{ in } [C^1(\overline{\Omega})]^n, \mbox{ and weakly in } [H^2(\Omega)]^n, \mbox{ as $ k \to \infty $.}
\end{align}
As a consequence of \eqref{mTh02-09}, \eqref{mTh02-11}, and \eqref{mTh02-12}, it is verified that:
\begin{align*}
    \mathcal{J}^{(\varepsilon, \delta)} (u^{**}) &= \frac{1}{2} \bigl| \sqrt{M_w} (w^{**}-w^\mathrm{ad}) \bigr|_{\mathbb{X}}^{2} + \frac{1}{2} \bigl| \sqrt{M_u} u^{**}] \bigr|_{\mathbb{X}}^{2}
    \\
    & \leq \frac{1}{2} \lim_{k \to \infty} \bigl| \sqrt{M_w} (w^{(*, m_k)} -w^\mathrm{ad}) \bigr|_{\mathbb{X}}^{2} + \frac{1}{2}\varliminf_{k \to \infty} \bigl| \sqrt{M_u} u^{(*, m_k)} \bigr|_{\mathbb{X}}^{2}
    \\
    & = \varliminf_{k \to \infty}\mathcal{J}^{(\varepsilon_{m_k}, \delta_{m_k})}(u^{(*, m_k)}) \leq \lim_{k \to \infty}\mathcal{J}^{(\varepsilon_{m_k}, \delta_{m_k})}(\bar{u})
    \\
    & = \frac{1}{2}\lim_{k \to \infty} \bigl| \sqrt{M_w} (\bar{w}^{m_k}-w^\mathrm{ad})) \bigr|_{\mathbb{X}}^{2} + \frac{1}{2}\bigl| \sqrt{M_u} \bar{u} \bigr|_{\mathbb{X}}^{2}
    \\
    & = \mathcal{J}^{(\varepsilon, \delta)}(\bar{u}).
\end{align*}
Since the choice of $\bar{u} \in \mathbb{X}$ is arbitrary, we conclude that:
\begin{equation*}
    \mathcal{J}^{(\varepsilon, \delta)}(u^{**}) = \min_{u \in \mathbb{X}} \mathcal{J}^{(\varepsilon, \delta)}(u),
\end{equation*}
and complete the proof of Main Theorem \ref{mainTh02} (\hyperlink{III-B}{III-B}).
\qed

\section{Proof of Main Theorem \ref{mainTh03}}

Let $ \varepsilon \in (0, 1] $ and $\delta \in (0, 1]$ be fixed constants, and let $ w_0 \in Y $ be the initial data, satisfying $\hat{K}^\delta(w_0) \in L^1(\Omega)$. Let us take any forcing term $ u = [u_1, \ldots , u_n] \in \mathbb{X}$, and take the unique solution $ w = [w_1, \ldots , w_n] \in \mathbb{X}$ to the state-system \hyperlink{(AC)$^{(\varepsilon, \delta)}$}{(AC)$^{(\varepsilon, \delta)}$}. Also, let us take any constant $ \lambda \in (0, 1) $ and any function $ h = [h_1, \ldots , h_n] \in \mathbb{X}$, and consider another solution $ w^\lambda = [w_1^\lambda, \ldots , w_n^\lambda] \in \mathbb{X} $ to the state-system \hyperlink{(AC)$^{(\varepsilon, \delta)}$}{(AC)$^{(\varepsilon, \delta)}$} for the perturbed forcing term $ u + \lambda h = [u_1 + \lambda h_1, \ldots , u_n + \lambda h_n] $ and initial data $ w_0 $. On this basis, we consider a sequence of function $ \{\chi^\lambda \}_{\lambda \in (0, 1)} = \{  [\chi_1^\lambda, \ldots , \chi^\lambda_n ] \}_{\lambda \in (0, 1)} \subset \mathbb{X}$, defined as:
\begin{equation}\label{set00}
    \chi^\lambda = [\chi_1^\lambda, \ldots , \chi^\lambda_n] := \frac{w^\lambda - w}{\lambda} =   \left[ \frac{w_1^\lambda -w_1}{\lambda}, \ldots , \frac{w_n^\lambda -w_n}{\lambda} \right] \in \mathbb{X}, \mbox{ for $ \lambda \in (0, 1)$.}
\end{equation}
This sequence acts a key-role in the computation of G\^{a}teaux differential of the cost function $ \mathcal{J}^{(\varepsilon, \delta)} $, for $ \varepsilon \in (0, 1] $ and $\delta \in (0, 1]$. 
\begin{rem}\label{Rem.GD01}
    Note that for any  $ \lambda \in (0, 1) $, the function $ w^\lambda = [w_1^\lambda, \ldots , w_n^\lambda] \in \mathbb{X} $ fulfills the following variational forms:
\begin{align*}
    \frac{1}{\tau}( \chi^\lambda_i -&\, \chi_{i-1}^\lambda , \varphi)_{X} + \nu^2(\partial_x \chi^\lambda_i, \partial_x \varphi)_{X} + \int_{\Omega} \left( \int_{0}^{1} (f^{\varepsilon})''(\partial_x w_i +\varsigma \lambda \partial_x \chi^\lambda_i) \, d\varsigma \right) \partial_x \chi^\lambda_i \partial_x \varphi \, dx
    \\
    & +\int_{\Omega} \left(\int_{0}^{1} g'(w_i +\varsigma \lambda \chi^\lambda_i) \, d \varsigma \right) \chi^\lambda_i \varphi \, dx +\int_{\Omega}\left( \int_{0}^{1} (K^{\delta})'(w_i +\varsigma \lambda \chi^\lambda_i) \, d\varsigma \right) \chi^\lambda_i \varphi \, dx
    \\
    = & (M_u h_i, \varphi)_{X}, \mbox{ for any $ \varphi \in Y $ and $  i = 1, 2, 3,  \ldots, n$,} \mbox{ subject to $ \chi^\lambda_0 = 0 $ in $ X $.}
\end{align*}
    In fact, this variational form is obtained by taking the difference between respective two variational forms for $ w^\lambda = [w_1^\lambda, \ldots, w_n^\lambda]$ and $ w = [w_1, \ldots , w_n]$, as in Main Theorem \ref{mainTh01}, and by using the following linearization formulas: 

\begin{align*}
    & \frac{1}{\lambda} \bigl( (f^\varepsilon)'(\partial_x w_i^\lambda)- (f^\varepsilon)'(\partial_x w_i) \bigr) = \left( \int_{0}^{1}(f^\varepsilon)''(\partial_x w_i + \varsigma \lambda \partial_x \chi^\lambda_i) \, d\varsigma \right) \partial_x \chi^{\lambda}_i \mbox{ in $ X $,}
\end{align*}
\begin{align*}
    & \frac{1}{\lambda} \bigl( g(w_i^\lambda)-g(w_i) \bigr) = \left( \int_{0}^{1}g'(w_i + \varsigma \lambda \chi^\lambda_i) \, d\varsigma \right) \chi^{\lambda}_i \mbox{ in $ X $,}
\end{align*}
and
\begin{align*}
    & \frac{1}{\lambda} \bigl( K^\delta(w_i^\lambda)-K^\delta(w_i) \bigr) = \left( \int_{0}^{1}(K^\delta)'(w_i + \varsigma \lambda \chi^\lambda_i) \, d \varsigma \right) \chi^{\lambda}_i \mbox{ in $ X $, for any $i = 1, 2, 3,  \ldots, n$.}
\end{align*}
Incidentally, the above linearization formulas can be verified as consequences of the assumptions (\hyperlink{A1l}{A1})--(\hyperlink{A6l}{A6}) and the mean-value theorem (cf. \cite[Theorem 5 in p. 313]{lang1968analysisI}).
\end{rem}
\medskip

Now, we prepare the following three Lemmas, for the proof of Main Theorem \ref{mainTh03}.

    \begin{lem}\label{Lem.GD00}%\textbf{(time-discrete version of Gronwall's inequality.)}
        Let $\tau \in (0, \tau^*)$ be as in the assumption (\hyperlink{A5l}{A5}), and let $ c \geq 0 $ be a fixed constant such that:
        \begin{equation}\label{axLem03-00}
            0 < c \tau < \frac{1}{2}.
        \end{equation}
        Let $ \{A_i\}_{i = 0}^{n} \subset [0, \infty) $, $ \{ B_i \}_{i = 0}^{n} \subset [0, \infty) $ , and $ \{ C_i \}_{i = 1}^{n} \subset [0, \infty) $ be sequences such that:
        \begin{equation}\label{axLem03-01}
            \frac{1}{\tau} (A_i -A_{i -1} + \tau B_i) \leq c A_i +C_i, ~ i = 1, 2, 3, \ldots, n.
        \end{equation}
        Then, it is estimated that:
        \begin{align}\label{axLem03-02}
            A_i + \tau B_i  \leq 2^n \left( A_0 + \tau B_0 +\tau \sum_{j = 1}^{n} C_j  \right), ~ i = 1, 2, 3, \dots, n.
        \end{align}
    \end{lem}
    \paragraph{Proof.}{
        From the assumptions \eqref{axLem03-00} and \eqref{axLem03-01}, it is easily derived that:
        \begin{align*}
            \frac{1}{2} (A_i + \tau B_i) & \leq A_{i-1} + \tau C_i
            \\
                                         & \leq (A_{i-1} + \tau B_{i-1}) + \tau C_i, ~ i = 1, 2, 3, \dots, n.
        \end{align*}
        Next, we put: 
        \begin{align*}
         P_i := (A_i + \tau B_i), i = 0, 1, 2, \ldots, n.
        \end{align*}
        On this basis, we observe that:
        \begin{align*}
            P_1 ~& \leq 2 P_0 + 2\tau C_1,
            \\
            P_2 ~& \leq 2^2 P_0 + 2^2 \tau C_1 +2 \tau C_2,
            \\
            P_3 ~& \leq 2^3 P_0 + 2^3 \tau C_1 +2^2 \tau C_2 + 2 \tau C_3,
        \end{align*}
        and in general,
        \begin{align}\label{axLem03-03}
            P_i  \leq &\, 2^i P_0 + \tau (2^i C_1 + 2^{i-1}C_2 + \cdots + 2^2 C_{i-1} + 2C_i) \nonumber
\\
  \leq &\, 2^i(P_0 + \tau (C_1 + \cdots + C_i)),   \mbox{ for } i = 1, 2, 3, \dots, n.
        \end{align}
        The estimate \eqref{axLem03-02} is obtained as a straightforward consequence of \eqref{axLem03-03}. 
        \qed

\begin{lem}\label{Lem.GD01}
    Under the assumptions (\hyperlink{A1l}{A1})--(\hyperlink{A6l}{A6}), let us fix $ \varepsilon \in (0, 1] $ and $\delta \in (0, 1]$, and fix the initial data $w_0 \in Y$ satisfying $\hat{K}^\delta(w_0) \in L^1(\Omega)$. Then, for any $ u = [u_1, \ldots , u_n] \in \mathbb{X}$, the cost function $\mathcal{J}^{(\varepsilon, \delta)}$ admits the G\^{a}teaux derivative $ \left(\mathcal{J}^{(\varepsilon, \delta)} \right)'(u) \in \mathbb{X} ( = \mathbb{X}^*) $, such that: 
    \begin{align*}%\label{GD01}
        \bigl( (\mathcal{J}^{(\varepsilon, \delta)})'(u), & h \bigr)_{\mathbb{X}} = \bigl( M_w (w - w^\mathrm{ad}), \chi \bigr)_{\mathbb{X}} +\bigl( M_u u, h \bigr)_{\mathbb{X}}, 
        \\
        & \mbox{ for any $ h = [h_1, \ldots , h_n] \in \mathbb{X}$.} \nonumber
    \end{align*}
    In the context, $ w = [w_1, \ldots , w_n] \in \mathbb{X} $ is the solution to the state-system \hyperlink{(AC)$^{(\varepsilon, \delta)}$}{(AC)$^{(\varepsilon, \delta)}$} for the forcing term $ u $ and initial data $w_0$, and $ \chi = [\chi_1, \ldots , \chi_n] \in \mathbb{X}$ is a unique solution to the following linearization system:
    \begin{equation}\label{set01}
        \begin{cases}
            \displaystyle \frac{1}{\tau}(\chi_i - \chi_{i-1}) - \partial_x ((f^\varepsilon)'' (\partial_x w_i)\partial_x \chi_i + \nu^2 \partial_x \chi_i ) + g'(w_i)\chi_i + (K^\delta)'(w_i)\chi_i
            \\[0.5ex]
             \qquad = M_u h_i, \mbox{ in } \Omega,
            \\[1ex]
           \partial_x \chi_i (\pm L) = 0 , \mbox{ for any } i = 1, 2, 3,  \ldots, n,
            \\[1ex]
            \chi_0 = 0 \mbox{ in } X.
            \end{cases}
    \end{equation}
\end{lem}
\paragraph{Proof.}{
We prove this Lemma in accordance with the following two Steps.
    \begin{description}
        \item[\textbf{\boldmath First Step:}]The linearization system \eqref{set01} admits a unique solution.
        \item[\textbf{\boldmath Second Step:}]The cost function $\mathcal{J}^{(\varepsilon, \delta)}$ admits the G\^{a}teaux derivative $ \left(\mathcal{J}^{(\varepsilon, \delta)} \right)'(u) \in \mathbb{X} $, for any $u = [u_1, \ldots , u_n] \in \mathbb{X}$. 
    \end{description}

\paragraph{\boldmath\underline{Verification of First Step.}}{At first, we verify that the linearization system \eqref{set01} admits a solution $\chi = [\chi_1, \ldots , \chi_n]$. 
Let us fix $i \in \{ 1, 2, 3,  \ldots, n\}$ and let us define a functional $\mathcal{G}: X \longrightarrow(-\infty, \infty]$, by letting: 
\begin{align*}%\label{set01-01}
    z & \in X \mapsto \mathcal{G}(z) 
    \nonumber
    \\
    &:= \left\{
        \begin{array}{lll}
            \multicolumn{2}{l}{\displaystyle \frac{1}{2\tau}\int_\Omega |z - \chi_{i-1}|^2\, dx + \int_{\Omega}(f^\varepsilon)''(\partial_x w_i)|\partial_x z|^2\, dx + \frac{\nu^2}{2} \int_\Omega |\partial_x z|^2 \, dx}
            \\[2ex]
            \multicolumn{2}{l}{\displaystyle + \frac{1}{2}\int_{\Omega}g'(w_i)|z|^2\, dx +\frac{1}{2}\int_\Omega (K^\delta)'(w_i)|z|^2\, dx  - (M_u h_i, z)_X,}
            \\[2ex]
            & \mbox{ if $ z \in Y$,}
            \\[2ex]
            \infty, & \mbox{otherwise.}
        \end{array}
    \right. 
\end{align*}
From (\hyperlink{A5l}{A5}) and Remark \ref{g.prime}, it is observed that: 
\begin{align*}
 \ds \frac{1}{\tau} - g'(w_i) \geq \frac{1}{\tau} - [g']^- (w_i) \geq \frac{1}{\tau} - C_g > 0,
\end{align*}
and this enables us to say that the functional 
\begin{align}\label{g.conv}
 \ds z \in X \mapsto \frac{1}{2\tau}\int_\Omega |z - \chi_{i-1}|^2\, dx + \frac{1}{2} \int_\Omega g'(w_i)|z|^2\, dx
\end{align}
is strictly convex on $X$. 
From the assumptions (\hyperlink{A1l}{A1})--(\hyperlink{A5l}{A5}), Remark \ref{g.prime}, and \eqref{g.conv}, it is easily checked that $\mathcal{G}$ is a proper, l.s.c., and strictly convex function on $X$, such that:
\begin{align}\label{set01-02}
 \ds \mathcal{G}(z) &\, \geq \frac{1}{4\tau}|z|^2_X - \frac{1}{2\tau}|\chi_{i-1}|^2_X + \frac{\nu^2}{2}|\partial_x z|_X^2 \nonumber
 \\
 &\, \quad -\frac{C_g}{2}|z|_X^2 - M_u |h_i|_X|z|_X \nonumber
 \\
 &\, \geq \frac{1}{4\tau}|z|^2_X - \frac{1}{2\tau}|\chi_{i-1}|^2_X + \frac{\nu^2}{2}|\partial_x z|_X^2 \nonumber
 \\
 &\, \quad -\frac{C_g}{2}|z|^2_X - \frac{1}{16\tau}|z|_X^2 - 4\tau M_u^2|h_i|_X^2 \nonumber
 \\
 &\, \geq \frac{1}{8\tau}|z|^2_X - \frac{1}{2\tau}|\chi_{i-1}|^2_X + \frac{\nu^2}{2}|\partial_x z|_X^2 - 4\tau M_u^2|h_i|_X^2, \nonumber
 \\
 &\, \qquad \mbox{ for any } z \in Y.
\end{align} 
\eqref{set01-02} implies that $\mathcal{G}$ is coercive. 

Now, applying \cite[Proposition 1.2, Chapter II]{MR1727362}, we find a unique minimizer $\tilde{z} \in Y$ of $\mathcal{G}$, and hence, we can see that: 
\begin{align}\label{set01-03}
 \ds 0 &\, \leq \frac{1}{\lambda}\bigl( \mathcal{G}(\tilde{z} + \lambda \varphi) - \mathcal{G}(\tilde{z})\bigr) \nonumber
 \\
     &\, = \frac{1}{\tau}\int_\Omega \varphi (\tilde{z} - \chi_{i-1})\, dx + \frac{\lambda}{2\tau} \int_\Omega |\varphi|^2\, dx \nonumber
     \\
     &\, \quad + \int_\Omega (f^\varepsilon)''(\partial_x w_i)\partial_x \tilde{z}\partial_x \varphi\, dx + \frac{\lambda}{2}\int_\Omega (f^\varepsilon)''(\partial_x w_i)|\partial_x \varphi|^2\, dx \nonumber
     \\
     &\, \quad + \nu^2\int_\Omega \partial_x \tilde{z} \partial_x \varphi\, dx + \frac{\nu^2\lambda}{2}\int_\Omega |\partial_x \varphi|^2\, dx + \int_\Omega g'(w_i)\tilde{z}\varphi\, dx + \frac{\lambda}{2} \int_\Omega g'(w_i)|\varphi|^2\, dx \nonumber
     \\
     &\, \quad + \int_\Omega (K^\delta)'(w_i)\tilde{z}\varphi\, dx + \frac{\lambda}{2} \int_\Omega(K^\delta)'(w_i) |\varphi|^2\, dx - (M_u h_i, \varphi)_X, \nonumber
     \\
     &\, \qquad \mbox{ for any } \varphi \in Y, \mbox{ and } \lambda \in (0, 1).  
\end{align}
Taking $\lambda \downarrow 0$ in \eqref{set01-03}, it is inferred that: 
\begin{align}\label{set01-04}
 \ds 0 \leq &\, \frac{1}{\tau}(\tilde{z} - \chi_{i-1}, \varphi)_X + \bigl((f^\varepsilon)''(\partial_x w_i)\partial_x \tilde{z} + \nu^2 \partial_x \tilde{z}, \partial_x \varphi \bigr)_X \nonumber
 \\
            &\, + \bigl(g'(w_i)\tilde{z}, \varphi \bigr)_X + \bigl((K^\delta)'(w_i)\tilde{z}, \varphi \bigr)_X - (M_u h_i, \varphi)_X, \mbox{ for any } \varphi \in Y,
\end{align}
i.e.
\begin{align}\label{set01-05}
 &\, \frac{1}{\tau}(\tilde{z} - \chi_{i-1}, \varphi)_X + \bigl((f^\varepsilon)''(\partial_x w_i)\partial_x \tilde{z} + \nu^2 \partial_x \tilde{z}, \partial_x \varphi \bigr)_X \nonumber
 \\
 &\, \quad + \bigl(g'(w_i)\tilde{z}, \varphi \bigr)_X + \bigl((K^\delta)'(w_i)\tilde{z}, \varphi \bigr)_X = (M_u h_i, \varphi)_X, \mbox{ for any } \varphi \in Y.
\end{align}
In particular, taking any $\varepsilon_0 \in H_0^1(\Omega)$ and putting $\varphi = \varphi_0$ in \eqref{set01-05}, 
\begin{align*}%\label{set01-05-01}
 \ds &\, \left( M_u h_i - \frac{1}{\tau}(\tilde{z} - \chi_{i-1}) - g'(w_i)\tilde{z} - (K^\delta)'(w_i)\tilde{z}, \varphi_0 \right)_X \nonumber
 \\
 &\, = \int_\Omega \bigl( (f^\varepsilon)''(\partial_x w_i)\partial_x \tilde{z} + \nu^2 \partial_x \tilde{z} \bigr) \partial_x \varphi_0\, dx, \mbox{ for any } \varphi_0 \in H_0^1(\Omega),
\end{align*}
which implies: 
\begin{align}\label{set01-05-02}
\ds -\partial_x \bigl( (f^\varepsilon)''(\partial_x w_i)\partial_x \tilde{z} + \nu^2 \partial_x \tilde{z} \bigr)  = M_u h_i - \frac{1}{\tau}(\tilde{z} - \chi_{i-1}) - g'(w_i)\tilde{z} - (K^\delta)'(w_i)\tilde{z} \in X, \mbox{ in } \mathscr{D}'(\Omega).
\end{align} 
Additionally, having in mind Remark \ref{rem.f.lip}, (\hyperlink{A2l}{A2}), \eqref{set01-05-02}, and $w_i \in H^2(\Omega) \subset C^1(\overline{\Omega})$, we infer that 
\begin{subequations}\label{set01-06}
 \begin{align}\label{set01-06a}
 (f^\varepsilon)''(\partial_x w_i)\partial_x \tilde{z} + \nu^2 \partial_x \tilde{z} \in Y, \mbox{ and } |\partial_x \tilde{z}(\pm L)| \leq \frac{1}{\nu^2}\bigl| \bigl((f^\varepsilon)''(\partial_x w_i)\partial_x \tilde{z} + \nu^2 \partial_x \tilde{z} \bigr)(\pm L)\bigr| = 0
 \end{align}
i.e.
 \begin{align}\label{set01-06b}
 (f^\varepsilon)''(\partial_x w_i)\partial_x \tilde{z} + \nu^2 \partial_x \tilde{z} \in H^1_0(\Omega).
 \end{align}
\end{subequations}
As a consequence of \eqref{set01-05}--\eqref{set01-06}, we obtain the existence and uniqueness of the solution to the linearization system \eqref{set01}. 
}
\medskip

\paragraph{\boldmath\underline{Verification of Second Step.}}{
Let us fix any $ u = [u_1, \ldots , u_n] \in \mathbb{X} $, and take any $ \lambda \in (0, 1)$ and any $ h = [h_1, \ldots , h_n] \in \mathbb{X} $. 
Then, it is easily seen that: 
\begin{equation*}
    M_u(u +\lambda h) \to M_u u \mbox{ in $ \mathbb{X}$, as $ \lambda \downarrow 0 $,}
\end{equation*}
and also, as a consequence of Main Theorem \ref{MainTh02}, it is observed that:
\begin{align}\label{ken01}
    w^\lambda &\, = [w_1^\lambda, \ldots , w_n^\lambda] \to w = [w_1, \ldots , w_n] \nonumber
    \\
    &\, \mbox{ in $ \mathbb{Y} $, in $[C^1(\overline{\Omega})]^n$, and weakly in $[H^2(\Omega)]^n$, as $ \lambda \downarrow 0 $.}
\end{align}
Here, in the light of \eqref{set00}, \eqref{ken01}, (\hyperlink{A2l}{A2})--(\hyperlink{A4l}{A4}), we can find a constant $R_1 > 0$, independent of $\lambda \in (0, 1)$, such that: 
\begin{subequations}\label{ken02}
 \begin{align}\label{ken02a}
  |w|_{[C^1(\overline{\Omega})]^n \cap [H^2(\Omega)]^n} \vee \sup_{\lambda \in (0, 1)}|w^\lambda|_{[C^1(\overline{\Omega})]^n \cap [H^2(\Omega)]^n} \leq R_1,
 \end{align}
 and
 \begin{align}\label{ken02b}
 \max_{\substack{ 1 \leq i \leq n\\ \varsigma \in [0, 1]}}\Bigl\{ &\, |(f^\varepsilon)''(\partial_x w_i + \varsigma \lambda \partial_x \chi_i^\lambda)|_{C(\overline{\Omega})}, |g'(w_i + \varsigma \lambda \chi_i^\lambda)|_{C(\overline{\Omega})}, \nonumber
 \\
 &\, \quad  |(K^\delta)'(w_i + \varsigma \lambda \chi_i^\lambda)|_{L^\infty(\Omega)}\Bigr\} \leq R_1, 
 \nonumber 
 \\
 &\, \mbox{ for all } 0 < \lambda < 1.
 \end{align}
\end{subequations}
Also, taking a subsequence if necessary, we see from the assumptions (\hyperlink{A2l}{A2})--(\hyperlink{A4l}{A4}) that: 
    \begin{align}\label{set01-11}
      \ds  &\, \begin{cases}
            \ds \overline{f}_\lambda^\varepsilon := \left(\int_0^1 (f^\varepsilon)''(\partial_x w_i + \varsigma \lambda \partial_x \chi^\lambda_i) \, d\varsigma \right) \to (f^\varepsilon)''(\partial_x w_i),
            \\[2.5ex]
            \ds \overline{g}_\lambda := \left(\int_0^1 g'(w_i + \varsigma \lambda \chi_i^\lambda)\, d\varsigma \right) \to g'(w_i),
            \\[2.5ex]
           \ds \overline{K}^\delta_\lambda := \left(\int_0^1 (K^\delta)'( w_i + \varsigma \lambda \chi_i^\lambda)\, d\varsigma \right) \to (K^\delta)'(w_i), 
        \end{cases} 
        \\
        &\, \quad \quad \quad \quad \quad \quad \quad \mbox{ in the pointwise sense a.e. in $ \Omega $, as $ \lambda \to 0 $.} \nonumber
    \end{align}

In the meantime, it is easily seen that:
\begin{align}\label{GD02}
    \frac{1}{\lambda} \bigl(  \mathcal{J}^{(\varepsilon, \delta)} &( u +\lambda h) -\mathcal{J}^{(\varepsilon, \delta)}(u) \bigr)
    \nonumber
    \\
    = & \left( \frac{M_w}{2} (w^\lambda +w -2 w^\mathrm{ad}), \chi^\lambda \right)_{\hspace{-0.5ex}\mathbb{X}}  +\left( \frac{M_u}{2}(2u +\lambda h), h \right)_{\hspace{-0.5ex}\mathbb{X}}. 
\end{align}

Now, we fix $i \in \{ 1, 2, 3,  \ldots, n\}$. 
By using assumptions (\hyperlink{A1l}{A1})--(\hyperlink{A4l}{A4}), and Remark \ref{g.prime}, and by choosing $\varphi = \chi_i^\lambda$ in Remark \ref{Rem.GD01}, we can deduce that: 
\begin{subequations}\label{chi.est}
\begin{align}
\displaystyle \frac{1}{2\tau} (|\chi_i^\lambda|_X^2 - |\chi^\lambda_{i-1}|_{X}^2) + I_A + I_B \leq I_C + I_D, 
\end{align}
via 
\begin{align}
 \displaystyle \frac{1}{\tau} (\chi_i^\lambda - \chi_{i-1}^\lambda, \chi_i^\lambda)_X &\, \geq \frac{1}{\tau}|\chi_i^\lambda|_X^2 - \frac{1}{\tau} (\chi_{i-1}^\lambda, \chi_i^\lambda)_X \nonumber
 \\
 &\, \geq \frac{1}{2\tau} (|\chi_i^\lambda|_X^2 - |\chi_{i-1}^\lambda|_X^2),
\end{align}
\begin{align}
 \displaystyle I_A := &\, \left( -\partial_x  \left(\left( \int_{0}^{1}(f^\varepsilon)''(\partial_x w_i + \varsigma \lambda \partial_x \chi^\lambda_i) \, d\varsigma \right) \partial_x \chi^{\lambda}_i + \nu^2 \partial_x\chi_i^\lambda \right), \chi_i^\lambda \right)_X \nonumber
 \\
 &\, = \int_\Omega \left( \int_{0}^{1}(f^\varepsilon)''(\partial_x w_i + \varsigma \lambda \partial_x \chi^\lambda_i) \, d\varsigma \right) |\partial_x \chi^{\lambda}_i|^2 dx + (\nu^2 \partial_x \chi^\lambda_i, \partial_x \chi_i^\lambda)_X \nonumber
 \\
 &\, \geq \nu^2|\partial_x \chi_i^\lambda|_X^2,
\end{align}
\begin{align}
 \displaystyle I_B &\,: = \left( \left( \int_{0}^{1}(K^\delta)'(w_i + \varsigma \lambda \chi^\lambda_i) \, d \varsigma \right) \chi^{\lambda}_i, \chi_i^\lambda \right)_X \nonumber
 \\
 &\, = \int_\Omega \left( \int_{0}^{1}(K^\delta)'(w_i + \varsigma \lambda \chi^\lambda_i) \, d \varsigma \right) |\chi^{\lambda}_i|^2 dx \nonumber
 \\
 &\, \geq 0,
\end{align}
\begin{align}
\displaystyle  I_C &\, := - \left(\left( \int_{0}^{1}g'(w_i + \varsigma \lambda \chi^\lambda_i) \, d\varsigma \right) \chi^{\lambda}_i, \chi_i^\lambda \right)_X \nonumber
 \\
 &\, \leq - \int_\Omega \left( \int_{0}^{1}-[g']^{-}(w_i + \varsigma \lambda \chi^\lambda_i) \, d\varsigma \right) |\chi_i^\lambda|^2 dx \nonumber
 \\
 &\, \leq C_g|\chi_i^\lambda|_X^2,
\end{align}
and
\begin{align}
 I_D &\, := (M_u h_i, \chi_i^\lambda)_X \leq \frac{M_u^2}{2}|h_i|^2 + \frac{1}{2}|\chi_i^\lambda|_X^2.
\end{align}
\end{subequations}
Based on the \eqref{chi.est}, we compute that: 
\begin{align*}
 \ds \frac{1}{\tau}(|\chi_i^\lambda|_X^2 - |\chi_{i-1}^\lambda|_X^2 + 2\tau \nu^2 |\partial_x \chi_i^\lambda|_X^2) \leq (1 + 2C_g)|\chi_i^\lambda|_X^2 + M_u^2 |h_i|^2_X.
\end{align*} 

Now, let us set:
\begin{align*}%\label{risan.gron}
    \begin{cases}
        \ds A_i := |\chi_i^\lambda|_X^2, \quad B_i := 2\nu^2|\partial_x \chi_i^\lambda |_X^2 \mbox{ with } B_0 := B_1,
        \\[2ex]
        \ds c := 1 + 2C_g, \quad C_i := M_u^2 |h_i|_X^2.
    \end{cases} 
\end{align*}
Then, in the light of Lemma \ref{Lem.GD00}, one can say that:
\begin{description}
    \item[\textmd{$(\hypertarget{star1}{\star\,1})$}]the sequence $\{ \chi^\lambda\}_{\lambda \in (0, 1)} = \{ [\chi_1^\lambda, \ldots , \chi_n^\lambda\}_{\lambda \in (0, 1)} $ is bounded in $ \mathbb{Y} $, and compact in $\mathbb{X}$.
\end{description}

As consequences of $(\hyperlink{star1}{\star\,1})$, \eqref{ken02}, \eqref{set01-11}, (\hyperlink{A2l}{A2})--(\hyperlink{A4l}{A4}), and Lebesgue's dominated convergence theorem, one can find a sequence $ \{ \lambda_m \}_{m \in \mathbb{N}} \subset \{\lambda \}_{\lambda \in (0, 1)} $ and the function $\chi = [\chi_1, \ldots , \chi_n] \in \mathbb{X}$ such that:

\begin{subequations}\label{convs3-01}
    \begin{align}\label{convs3-00a}
        0 < |\lambda_m| < 1, \mbox{ and } \lambda_m \to 0, \mbox{ as $ m \to \infty $,}
    \end{align}
    \begin{align}\label{convs3-01a}
           & \chi_i^{\lambda_m} \to \chi_i \mbox{ in $X$, and weakly in $Y$,} \mbox{ as $ m \to \infty $, for any $i = 1, 2, 3,  \ldots, n$,}
    \end{align}
    \begin{align}\label{convs3-01b}
        \ds  \overline{f}_{\lambda_m}^\varepsilon \partial_x \varphi &\, \to (f^\varepsilon)''(\partial_x w_i)\partial_x \varphi \mbox{ in $ X $,} \nonumber 
        \\
        &\, \mbox{ for any $\varphi \in Y$, as $ m \to \infty $, for any $i = 1, 2, 3,  \ldots, n$,}
    \end{align}
    \begin{align}\label{convs3-01c}
        \ds &\, \overline{g}_{\lambda_m} \varphi \to g'(w_i) \varphi \mbox{ in $ X $,}
\mbox{ for any $\varphi \in Y$, as $ m \to \infty $, for any $i = 1, 2, 3,  \ldots, n$,}
    \end{align}
    and
    \begin{align}\label{convs3-01d}
        \ds  \overline{K}_{\lambda_m}^\delta \varphi &\, \to (K^\delta)'(w_i) \varphi \mbox{ in $ X $,} \nonumber
        \\
        &\, \mbox{ for any $\varphi \in Y$, as $ m \to \infty $, for any $i = 1, 2, 3,  \ldots, n$.}
    \end{align}
\end{subequations}

On account of \eqref{convs3-01}, Remark \ref{Rem.GD01}, and First Step, we will verify that $\chi$ is the unique solution to the  linearization system \eqref{set01}. 

Now, taking into account \eqref{ken01}, \eqref{GD02}, and \eqref{convs3-01}, and the uniqueness of the limit $\chi$, we can compute the directional derivative $D_h \mathcal{J}^{(\varepsilon, \delta)}(u) \in \mathbb{R}$, as follows:
\begin{align}\label{sho-10}
 \ds &\,  D_h\mathcal{J}^{(\varepsilon, \delta)}(u) := \lim_{\lambda \to 0} \frac{1}{\lambda} \bigl(  \mathcal{J}^{(\varepsilon, \delta)} ( u +\lambda h) -\mathcal{J}^{(\varepsilon, \delta)}(u) \bigr) = ( M_w (w - w^\mathrm{ad}), \chi )_{\mathbb{X}}  + ( M_u u , h )_{\mathbb{X}}, \nonumber
 \\
 &\, \quad \mbox{ for any } u = [u_1, \ldots , u_n] \in \mathbb{X}, \mbox{ and any direction } h = [h_1, \ldots, h_n] \in \mathbb{X}.
\end{align}
Moreover, with Riesz's theorem in mind, we deduce the existence of the G\^{a}teaux derivative $ \left(\mathcal{J}^{(\varepsilon, \delta)} \right)'(u) \in \mathbb{X}^* (= \mathbb{X})$ at $u = [u_1, \ldots , u_n] \in \mathbb{X}$, i.e.
\begin{align}\label{sho-11}
 \left(\left(\mathcal{J}^{(\varepsilon, \delta)} \right)'(u), h \right)_{\hspace{-0.5ex}\mathbb{X}} &\, = D_h\mathcal{J}^{(\varepsilon, \delta)}(u), \nonumber 
 \\
 &\, \mbox{ for every } u = [u_1, \ldots , u_n], h = [h_1, \ldots, h_n] \in \mathbb{X}.
\end{align} 
On account of \eqref{sho-10} and \eqref{sho-11}, we conclude this Lemma. 
\hfill \qed
}
}

\begin{lem}\label{Lem.GD02}
    Under the assumptions (\hyperlink{A1l}{A1})--(\hyperlink{A6l}{A6}), let us fix $\varepsilon \in (0, 1]$ and $\delta \in (0, 1]$, and fix the initial data $w_0 \in Y$ satisfying $\hat{K}^\delta(w_0) \in L^1(\Omega)$. Also, let $ u^* = [u_1^*, \ldots , u_n^*] \in \mathbb{X} $ be an optimal control of the problem \hyperlink{(OP)$^{(\varepsilon, \delta)}$}{(OP)$^{(\varepsilon, \delta)}$}, and let $ w^* =  [w_1^*, \ldots , w_n^*]$ be the solution to the state-system \hyperlink{(AC)$^{(\varepsilon, \delta)}$}{(AC)$^{(\varepsilon, \delta)}$} for the forcing term $ u^* $ and initial data $ w_0 $. 
Then, we can get the following equation:
    \begin{align*}
        \bigl( p , h \bigr)_{\mathbb{X}} = \bigl( v, \chi \bigr)_{\mathbb{X}}, \mbox{ for all $ h = [h_1, \ldots , h_n], v = [v_1, \ldots , v_n] \in \mathbb{X} $.}
    \end{align*}
    In the context, $\chi = [\chi_1, \ldots , \chi_n] \in \mathbb{X}$ and $p = [p_1, \ldots , p_n] \in \mathbb{X}$ are unique solutions for the following linearization system:
    \begin{equation}\label{set02}
        \begin{cases}
            \displaystyle \frac{1}{\tau}(\chi_i - \chi_{i-1}) - \partial_x ((f^\varepsilon)'' (\partial_x w_i^*)\partial_x \chi_i + \nu^2 \partial_x \chi_i ) + g'(w_i^*)\chi_i + (K^\delta)'(w_i^*)\chi_i
            \\[0.5ex]
             \qquad =  h_i, \mbox{ in } \Omega,
            \\[1ex]
           \partial_x \chi_i (\pm L) = 0 , \mbox{ for any } i = 1, 2, 3,  \ldots, n,
            \\[1ex]
            \chi_0 = 0 \mbox{ in } X;
            \end{cases}
    \end{equation}
     and the adjoint system:
    \begin{equation}\label{set03}
        \begin{cases}
            \displaystyle \frac{1}{\tau}(p_i - p_{i+1}) - \partial_x ((f^\varepsilon)'' (\partial_x w_i^*)\partial_x p_i + \nu^2 \partial_x p_i ) + g'(w_i^*)p_i + (K^\delta)'(w_i^*)p_i
            \\[0.5ex]
            \qquad =  v_i, \mbox{ in } \Omega,
            \\[1ex]
           \partial_x p_i (\pm L) = 0 , \mbox{ for any } i = n, \ldots , 3, 2, 1,
            \\[1ex]
            p_{n+1} = 0 \mbox{ in } X;
            \end{cases}
    \end{equation}
    respectively.
\end{lem}

\paragraph{Proof.}{
Let us fix arbitrary functions $ h = [h_1, \ldots , h_n],  v = [v_1, \ldots , v_n] \in \mathbb{X} $. 

First, we verify \eqref{set02} and \eqref{set03} admit  unique solutions $\chi = [\chi_1, \ldots , \chi_n] $ and $p = [p_1, \ldots , p_n]$, respectively. 
The existence and uniqueness of the solution $\chi$ to \eqref{set02} will be verified, immediately, by applying the same argument as in the First Step in the proof of Lemma \ref{Lem.GD01}, to the special case when $M_u = 1$. 

Now, we prove the existence and uniqueness for the system \eqref{set03}. 
Let us fix $i \in \{ 1, 2, 3, \ldots, n\}$ and we define a functional $\mathcal{L} : X \longrightarrow (-\infty, \infty]$, by letting: 
\begin{align*}%\label{set03-01}
    z & \in X \mapsto \mathcal{L}(z) 
    \nonumber
    \\
    &:= \left\{
        \begin{array}{lll}
            \multicolumn{2}{l}{\displaystyle \frac{1}{2\tau}\int_\Omega |z - p_{i+1}|^2\, dx + \int_{\Omega}(f^\varepsilon)''(\partial_x w_i^*)|\partial_x z|^2\, dx + \frac{\nu^2}{2} \int_\Omega |\partial_x z|^2 \, dx}
            \\[2ex]
            \multicolumn{2}{l}{\displaystyle + \frac{1}{2}\int_{\Omega}g'(w_i^*)|z|^2\, dx +\frac{1}{2}\int_\Omega (K^\delta)'(w_i^*)|z|^2\, dx  - (u_i, z)_X,}
            \\[2ex]
            & \mbox{ if $ z \in Y$,}
            \\[2ex]
            \infty, & \mbox{otherwise.}
        \end{array}
    \right. 
\end{align*}
Then, by using the assumptions (\hyperlink{A1l}{A1})--(\hyperlink{A5l}{A5}), and Remark \ref{g.prime}, and by applying similar arguments as in \eqref{g.conv}--\eqref{set01-06}, we can see the existence of the solution to the adjoint system \eqref{set03}. 

Next, having in mind (\hyperlink{A2l}{A2})--(\hyperlink{A4l}{A4}), \eqref{set02}, \eqref{set03}, and Main Theorem \ref{mainTh01}, it is deduced that: 
\begin{subequations}\label{GD02-01}
\begin{align}\label{GD02-01a}
\ds \frac{1}{\tau}&\, \sum_{i=1}^n (\chi_i - \chi_{i-1}, p_i)_X \nonumber
\\
&\, = \frac{1}{\tau}\bigl( (\chi_1 - \chi_0, p_1)_X + \cdots + (\chi_n - \chi_{n-1}, p_n)_X - (\chi_n, p_{n+1})_X \bigr) \nonumber
\\
&\, = \frac{1}{\tau}\bigl( (p_1 - p_2, \chi_1)_X + \cdots + (p_{n-1} - p_n, \chi_{n-1})_X + (p_n - p_{n+1}, \chi_n)_X \bigr) \nonumber
\\
&\, = \frac{1}{\tau}\sum_{i=1}^n (p_i - p_{i+1}, \chi_i)_X,
\end{align}
\begin{align}\label{GD02-01b}
 \bigl( (f^\varepsilon)''(\partial_x w_i^*)\partial_x \chi_i, \partial_x p_i \bigr)_X = \bigl( (f^\varepsilon)''(\partial_x w_i^*)\partial_x p_i, \partial_x \chi_i \bigr)_X,
\end{align}
\begin{align}\label{GD02-01c}
 (\nu^2 \partial_x \chi_i, \partial_x p_i)_X = (\nu^2 \partial_x p_i, \partial_x \chi_i)_X,
\end{align}
\begin{align}\label{GD02-01d}
 (g'(w_i^*)\chi_i, p_i)_X = (g'(w_i^*)p_i, \chi_i)_X, 
\end{align}
and
\begin{align}\label{GD02-01e}
 ((K^\delta)'(w_i^*)\chi_i, p_i)_X = ((K^\delta)'(w_i^*)p_i, \chi_i)_X, \mbox{ for any } i = 1, 2, 3,  \ldots, n.
\end{align}
\end{subequations}
Here, invoking \eqref{set02}, \eqref{set03}, and \eqref{GD02-01}, we compute that: 
\begin{align*}
 \ds \bigl( p ,&\, h \bigr)_{\mathbb{X}} = \sum_{i=1}^{n} (p_i, h_i)_X = \sum_{i=1}^{n} \langle h_i, p_i \rangle_Y
 \\
 &\, = \frac{1}{\tau}\sum_{i=1}^n (\chi_i - \chi_{i-1}, p_i)_X + \sum_{i=1}^n \bigl( (f^\varepsilon)''(\partial_x w_i^*)\partial_x \chi_i, \partial_x p_i \bigr)_X
 \\
 &\, \quad + \sum_{i=1}^n (\nu^2 \partial_x \chi_i, \partial_x p_i)_X + \sum_{i=1}^n (g'(w_i^*)\chi_i, p_i)_X + \sum_{i=1}^n ((K^\delta)'(w_i^*)\chi_i, p_i)_X
 \\
 &\, = \frac{1}{\tau}\sum_{i=1}^n (p_i - p_{i+1}, \chi_i)_X + \sum_{i=1}^n \bigl( (f^\varepsilon)''(\partial_x w_i^*)\partial_x p_i, \partial_x \chi_i \bigr)_X
 \\
 &\, \quad + \sum_{i=1}^n (\nu^2 \partial_x p_i, \partial_x \chi_i)_X + \sum_{i=1}^n (g'(w_i^*)p_i, \chi_i)_X + \sum_{i=1}^n ((K^\delta)'(w_i^*)p_i, \chi_i)_X
 \\
 &\, = \sum_{i=1}^{n} ( v_i, \chi_i )_X = (v, \chi)_{\mathbb{X}}.
\end{align*}
Thus, we conclude this Lemma.
\qed

}

\bigskip

Now, we are on the stage to prove the Main Theorem \ref{mainTh03}.

\paragraph{Proof of Main Theorem \ref{mainTh03}.}{
    Let $ u^* = [u_1^*, \ldots , u_n^*] \in \mathbb{X} $ be the optimal control of \hyperlink{(OP)$^{(\varepsilon, \delta)}$}{(OP)$^{(\varepsilon, \delta)}$}, with the solution $ w^* = [w_1^*, \ldots , w_n^*] \in \mathbb{X} $ to the state-system \hyperlink{(AC)$^{(\varepsilon, \delta)}$}{(AC)$^{(\varepsilon, \delta)}$} for the forcing term $ u^* $ and initial data $ w_0 \in Y$ satisfying $\hat{K}^\delta(w_0) \in L^1(\Omega)$. 
    Also, let us fix $h = [h_1, \ldots , h_n] \in \mathbb{X}$ and let $\chi^* = [\chi_1^*, \ldots , \chi_n^*] \in \mathbb{X}$ be the solution to the linearization system \eqref{set02} for the forcing term $M_u h$, and let $p^* = [p_1^*, \ldots , p_n^*] \in \mathbb{X}$ be the solution to the adjoint system \eqref{set03} for the forcing term $M_w (w^* - w^{\mathrm{ad}})$. 
Then, by applying Lemma \ref{Lem.GD01}, and by applying Lemma \ref{Lem.GD02} to the case when $h = M_u h$ and $v = M_w (w^* -w^{\mathrm{ad}})$, we compute that:
    \begin{align*}
        0 = &~ \left((\mathcal{J}^{(\varepsilon, \delta)})'(u^*), h \right)_{\mathbb{X}} 
        \\[1ex]
        = &~ \lim_{\lambda \to 0} \frac{1}{\lambda} \bigl( \mathcal{J}^{(\varepsilon, \delta)}(u^* +\lambda h) -\mathcal{J}^{(\varepsilon, \delta)}(u^*) \bigr) 
        \\[1ex]
        = &~ \bigl( M_w (w^* -w^{\mathrm{ad}}), \chi \bigr)_{\mathbb{X}} + \bigl( M_u u^*, h \bigr)_{\mathbb{X}}
        \\[1ex] 
        = &~ \bigl( p^*, M_u h \bigr)_{\mathbb{X}} + \bigl( M_u u^*, h \bigr)_{\mathbb{X}}
        \\[1ex]
        = &~ \bigl(M_u(p^* + u^*), h \bigr)_{\mathbb{X}}, \mbox{ for any } h = [h_1, \ldots , h_n] \in \mathbb{X}.
    \end{align*}

    Thus, we conclude  Main Theorem \ref{mainTh03}. \qed
}

\section{Proof of Main Theorem \ref{mainTh04}}

At first, we prepare the notations used throughout this Section. 
For every $\varepsilon \in (0, 1]$ and $\delta \in (0, 1]$, we denote by $u^{(*, (\varepsilon, \delta))} = [u_1^{(*, (\varepsilon, \delta))}, \ldots , u_n^{(*, (\varepsilon, \delta))}] \in \mathbb{X}$ the optimal control of the approximating problem \hyperlink{(OP)$^{(\varepsilon, \delta)}$}{(OP)$^{(\varepsilon, \delta)}$}, together with the solution $w^{(*, (\varepsilon, \delta))} = [w_1^{(*, (\varepsilon, \delta))}, \ldots , w_n^{(*, (\varepsilon, \delta))}] $ $\in \mathbb{X}$ to the state-system \hyperlink{(AC)$^{(\varepsilon, \delta)}$}{(AC)$^{(\varepsilon, \delta)}$} for the forcing term $u^{(*, (\varepsilon, \delta))}$ and initial data $w_0 \in Y$, satisfying $|w_0| \leq 1$, a.e. in $\Omega$. 
Note that: 
\begin{align*}
 \hat{K}^{\tilde{\delta}}(w_0) &\, \equiv 0 \mbox{ on } \overline{\Omega}, \mbox{ and hence } 
 \\
 &\, |\hat{K}^{\tilde{\delta}}(w_0)|_{L^1(\Omega)} = 0, \mbox{ for all } \tilde{\delta} \in [0, 1].
\end{align*}
In addition, for every $\varepsilon \in (0, 1]$ and $\delta \in (0, 1]$, we denote by $p^{(*, (\varepsilon, \delta))} = [p_1^{(*, (\varepsilon, \delta))}, \ldots , p_n^{(*, (\varepsilon, \delta))}] $ $\in \mathbb{X}$ the solution to the adjoint system \eqref{set03} for the forcing term $M_w (w^{(*, (\varepsilon, \delta))} - w^{\mathrm{ad}})$. 

Now, by Main Theorem \ref{mainTh02} (III-B), we find an optimal control $u^\circ = [u_1^\circ, \ldots , u_n^\circ] \in \mathbb{X}$ of (OP)$^{(0, 0)}$, with a zero-convergent sequences $\{ \varepsilon_m \}_{m \in \mathbb{N}} \subset (0, 1]$ and $\{ \delta_m\}_{m \in \mathbb{N}} \subset (0, 1]$, such that: 
\begin{subequations}\label{main04-01}
\begin{align}\label{main04-01-01}
 u^{(*, m)} := u^{(*, (\varepsilon_m, \delta_m))} \to u^\circ \mbox{ weakly in } \mathbb{X}, \mbox{ as } m \to \infty.
\end{align}
Let $w^\circ = [w_1^\circ, \ldots , w_n^\circ] \in \mathbb{X}$ be the solution to (AC)$^{(0, 0)}$ for the forcing term $u^\circ$ and initial data $w_0$. 
Then, having in mind (\hyperlink{A3l}{A3}) and Main Theorem \ref{MainTh02}, we can find subsequences of $\{ \varepsilon_m\}_{m \in \mathbb{N}}$ and $\{ \delta_m\}_{m \in \mathbb{N}}$ (not relabeled) such that: 
\begin{align}\label{main04-01-02}
 w^{(*, m)}:= w^{(*, (\varepsilon_m, \delta_m))} &\, \to w^\circ \mbox{ in } \mathbb{Y}, \mbox{ in } [C^1(\overline{\Omega})]^n, \nonumber
 \\
 &\, \mbox{ and weakly in } [H^2(\Omega)]^n, \mbox{ as } m \to \infty,
\end{align}
\begin{align}\label{main04-01-03}
 \partial_x w^{(*, m)} &\, := \partial_x w^{(*, (\varepsilon_m, \delta_m))} \to \partial_x w^\circ \mbox{ in } [C(\overline{\Omega})]^n, \mbox{ as } m \to \infty,
\end{align}
and 
\begin{align}\label{main04-01-04}
 g'(w_i^{(*, m)}) &\, := g'(w_i^{(*, (\varepsilon_m, \delta_m))})  \to g'(w_i^\circ) \mbox{ in } C(\overline{\Omega}), \nonumber
                  \\
                  &\, \mbox{ for any } i = 1, 2, 3, \ldots, n, \mbox{ as } m \to \infty. 
\end{align}
\end{subequations}
Here, in the light of (\hyperlink{A3l}{A3}) and \eqref{main04-01}, we can find a constant $R_2 > 0$, independent of $m \in \mathbb{N}$, such that:
\begin{align}\label{sho01}
 \max_{1 \leq i \leq n}|g'(w_i^{(*, m)})|_{C(\overline{\Omega})} \leq R_2, \mbox{ for all } m \in \mathbb{N}.
\end{align}

Next, let us fix $m \in \mathbb{N}$. 
Then, by \eqref{set03}, $p^{(*, m)} := p^{(*, (\varepsilon_m, \delta_m))}$ satisfies the following equation: 
\begin{align}\label{main04-02}
   \ds  &\, \frac{1}{\tau} ( p^{(*, m)}_i - p^{(*, m)}_{i + 1},  \varphi )_X + \bigl( (f^{\varepsilon_m})''(\partial_x w_i^{(*, m)})\partial_x p_i^{(*, m)}, \partial_x \varphi \bigr)_X 
    \nonumber
    \\[0.5ex]
    &\,~ + \nu^2(\partial_x p^{(*, m)}_i, \partial_x \varphi )_X  +\bigl( g'(w_i^{(*, m)}) p^{(*, m)}_i, \varphi \bigr)_{X} \nonumber
    \\[0.5ex]
    &\,~ +\bigl( (K^\delta)'(w_i^{(*, m)})p_i^{(*, m)} , \varphi \bigr)_{X} = \bigl( M_w (w^{(*, m)}_i - w_i^{\mbox{\scriptsize ad}}), \varphi \bigr)_{X}, \nonumber
    \\[0.5ex]
    &\,~ \mbox{ for any $ \varphi \in Y $, and $ i = n, \ldots, 3, 2, 1 $.}
\end{align}
Let us fix $i \in \{1, 2, 3, \ldots, n\}$. 
Then, by applying the assumptions (\hyperlink{A1l}{A1})--(\hyperlink{A4l}{A4}), and Remark \ref{g.prime}, and by choosing  $\varphi = p_i^{(*, m)}$ in \eqref{main04-02}, we can deduce that:
\begin{subequations}\label{main04-03}
\begin{align}\label{main04-03-01}
 \ds \frac{1}{2\tau}(|p_i^{(*, m)}|_X^2 - |p_{i+1}^{(*, m)}|_X^2) + \nu^2|\partial_x p_i^{(*, m)}|_X^2 + L_A + L_B \leq L_C + L_D, 
\end{align}
with
\begin{align}\label{main04-03-02}
 \ds L_A  := &\, \bigl( (f^{\varepsilon_m})''(\partial_x w_i^{(*, m)})\partial_x p_i^{(*, m)}, \partial_x p_i^{(*, m)} \bigr)_X \nonumber 
 \\
      = &\, \int_\Omega (f^{\varepsilon_m})''(\partial_x w_i^{(*, m)})|\partial_x p_i^{(*, m)}|^2\, dx \geq 0,
\end{align}
\begin{align}\label{main04-03-03}
 \ds L_B  := &\, \bigl( (K^{\delta_m})'(w_i^{(*, m)})p_i^{(*, m)} , p_i^{(*, m)} \bigr)_{X} \nonumber
 \\
  = &\, \int_\Omega (K^{\delta_m})'(w_i^{(*, m)}) | p_i^{(*, m)}|^2\, dx \geq 0,
\end{align}
\begin{align}\label{main04-03-04}
 \ds L_C  := &\, - \bigl( g'(w_i^{(*, m)}) p^{(*, m)}_i, p^{(*, m)}_i \bigr)_{X} \nonumber
 \\
 &\, \leq  \int_{\Omega} |[g']^{-}|_{L^\infty(\mathbb{R})}|p^{(*, m)}_i|^2\, dx \nonumber
 \\
 &\, \leq C_g|p^{(*, m)}_i|_X^2,
\end{align}
and 
\begin{align}\label{main04-03-05}
\ds L_D  := &\, \bigl( M_w (w^{(*, m)}_i - w_i^{\mbox{\scriptsize ad}}), p^{(*, m)}_i \bigr)_{X} \nonumber
\\
         \leq &\, \frac{1}{2}| p^{(*, m)}_i|_X^2 + \frac{M_w^2}{2}| w^{(*, m)}_i - w_i^{\mbox{\scriptsize ad}}|_X^2.
\end{align}
\end{subequations}
Based on the \eqref{main04-03}, we compute that: 
\begin{align*}%\label{main04-04}
 \ds \frac{1}{\tau}\bigl((|p_i^{(*, m)}|_X^2 &\, - |p_{i+1}^{(*, m)}|_X^2) + 2\tau \nu^2|\partial_x p_i^{(*, m)}|_X^2 \bigr) \nonumber
 \\
 &\, \leq (1 + 2C_g)| p^{(*, m)}_i|_X^2 + M_w^2| w^{(*, m)}_i - w_i^{\mbox{\scriptsize ad}}|_X^2.
\end{align*}
Here, let us apply Lemma \ref{Lem.GD00} to the case when:
\begin{align*}%\label{main04-05}
    \begin{cases}
        \ds A_i := |p_{n-i+1}^{(*, m)}|_X^2, \quad B_i := 2\nu^2 |\partial_x p_{n-i+1}^{(*, m)}|_X^2 \mbox{ with } B_0 := B_1,
        \\[2ex]
        \ds c := 1 + 2C_g, \quad C_i := M_w^2| w^{(*, m)}_{n-i+1} - w_{n-i+1}^{\mbox{\scriptsize ad}}|_X^2.
    \end{cases} 
\end{align*}
Then, in the light of Lemma \ref{Lem.GD00} and \eqref{main04-01-02}, one can say that:
\begin{description}
    \item[\textmd{$(\hypertarget{star2}{\star\,2})$}]the sequence $\{ p^{(*, m)} \}_{m \in \mathbb{N}} = \{ [p_1^{(*, m)}, \ldots , p_n^{(*, m)}] \}_{m \in \mathbb{N}}$ is bounded in $ \mathbb{Y} $, and compact in $\mathbb{X}$.
\end{description}

Now, let us fix $m \in \mathbb{N}$ and $i \in \{ 1, 2, 3, \ldots, n\}$, and define a bounded and linear functional $\zeta_i^m \in Y^*$ on $Y$, by putting:
\begin{align}\label{main04-06}
 \ds \langle \zeta_i^m, \varphi \rangle_Y := &\, \bigl( (f^{\varepsilon_m})''(\partial_x w_i^{(*, m)})\partial_x p_i^{(*, m)}, \partial_x \varphi \bigr)_X \nonumber
 \\
 &\, \quad + \bigl( (K^{\delta_m})'(w_i^{(*, m)})p_i^{(*, m)} , \varphi \bigr)_{X}, \mbox{ for any } \varphi \in Y.
\end{align}
Then, on account of \eqref{sho01}, \eqref{main04-02}, and \eqref{main04-06}, we can estimate that: 
\begin{align}\label{main04-07}
 \ds |\langle \zeta_i^m, \varphi \rangle_Y | &\, \leq \frac{1}{\tau}| p_i^{(*, m)} - p_{i+1}^{(*, m)}|_X|\varphi|_X + \nu^2 |\partial_x p_i^{(*, m)}|_X|\partial_x \varphi|_X \nonumber 
 \\
 &\, \quad + \int_\Omega |g'(w_i^{(*, m)}) p^{(*, m)}_i\varphi |\, dx +  M_w| w^{(*, m)}_i - w_i^{\mbox{\scriptsize ad}}|_X|\varphi|_X \nonumber
 \\
 &\, \leq \sup_{m \in \mathbb{N}}\Bigl( \frac{1}{\tau}| p_i^{(*, m)} - p_{i+1}^{(*, m)}|_X| + \nu^2 |\partial_x p_i^{(*, m)}|_X \nonumber
 \\
 &\, \qquad \qquad + R_2| p_i^{(*, m)}|_X + M_w| w^{(*, m)}_i - w_i^{\mbox{\scriptsize ad}}|_X \Bigr)|\varphi|_Y, \nonumber 
 \\
 &\, \qquad \quad \mbox{ for any } \varphi \in Y.
\end{align}
From \eqref{main04-01-01}, \eqref{main04-07}, and $(\hyperlink{star2}{\star\,2})$, we can see that: 
\begin{description}
    \item[\textmd{$(\hypertarget{star3}{\star\,3})$}]the sequence $ \{ \zeta^{m} \}_{m \in \mathbb{N}} = \{ [\zeta_1^{m}, \ldots , \zeta_n^m] \}_{m \in \mathbb{N}}$ is bounded in $ \mathbb{Y}^* $.
\end{description}
Having in mind \eqref{main04-01}, $(\hyperlink{star2}{\star\,2})$, and $(\hyperlink{star3}{\star\,3})$, we can find subsequences of $\{w^{(*, m)}\}_{m \in \mathbb{N}} $ $\subset \mathbb{X}$, $\{ p^{(*, m)}\}_{m \in \mathbb{N}} \subset \mathbb{X}$, and $ \{ \zeta^{m} \}_{m \in \mathbb{N}} = \{ [\zeta_i^m, \ldots \zeta_n^m ] \}_{m \in \mathbb{N}} \subset \mathbb{Y}^*$ (not relabeled), together with $w^\circ \in \mathbb{X}$, $p^\circ = [p_1^\circ, \ldots , p_n^\circ] \in \mathbb{X}$, and $\zeta^\circ = [\zeta_1^\circ, \ldots , \zeta_n^\circ] \in \mathbb{Y}^*$, such that: 
\begin{subequations}\label{main04-08}
\begin{align}\label{main04-08-04}
 p^{(*, m)} &\, \to p^{\circ} \mbox{ in } \mathbb{X}, \mbox{ weakly in } \mathbb{Y}, \mbox{ as } m \to \infty, 
\end{align}
and 
\begin{align}\label{main04-08-05}
 \zeta^m &\, \to \zeta^\circ \mbox{ weakly in } \mathbb{Y}^*, \mbox{ as } m \to \infty. 
\end{align}
\end{subequations}

Now, the properties \eqref{Thm.5-10}--\eqref{Thm.5-11} will be verified through the limiting observations for \eqref{main04-02}, as $m \to \infty$, with use of \eqref{main04-01} and \eqref{main04-08}.
\medskip

Next, let us fix $\rho > 0$ and define $\gamma_\rho$, as follows: 
        \begin{align}\label{Thm.5-12}
          \gamma_\rho(r) := & \begin{cases} 
                \gamma_0(r - \rho), & \mbox{ if } r \geq \rho,
                \\[0.5ex]
                \gamma_0(r + \rho), & \mbox{ if } r \leq -\rho,
                \\[0.5ex]
                0, & \mbox{ if } |r| < \rho.
            \end{cases}
        \end{align}
Also, we fix $i \in \{1, 2, 3, \ldots, n\}$, and we define 
\begin{align}\label{main04-09}
 M_\rho := \{ |\partial_x w_i^\circ| \geq \rho\}. 
\end{align}
In the light of Main Theorem \ref{mainTh01}, \eqref{main04-01}, and \eqref{main04-09}, there exists $m_\rho \in \mathbb{N}$ such that: 
\begin{align}\label{main04-10}
 \ds |\partial_x w_i^{(*, m)}| \geq \frac{\rho}{2} \mbox{ uniformly on } M_\rho, \mbox{ for all } m \geq m_\rho,
\end{align} 
and
\begin{align*}%\label{main04-13}
 w_i^\circ(M_\rho) \cup w_i^{(*, m)}(M_\rho) \subset (-1, 1) \mbox{ for all } m \geq m_\rho,
\end{align*}
i.e.
\begin{align}\label{main04-14}
 (K^{\delta_m})'(w_i^{(*, m)}) \equiv 0 \mbox{ on } M_\rho, \mbox{ for all } m \geq m_\rho.
\end{align}
By \eqref{A2.loc} as in the assumption (\hyperlink{A2l}{A2}), \eqref{main04-01}, and \eqref{main04-08}--\eqref{main04-10}, we can obtain that: 
\begin{align}\label{main04-11}
 \ds &\, \int_{M_{\rho}} (f^{\varepsilon_m})''(\partial_x w_i^{(*, m)})\partial_x p_i^{(*, m)} \cdot \nonumber
 \\
 &\, \qquad \qquad \qquad \cdot \bigl(\gamma_\rho'(\partial_x w_i^{(*, m)})\partial_x^2 w_i^{(*, m)}\varphi + \gamma_\rho(\partial_x w_i^{(*, m)})\partial_x \varphi \bigr)\, dx \nonumber 
 \\
 &\, \qquad \to 0 , \mbox{ for any } \varphi \in Y, \mbox{ as } m \to \infty.
\end{align} 
Moreover, from Main Theorem \ref{mainTh01}, 
\begin{align}\label{main04-11-1}
\partial_x w_i^{(*, m)}(\pm L) = 0, \mbox{ for all } m \geq m_\rho.
\end{align} 
Hence, as a consequence of $\gamma_\rho(0) = \gamma_\rho'(0) = 0$, \eqref{main04-01}, \eqref{main04-11}, and \eqref{main04-11-1}, we can compute that: 
\begin{align}\label{main04-12}
 \ds &\, \langle - \gamma_\rho(\partial_x w_i^{(*, m)}) \partial_x \bigl((f^{\varepsilon_m})''(\partial_x w_i^{(*, m)})\partial_x p_i^{(*, m)} + \nu^2 \partial_x p_i^{(*, m)} \bigr), \varphi \rangle_{Y} \nonumber
 \\
 &\, = \int_\Omega \bigl((f^{\varepsilon_m})''(\partial_x w_i^{(*, m)})\partial_x p_i^{(*, m)} + \nu^2 \partial_x p_i^{(*, m)} \bigr)\partial_x\bigl(\gamma_\rho(\partial_x w_i^{(*, m)})\varphi\bigr)\, dx \nonumber
 \\
 &\, = \int_{M_\rho} (f^{\varepsilon_m})''(\partial_x w_i^{(*, m)})\partial_x p_i^{(*, m)} \bigl( \gamma_\rho'(\partial_x w_i^{(*, m)})\partial_x^2 w_i^{(*, m)} \varphi + \gamma_\rho(\partial_x w_i^{(*, m)})\partial_x \varphi \bigr)\, dx \nonumber
 \\
 &\, \quad + \int_\Omega \nu^2 \partial_x p_i^{(*, m)} \bigl( \gamma_\rho'(\partial_x w_i^{(*, m)})\partial_x^2 w_i^{(*, m)} \varphi + \gamma_\rho(\partial_x w_i^{(*, m)})\partial_x \varphi \bigr)\, dx \nonumber
 \\
 &\, \to \int_\Omega \nu^2 \partial_x p_i^\circ \bigl( \gamma_\rho'(\partial_x w_i^\circ)\partial_x^2 w_i^\circ \varphi + \gamma_\rho(\partial_x w_i^\circ)\partial_x \varphi \bigr)\, dx, \mbox{ for any } \varphi \in Y, \mbox{ as } m \to \infty.
\end{align}
Meanwhile, on account of \eqref{Thm.5-12}, \eqref{main04-09}, and \eqref{main04-14}, we can easily check that: 
\begin{align}\label{main04-15}
 \gamma_\rho(\partial_x w_i^{(*, m)})(K^{\delta_m})'(w_i^{(*, m)}) \equiv 0, \mbox{ for all } m \geq m_\rho.
\end{align}

Here, we invoke that $p^{(*, m)}$ is the solution to the adjoint system \eqref{set03} for the forcing term $M_w (w^{(*, m)} - w^{\mathrm{ad}})$, i.e.
\begin{align}\label{main04-16}
 \displaystyle &\, \frac{1}{\tau}(p_i^{(*, m)} - p_{i+1}^{(*, m)}) - \partial_x ((f^{\varepsilon_m})'' (\partial_x w_i^{(*, m)})\partial_x p_i^{(*, m)} + \nu^2 \partial_x p_i^{(*, m)} ) \nonumber 
 \\
 &\, \quad + g'(w_i^{(*, m)})p_i^{(*, m)} + (K^{\delta_m})'(w_i^{(*, m)})p_i^{(*, m)} \nonumber
 \\
 &\,  =  M_w (w_i^{(*, m)} - w_i^{\mathrm{ad}}) \mbox{ in } \Omega, \mbox{ for any } i = n, \ldots, 3, 2, 1, \mbox{ and any } m \in \mathbb{N}.
\end{align} 
On this basis, we multiply $\gamma_\rho(\partial_x w_i^{(*, m)})$ and passing the limit of the both sides of \eqref{main04-16}, as $m \to \infty$. 
Then, from \eqref{main04-01}, \eqref{main04-08}, \eqref{main04-12}, and \eqref{main04-15}, one can see that: 
\begin{align}\label{main04-17}
 \ds &\, \frac{1}{\tau}\int_\Omega \gamma_\rho(\partial_x w_i^\circ)(p_i^\circ - p_{i+1}^\circ)\varphi\, dx + \int_\Omega \nu^2 \partial_x p_i^\circ \bigl( \gamma_\rho'(\partial_x w_i^\circ)\partial_x^2 w_i^\circ \varphi + \gamma_\rho(\partial_x w_i^\circ)\partial_x \varphi \bigr)\, dx, \nonumber
 \\
 &\, \quad + \int_\Omega \gamma_\rho(\partial_x w_i^\circ) g'(w_i^\circ)p_i^\circ \varphi\, dx = \int_\Omega M_w\gamma_\rho(\partial_x w_i^\circ)(w_i^\circ - w_i^\mathrm{ad})\varphi\, dx, \nonumber
 \\
 &\, \qquad \mbox{ for all } \varphi \in Y, \mbox{ and } i = n, \ldots , 3, 2, 1.
\end{align} 
Then, taking $\rho \to 0$ in \eqref{main04-17}, we can deduce that: 
\begin{align}\label{main04-17-1}
 \ds &\, \frac{1}{\tau}\int_\Omega \gamma_0(\partial_x w_i^\circ)(p_i^\circ - p_{i+1}^\circ)\varphi\, dx + \int_\Omega \nu^2 \partial_x p_i^\circ \bigl( \gamma_0'(\partial_x w_i^\circ)\partial_x^2 w_i^\circ \varphi + \gamma_0(\partial_x w_i^\circ)\partial_x \varphi \bigr)\, dx \nonumber
 \\
 &\, \quad + \int_\Omega \gamma_0(\partial_x w_i^\circ) g'(w_i^\circ)p_i^\circ \varphi\, dx = \int_\Omega M_w\gamma_0(\partial_x w_i^\circ)(w_i^\circ - w_i^\mathrm{ad})\varphi\, dx, \nonumber 
 \\
 &\, \qquad \mbox{ for all } \varphi \in Y, \mbox{ and } i = n, \ldots , 3, 2, 1.
\end{align}
Here, in the light of Main Theorem \ref{mainTh01} and $\gamma_0(0) = \gamma_0'(0) = 0$, we compute that: 
\begin{align}\label{main04-17-2}
\int_\Omega \nu^2 \partial_x p_i^\circ \bigl( \gamma_0'(\partial_x w_i^\circ)\partial_x^2 w_i^\circ \varphi &\, + \gamma_0(\partial_x w_i^\circ)\partial_x \varphi \bigr)\, dx \nonumber
\\
&\, = \int_\Omega \nu^2 \partial_x^2 p_i^\circ \gamma_0(\partial_x w_i^\circ)\varphi\, dx, \mbox{ for any } \varphi \in Y.
\end{align}
From \eqref{main04-17-1} and \eqref{main04-17-2}, we deduce that: 
\begin{align}\label{main04-18}
 \ds \left(\frac{1}{\tau}(p_i^\circ - p_{i+1}^\circ) - \nu^2\partial_x^2 p_i^\circ + g'(w_i^\circ) - M_w(w_i^\circ - w_i^\mathrm{ad}) \right)\gamma_0(\partial_x w_i^\circ) = 0 \mbox{ in } Y^*,
\end{align}
and furthermore,  
\begin{align}\label{main04-19}
 \ds &\, \langle \nu^2\partial_x^2 p_i^\circ \gamma_0(\partial_x w_i^\circ), \varphi \rangle_{Y} \nonumber
 \\
 &\, \qquad = \left(\left(\frac{1}{\tau}(p_i^\circ - p_{i+1}^\circ) + g'(w_i^\circ) - M_w(w_i^\circ - w_i^\mathrm{ad}) \right)\gamma_0(\partial_x w_i^\circ), \varphi \right)_X \nonumber 
 \\
 &\, \qquad \qquad  \mbox{ for any } \varphi \in Y.
\end{align}
\eqref{main04-18} and \eqref{main04-19} immediately leads to the required equation \eqref{Thm.5-13}.

Thus, we complete the proof of Main Theorem \ref{mainTh04}. \qed

\begin{rem}\label{main.del}
Especially, if we assume that: 
\begin{align}\label{Main04-00}
 \varepsilon_0 \in (0, 1] \mbox{ and } (f^\varepsilon)'' \to (f^{\varepsilon_0})'' \mbox{ in } C_{\mathrm{loc}}(\mathbb{R}), \mbox{ as } \varepsilon \to \varepsilon_0;
\end{align}
then we can obtain the characterization \eqref{Thm.5-13}, more precisely. 
In fact, applying \eqref{main04-01} under \eqref{Main04-00} yields that: 
\begin{align*}%\label{Main04-08-05}
 &\, (f^{\varepsilon_m})''(\partial_x w^{(*, m)}) \to (f^{\varepsilon_0})''(\partial_x w_i^{\circ}) \nonumber
 \\
 &\, \quad \mbox{ in } C(\overline{\Omega}), \mbox{ for any } i = n, \ldots, 3, 2, 1, \mbox{ as } m \to \infty,
\end{align*}
and this uniform convergence enables to improve, slightly, the characterization \eqref{Thm.5-13} as follows:  
\begin{align*}%\label{Main04-20}
 \ds \gamma_0(\partial_x w_i^\circ) &\, \Bigl( \frac{1}{\tau} (p_i^\circ - p_{i+1}^\circ) - \partial_x\bigl( (f^{\varepsilon_0})''(\partial_x w_i^{\circ})\partial_x p_i^{\circ} + \nu^2 \partial_x p_i^\circ \bigr)  \nonumber 
 \\
 &\, \quad + g'(w_i^\circ)p_i^\circ - M_w(w_i^\circ - w_i^{\mathrm{ad}}) \Bigr) = 0 \mbox{ in } X, \mbox{ for any } i = n, \ldots , 3, 2, 1. 
\end{align*}
\end{rem}

%%%%%%%%%%%%%%%%%%%%%%%%%%%%%%%%%%%%%%%%%%%%%%%%%%%%%%%%%%%%%%%%%%%%%%%%%%%%%%%%%%%%%%%%%%%%%%%%%%%%
\end{document}